\documentclass[11pt,amsfonts]{amsart}
\input fullpage.sty

\newtheorem{theorem}{Theorem}[section]

\newtheorem{lemma}[theorem]{Lemma}
\newtheorem{proposition}[theorem]{Proposition}
\newtheorem{definition}[theorem]{Definition}

\theoremstyle{remark}
\newtheorem{remark}[theorem]{\bf Remark}

\def\endproof{\qed \medskip}
\def\blacksquare{\hbox to .60em{\vrule width .60em height .60em}}

\begin{document}

\title[Geometric Structures on 3-Manifolds]{Scalar Curvature and the 
Existence of Geometric Structures on 3-Manifolds, I.}

\author[M.T. Anderson]{Michael T. Anderson}

\thanks{Partially supported by NSF Grants DMS 9505744 and 9802722}

\thanks{2000 {\it Mathematics Subject Classification}. Primary 58J60, 
57M50; Secondary 58E11, 53C20}

\maketitle

\contentsline {section}{\tocsection {}{0}{Introduction.}}{1}
\contentsline {section}{\tocsection {}{1}{Background Material.}}{7}
\contentsline {section}{\tocsection {}{2}{Geometrization of Tame 3-Manifolds.}}{12}
\contentsline {section}{\tocsection {}{3}{Metric Surgery on Spheres in Asymptotically Flat Ends.}}{30}
\contentsline {section}{\tocsection {}{4}{Asymptotically Flat Ends and Annuli.}}{36}
\contentsline {section}{\tocsection {}{}{References}}{43}

\abstract
This paper analyses the convergence and degeneration of sequences of 
metrics on a 3-manifold, and relations of such with Thurston's 
geometrization conjecture. The sequences are minimizing sequences for a 
certain (optimal) scalar curvature-type functional and their 
degeneration is related to the sphere and torus decompositions of the 
3-manifold under certain conditions.
\endabstract

\setcounter{section}{-1}

\section{Introduction.}
\setcounter{equation}{0}

 This paper and its sequel are concerned with the limiting behavior of 
minimizing sequences for certain curvature integrals on the space of 
metrics on a 3-manifold $M$, and the relations of such behavior with 
the geometrization conjecture of Thurston [37]. First, recall the 
statement of Thurston's conjecture, in the case of closed, oriented 
3-manifolds.

\medskip

\noindent
{\bf Geometrization Conjecture} (Thurston).

 {\it Let $M$ be a closed, oriented 3-manifold. Then $M$ admits a 
canonical decomposition into domains, each of which carries a canonical 
geometric structure.}

\medskip

 A geometric structure on a 3-manifold is a complete, locally 
homogeneous Riemannian metric. There are exactly eight such structures, 
namely the three constant curvature geometries, together with two 
further product geometries, $H^{2}\times {\mathbb R}$, $S^{2}\times {\mathbb 
R} $ and three twisted product geometries, $SL(2,{\mathbb R})$, Nil and 
Sol, c.f. [37], [33].

 The decomposition of $M$ is along certain essential 2-spheres and 
essential tori embedded in $M$. Recall that a 3-manifold $M$ is {\it  
irreducible}  if every embedded 2-sphere $S^{2}$ in $M$ bounds a 3-ball 
$B^{3}\subset M.$ 

 The sphere decomposition [22], [24] states that $M$ may be decomposed 
as a union of irreducible 3-manifolds, in that $M$ is diffeomorphic to 
a connected sum of closed, oriented 3-manifolds,
\begin{equation} \label{e0.1}
M = (M_{1}\# M_{2}\# \cdot   \cdot   \cdot  \# M_{p})\# (N_{1}\# 
N_{2}\# \cdot   \cdot   \cdot  \# N_{q})\# (\#_{1}^{r}(S^{2}\times 
S^{1})), 
\end{equation}
where each factor is irreducible, with the exception of $S^{2}\times 
S^{1}.$ The manifolds $M_{i}$ are defined to be those factors with 
infinite fundamental group, while $N_{j}$ are the factors with finite 
$\pi_{1}.$ Elementary 3-manifold topology implies that each manifold 
$M_{i}$ is a closed, 3-dimensional $K(\pi ,1),$ while each $N_{j}$ has 
universal cover given by a homotopy 3-sphere.

 Thus, the sphere decomposition allows one to understand the topology 
(and geometry) of 3-manifolds, in terms of the irreducible factors 
$M_{i}$ or $N_{j},$ since the topology and geometry of the standard 
factor $S^{2}\times S^{1}$ can be viewed as clear.

 In the present paper, we are only concerned with the structure of the 
$M_{i}$ factors, i.e. irreducible 3-manifolds with infinite fundamental 
group. The torus decomposition will be discussed later, c.f. Conjecture 
I and \S 2.

 The geometrization conjecture can be considered as a question 
concerning the existence of a ``best possible'' metric on a given 
3-manifold. This is a classical, albeit very difficult, question in 
Riemannian geometry. As mentioned above, we approach this existence 
question by seeking the minimum of a natural functional on the space of 
smooth metrics on $M$.

 Let ${\mathbb M} $ denote the space of all smooth Riemannian metrics on 
$M$ and consider the $L^{2}$ norm of the scalar curvature as a 
functional on ${\mathbb M} .$ This needs to be weighted with the 
appropriate power of the volume in order to obtain a scale-invariant 
functional. Thus, consider 
\begin{equation} \label{e0.2}
\mathcal{S}^{2}: {\mathbb M}  \rightarrow  {\mathbb R} ,\ \   \mathcal{S}^{2}(g) = 
\bigl(v^{1/3}\int_{M}s_{g}^{2}dV_{g}\bigr)^{1/2}, 
\end{equation}
where $s_{g}$ is the scalar curvature of the metric $g$, $dV_{g}$ is 
the volume form associated with $g$, and $v$ is the volume of $(M, g)$. 
As explained in detail in [5, \S 1], c.f. also Remark 3.2 below, it 
is technically advantageous to consider a slightly weaker functional 
than $\mathcal{S}^{2},$ namely
\begin{equation} \label{e0.3}
\mathcal{S}_{-}^{2}: {\mathbb M}  \rightarrow  {\mathbb R} , \ \  {\mathcal
S}_{-}^{2}(g) = \bigl(v^{1/3}\int_{M}(s_{g}^{-})^{2}dV_{g}\bigr)^{1/2}, 
\end{equation}
where $s^{-} =$ min $(s, 0)$ is the non-positive part of $s$. If $M$ 
carries a metric of positive scalar curvature, then it also carries 
scalar-flat metrics, so that inf $\mathcal{S}_{-}^{2} =$ inf $\mathcal{S}^{2} 
=$ 0. In this case, there is an infinite dimensional family of 
minimizers of $\mathcal{S}_{-}^{2}$ , so that there is no close relation 
between the geometry of such metrics and the topology of $M$. 

 On the other hand, let the Sigma constant $\sigma (M)$ be the supremum 
of the scalar curvatures of unit volume Yamabe metrics on $M$, c.f. 
[1], [30]. This is a topological invariant of the 3-manifold $M$, 
which should be thought of as an analogue of the Euler characteristic 
for surfaces. The set of closed 3-manifolds divides naturally into 
three topological classes, according to whether $\sigma (M)$ is 
negative, zero, or positive. 

 By the resolution to the Yamabe problem [29], $\sigma (M) \leq $ 0 if 
and only if $M$ carries no metrics of positive scalar curvature. It is 
not difficult to see, c.f. [5, Prop. 3.1], that if $\sigma (M) \leq $ 
0, then
\begin{equation} \label{e0.4}
\inf_{{\mathbb M} }\mathcal{S}_{-}^{2} = \inf_{{\mathbb M} }\mathcal{S}^{2}  = 
|\sigma (M)|. 
\end{equation}
Thus, for this paper we are exclusively interested in closed 
3-manifolds satisfying
\begin{equation} \label{e0.5}
\sigma (M) \leq  0. 
\end{equation}
We assume (0.5) holds throughout the paper. By a well-known result of 
Gromov-Lawson [15, Thm. 8.1], a closed 3-manifold $M$ satisfies 
$\sigma (M) \leq $ 0 if $M$ has at least one non-empty factor $M_{i}$ 
with infinite fundamental group in the sphere decomposition (0.1). Thus,
\begin{equation} \label{e0.6}
M_{i} \neq  \emptyset  , \ \  {\rm for \ some} \ \ i \ \ \Rightarrow  
\sigma (M) \leq  0. 
\end{equation}

 Further, any metric $g_{o}$ on $M$ realizing inf $\mathcal{S}  _{-}^{2},$ 
(or inf $\mathcal{S}^{2}),$ is necessarily an Einstein metric, and thus of 
constant sectional curvature, c.f. \S 2.

 Of course, an arbitrary closed, oriented 3-manifold does not admit an 
Einstein metric; this is the case for instance if $M$ has a non-trivial 
sphere decomposition (0.1). On the other hand, if $M$ is irreducible, 
we conjecture that there exist minimizing sequences $\{g_{i}\}$ for 
$\mathcal{S}_{-}^{2}$ which effectively implement the geometrization of 
$M$, in the sense that one can deduce the torus decomposition of $M$, 
and the geometrization of each canonical domain in $M$, from the 
limiting geometric behavior of $\{g_{i}\}.$ This is expressed in the 
following Conjectures, (equivalent to Conjectures I and II of [1], in 
the context of maximizing sequences of Yamabe metrics on $M$ in place 
of minimizing sequences of $\mathcal{S}_{-}^{2}$).

\medskip

\noindent
{\bf Conjecture I.}
  {\it Let $M$ be a closed, oriented, irreducible  3-manifold, with
$$\sigma (M) <  0. $$
Then there is a finite collection $\mathcal{T} $ of disjoint, embedded, 
incompressible tori $T_{i}^{2}\subset M,$ which separate $M$ into a 
union of two types of manifolds:
\begin{equation} \label{e0.7}
M\setminus \stackrel {\cdot} \cup T_{i}^{2} = \stackrel {\cdot} \cup 
H_{j} \ \cup \  \stackrel {\cdot} \cup G_{k}. 
\end{equation}
Each $H_{j}$ is a complete, connected hyperbolic manifold, of finite 
volume. The collection of boundary components of $\cup H_{j},$ i.e. the 
canonical tori in the hyperbolic cusps of $\{H_{j}\},$ forms exactly 
the collection $\mathcal{T} .$ Each $G_{k}$ is a connected graph manifold 
with toral boundary components, and the union of such boundary 
components again gives $\mathcal{T} .$ This decomposition of $M$ is unique 
up to isotopy of $M$.

 Let $vol_{-1}H_{j}$ denote the volume of $H_{j}$ in the hyperbolic 
metric. Then the Sigma constant of $M$ is given by 
\begin{equation} \label{e0.8}
|\sigma (M)| = (6\sum_{j} vol_{-1}H_{j})^{2/3}. 
\end{equation}
In particular, if $M$ is atoroidal, i.e. $M$ contains no ${\mathbb 
Z}\oplus{\mathbb Z}  \subset  \pi_{1}(M),$ then $M$ admits a hyperbolic 
structure which realizes the Sigma constant, i.e.}
\begin{equation} \label{e0.9}
|\sigma (M)| = (6vol_{-1}M)^{2/3}. 
\end{equation}

\medskip

 It is well-known that the condition that $M$ be irreducible and 
atoroidal is also necessary for $M$ to admit a hyperbolic metric.

\medskip

\noindent
{\bf Conjecture II.}
  {\it Let $M$ be a closed, oriented, irreducible 3-manifold, with
$$\sigma (M) = 0. $$
Then $M$ is a graph manifold satisfying 
$$|\pi_{1}(M)| = \infty . $$
The Sigma constant $\sigma (M)$ is realized by a smooth metric in 
${\mathbb M}_{1}$ if and only if $M$ is a flat 3-manifold.}

\medskip

 A {\it  graph manifold}  is a union of $S^{1}$ bundles over surfaces, 
or equivalently a union of Seifert fibered spaces, glued together by 
toral automorphisms along toral boundary components, c.f. [38,39]. This 
is also exactly the class of 3-manifolds which admit an F-structure, in 
fact a polarized F-structure, in the sense of Cheeger-Gromov [7,8], 
c.f. also [13, App.2]. 

 The geometrization of graph manifolds is relatively straightforward 
and well understood. Briefly, any irreducible graph manifold admits a 
further decomposition into domains, with toral boundary components, 
each of which admits a geometric structure modelled on one of the 
Seifert fibered geometries or it admits a Sol geometry. This 
decomposition and geometrization of the graph manifolds is obtained 
naturally from the proofs of the results below, c.f. \S 2.3 for further 
discussion.

 Topologically, these two conjectures imply that a closed, oriented, 
irreducible 3-manifold with $\sigma (M) \leq $ 0 is a union of 
hyperbolic manifolds and graph manifolds glued together along 
incompressible toral boundary components. The union of the graph 
manifolds in $M$ is also called the {\it  characteristic subvariety}  
of $M$. The case $\sigma (M) =$ 0 implies the absence of hyperbolic 
components, while the case $\sigma (M) < $ 0 implies their existence.

 Geometrically, Conjectures I and II are meant to describe the limiting 
behavior of a suitable minimizing sequence $\{g_{i}\}$ for 
$\mathcal {S}_{-}^{2}.$ Thus, in I, conjecturally the sequence $\{g_{i}\}$ 
converges on the manifolds $H_{j}$ in (0.7) to a complete, finite 
volume, constant curvature metric, (with scalar curvature $\sigma 
(M)),$ and collapses the graph manifold components $G_{k}$ in (0.7) 
along $S^{1}$ or $T^{2}$ fibers. In II, conjecturally the sequence 
$\{g_{i}\}$ fully collapses $M$ along $S^{1}$ or $T^{2}$ fibers, so 
that there is no limiting metric in general. In these cases of 
collapse, although the metrics are not converging, their degeneration 
implies the existence of a well defined topological structure, (a graph 
manifold structure), incompressibly embedded in $M$, which essentially 
describes how the degeneration is occuring. In effect, the limiting 
behavior of the sequence $\{g_{i}\}$ implements, or performs, the 
geometrization of $M$.

 Conjectures I and II, together with (0.6), are easily seen to imply 
the geometrization conjecture for closed, oriented irreducible 
3-manifolds with infinite fundamental group, c.f. [1, \S 4], or also 
Remark 2.11 and \S 2.3.

\medskip

 In this paper, we prove Conjectures I and II in case either one of two 
additional topological assumptions hold on $M$. Some further background 
is needed to explain these assumptions. For a given $\varepsilon  > $ 0 
small, consider the scale-invariant perturbation of $\mathcal{S}  
_{-}^{2}$ given by
\begin{equation} \label{e0.10}
I_{\varepsilon}^{~-} = \varepsilon v^{1/3}\int|z|^{2}dV + 
\bigl(v^{1/3}\int (s^{-})^{2}dV\bigr)^{1/2}. 
\end{equation}
The functional $I_{\varepsilon}^{~-}$ is a perturbation of 
$\mathcal S_{-}^{2}$ in the direction of the $L^{2}$ norm of the trace-free 
Ricci curvature $z$. While the existence of minimizers of 
$\mathcal S_{-}^{2}$ is difficult to prove - this is the essential geometric 
content of Conjectures I and II - it is not so difficult to prove the 
existence of minimizers $g_{\varepsilon}$ of $I_{\varepsilon}^{~-}.$ 
This is done in [2,5], where the following results are proved, (c.f. 
also Theorem 1.1 below for more details). For any closed oriented 
3-manifold $M$, (not necessarily irreducible), there exists a maximal 
domain $\Omega_{\varepsilon}$ and a $C^{2,\alpha}$ complete Riemannian 
metric $g_{\varepsilon}$ on $\Omega_{\varepsilon},$ for which the pair 
$(\Omega_{\varepsilon}, g_{\varepsilon})$ realizes $\inf_{{\mathbb M}_{1}} 
I_{\varepsilon}^{~-},$ in that
\begin{equation} \label{e0.11}
 I_{\varepsilon}^{~-}(g_{\varepsilon}) = 
\varepsilon\int_{\Omega_{\varepsilon}}|z_{g_{\varepsilon}}|^{2} + 
\bigr(\int_{\Omega_{\varepsilon}}(s_{g_{\varepsilon}}^{~-})^{2}\bigl)^{1
/2} = \inf_{{\mathbb M}_{1}} I_{\varepsilon}^{~-}, 
\end{equation}
and
\begin{equation} \label{e0.12}
vol_{g_{\varepsilon}}\Omega_{\varepsilon} = 1. 
\end{equation}
The domain $\Omega_{\varepsilon}$ weakly embeds in $M$, in the sense 
that any compact domain with smooth boundary embeds as such a domain in 
$M$. There is an exhaustion of $\Omega_{\varepsilon}$ by compact 
domains $K_{j}$ with each $\partial K_{j}$ given by a finite collection 
of smooth tori,  such that the complement $M\setminus K_{j}$ is a graph 
manifold in $M$ and so admits an F-structure $\mathcal{F}_{j}.$

 The domain $\Omega_{\varepsilon}$ is empty exactly when the closed 
3-manifold $M$ itself is a graph manifold. If $M$ is not a graph 
manifold, then $\Omega_{\varepsilon} \neq  \emptyset  ,$ and one may 
have $\Omega_{\varepsilon} = M$, or $\Omega_{\varepsilon}$ only weakly 
embedded in $M$.

\medskip

 It is proved in [5, (3.3)-(3.4)], that if $g_{\varepsilon}$ is a 
minimizer for $I_{\varepsilon}^{~-}$ as above, then as $\varepsilon  
\rightarrow $ 0,
\begin{equation} \label{e0.13}
\mathcal{S}_{-}^{2}(g_{\varepsilon}) \rightarrow  |\sigma (M)|, 
\end{equation}
and
\begin{equation} \label{e0.14}
\varepsilon\int|z_{g_{\varepsilon}}|^{2}dV_{g_{\varepsilon}} 
\rightarrow  0. 
\end{equation}

 Thus the family $(\Omega_{\varepsilon}, g_{\varepsilon})$ as 
$\varepsilon  \rightarrow $ 0, may be considered as a specific 
minimizing family for the functional $\mathcal{S}_{-}^{2}.$ This family is 
optimal in that it has the least amount of curvature in $L^{2}$ in the 
following sense: given any $(\Omega_{\varepsilon}, g_{\varepsilon})$ as 
above, for any smooth unit volume metric $\bar g$ on $M$ for which
$$\mathcal{S}_{-}^{2}(\bar g) \leq  \mathcal{S}_{-}^{2}(g_{\varepsilon}), $$
one has
$$\int_{M}|z_{\bar g}|^{2}dV_{\bar g} \geq  
\int_{\Omega_{\varepsilon}}|z_{g_{\varepsilon}}|^{2}dV_{g_{\varepsilon}}
. $$

 A pair $(\Omega_{\varepsilon}, g_{\varepsilon})$ as above is called a 
minimizing pair for $I_{\varepsilon}^{~-}.$ We point out that it is not 
known, (without resolution of Conjectures I and II), if a minimizing 
pair is unique. Thus neither the metrics $g_{\varepsilon},$ nor the 
topological type of the domains $\Omega_{\varepsilon},$ are known to be 
unique, in the sense that they depend only on the topology of $M$. The 
only exception to this is when $M$ is a graph manifold, when, as 
mentioned above, both $\Omega_{\varepsilon},$ and so $g_{\varepsilon},$ 
are empty. On the other hand, the collection of all minimizing pairs 
$(\Omega_{\varepsilon}, g_{\varepsilon}), \varepsilon  > $ 0, depends 
only on the topology of $M$.

 With this background, we now are able to state the first topological 
assumption on $M$.
\begin{definition} \label{d 0.1.}
 Let $M$ be a closed, oriented 3-manifold. Then $M$ is {\sf tame}  if 
there exists a sequence $\varepsilon  = \varepsilon_{i} \rightarrow $ 
0, and a sequence of minimizing pairs $(\Omega_{\varepsilon}, 
g_{\varepsilon})$ for $I_{\varepsilon}^{~-}$ on $M$ such that
\begin{equation} \label{e0.15}
\int_{\Omega_{\varepsilon}}|z_{\varepsilon}|^{2}dV_{g_{\varepsilon}} 
\leq  \Lambda , 
\end{equation}
for some $\Lambda  <  \infty .$ 
\end{definition}

 It is clear that the condition that $M$ is tame is a topological 
condition, i.e. depends only on the smooth, and hence topological, 
structure of $M$. In fact, from the definition of $g_{\varepsilon}$, 
$M$ is tame if and only if there exists a sequence $\{g_{i}\}$ of 
smooth metrics on $M$, which form a minimizing sequence for 
$\mathcal S_{-}^{2},$ and which have uniformly bounded $z$-curvature 
in $L^{2}.$

\begin{theorem} \label{t 0.2.}
  Let $M$ be a closed oriented 3-manifold with $\sigma (M) \leq $ 0. If 
$M$ is tame, then Conjectures I and II hold for $M$.
\end{theorem}

 We mention that Theorem 0.2 remains valid even if the tame manifold 
$M$ is not irreducible. Of course the hyperbolic part $\cup H_{j}$ of 
$M$ in (0.7) is irreducible, but the graph manifold part $\cup G_{k}$ 
of $M$, and hence $M$ itself, could be reducible. In particular, a tame 
3-manifold need not be irreducible, c.f. \S 1. 

  A proof of the analogue of Theorem 0.2 for maximizing sequences of 
Yamabe metrics was outlined in [1], without full details however. The 
proof given here may be adapted to such sequences of Yamabe metrics 
without difficulty. Theorem 0.2 is also analogous to a recent result of 
Hamilton [16] on the long-time behavior of the Ricci flow assuming a 
uniform $L^{\infty}$ bound on the curvature.

\medskip

 Next we turn to the more interesting and difficult case where $M$ is 
not tame, so that an arbitrary minimizing sequence for 
$\mathcal{S}_{-}^{2}$ has curvature diverging to infinity in $L^{2}.$ In 
particular, this is the case for any sequence of minimizing pairs 
$(\Omega_{\varepsilon}, g_{\varepsilon}),$ with $\varepsilon  = 
\varepsilon_{i}$ any sequence converging to 0. The metrics 
$g_{\varepsilon}$ are thus degenerating and one would like to relate 
this degeneration to the topology of $M$.

 In this case, it is proved in [5, Thm. B], c.f. also Theorem 1.3 
below, that for any sequence $\varepsilon_{i} \rightarrow $ 0 and 
sequence of minimizing pairs $(\Omega_{\varepsilon}, g_{\varepsilon}), 
\varepsilon  = \varepsilon_{i},$ there exist base points 
$y_{\varepsilon}\in\Omega_{\varepsilon}$ and scale factors $\rho 
(y_{\varepsilon}) \rightarrow $ 0 as $\varepsilon  \rightarrow $ 0, 
such that the rescaled or blow-up metrics
\begin{equation} \label{e0.16}
g_{\varepsilon}'  = \rho (y_{\varepsilon})^{-2}\cdot  g_{\varepsilon} 
\end{equation}
based at $\{y_{\varepsilon}\}$ have a subsequence converging to a 
complete non-flat Riemannian manifold $(N, g' , y)$ of non-negative 
scalar curvature and uniformly bounded curvature. The limit $(N, g' )$ 
minimizes the $L^{2}$ norm of $z$ over all comparison metrics $\bar g$ 
of non-negative scalar curvature on $N$ such that $vol_{\bar g} K \leq  
vol_{g}K$ and $\bar g|_{N\setminus K} = g|_{N\setminus K},$ for some 
arbitrary compact set $K \subset  N$. The base points $y_{\varepsilon} 
\in \Omega_{\varepsilon}$ are preferred in that the curvature of 
$g_{\varepsilon}$ is (locally) maximal near $y_{\varepsilon}$.

 The structure of these complete metrics $(N, g, y)$ models the small 
scale geometry of the degeneration of $(\Omega_{\varepsilon}, 
g_{\varepsilon}),$ near the base points $y_{\varepsilon}$ where the 
curvature blows up. The structure of the limit $(N, g')$ is discussed 
in more detail in \S 1. We point out that the base points 
$y_{\varepsilon},$ as well as the choice of subsequence above may not 
be unique. Thus it is possible apriori that there are many distinct 
blow-up limits $(N, g').$ 

\begin{definition} \label {d 0.3.}
 Let $M$ be a closed oriented 3-manifold. Then $M$ is {\sf spherically 
tame}  if there exists a complete non-flat blow-up limit $(N, g' , y)$ 
as above which has an {\sf  asymptotically flat end}  $E \subset  N$ in 
the following sense: there is a diffeomorphism $F$: ${\mathbb 
R}^{3}\setminus B \rightarrow  E$, for some compact ball $B\subset 
{\mathbb R}^{3},$ such that the metric $g'$ has the form
\begin{equation} \label{e0.17}
g_{ij}' = (1+\frac{2m}{r})\delta_{ij} + h, 
\end{equation}
in the chart $F$, where $h=O(r^{-2}), |D^{p}h| = O(r^{-2+p}), p =$ 1,2, 
$r(x) = |x|.$ The parameter $m$ is the mass of the end $E$ and is 
assumed to be positive.
\end{definition}

 Again this condition depends only on the smooth structure of the 
manifold $M$. Thus the class of non-tame 3-manifolds, (closed and 
oriented), divides into two classes, namely those which are spherically 
tame, and those which are not. Of course, this condition is expressed 
in terms of analysis and geometry, and so its relation with standard 
topological concepts may not be immediately clear. In this respect, we 
prove the following result in \S 3.

\begin{theorem} \label{t 0.4.}
  Let $M$ be a closed oriented 3-manifold with $\sigma (M) \leq $ 0, 
which is not tame but is spherically tame. Then $M$ is reducible, i.e. 
$M$ contains an essential embedded 2-sphere.
\end{theorem}

 In fact we prove that the ``natural'' 2-sphere $S^{2}$ in the 
asymptotically flat end $E \subset  N$ is an essential 2-sphere in $M$. 
The method of proof is by a cut and paste or comparison argument. If 
this $S^{2}$ bounds a 3-ball $B^{3}$ in $M$, we construct a specific 
metric glueing of such a $B^{3}$ onto $S^{2}$ so as to decrease the 
value of $I_{\varepsilon}^{~-}$ a definite amount below that of 
$I_{\varepsilon}^{~-}(g_{\varepsilon})$, for $\varepsilon$ sufficiently 
small, which contradicts the minimizing property of $g_{\varepsilon}$. 
An analogous cut and paste argument is used in the proof of Theorem 0.2 
to prove the tori $\mathcal{T}$ in (0.7) are incompressible.

  The fact that it is possible to carry out such cut and paste 
arguments on spheres and tori is a strong reason for preferring the 
functional $I_{\varepsilon}^{~-}$ over other candidates.

\medskip

  Theorem 0.4 relates the geometry and analysis of the metrics 
$g_{\varepsilon}$ with the sphere decomposition (0.1). In particular, 
it implies that if $M$ is irreducible and spherically tame, then $M$ 
must be tame. Hence, via Theorems 0.2 and 0.4, Conjectures I and II are 
reduced to the question of whether a closed oriented 3-manifold is 
automatically spherically tame, i.e. to the following:

\medskip

\noindent
{\bf Sphere conjecture.} 
{\it Any closed oriented 3-manifold $M$ which is not tame is spherically 
tame}.

\medskip

 The Sphere conjecture will be the focus of the sequel paper. This 
conjecture asserts that there exists a blow-up limit $(N, g' , y)$ of a 
non-tame minimizing sequence $(\Omega_{\varepsilon}, g_{\varepsilon}), 
\varepsilon  = \varepsilon_{i} \rightarrow $ 0, which has an 
asymptotically flat end. More loosely speaking, it claims that one can 
detect 2-spheres $S^{2}$ in $M$, (embedded in almost flat regions), 
from the degenerating geometry of a suitable minimizing sequence for 
$\mathcal{S}_{-}^{2}$.

  The existence of such an end might seem to be a strong condition. In 
fact, asymptotically flat ends are the most natural end structure for 
complete blow-up limits, and no situations are known where such a 
metric might have a non-asymptotically flat end. In this regard, 
[5, Thm.C], c.f. also Theorem 1.4 below, gives a relatively simple 
characterization of blow-up limits $(N, g' )$ which have a finite 
number of ends, each of which is asymptotically flat. To prepare for 
work to follow in the sequel, in \S 4 we extend this result by 
characterizing those blow-up limits $(N, g' )$ which have at least one 
asymptotically flat end, c.f. Proposition 4.1. Further, in Theorem 4.2 
the arguments proving Theorem 0.4 are shown to extend to limits $(N, g' 
)$ which have a suitable almost flat annulus in place of a full 
asymptotically flat end.

 Briefly, the paper is organized as follows; more detailed remarks on 
the contents are given at the beginning of each section or subsection. 
Following discussion of some necessary background material in \S 1, the 
proof of Theorem 0.2 is given in \S 2. Also, in \S 2.3, we discuss the 
geometrization of graph manifolds, as obtained from the limiting 
behavior of minimizing sequences for $\mathcal{S}_{-}^{2}$, c.f. Theorem 
2.14. The proof of Theorem 0.3 is given in \S 3, while \S 4 collects 
some results and techniques on the structure of blow-up limits $(N, g' 
)$ needed for the sequel paper.

  While this paper relies to a certain extent on results from earlier 
papers [2]-[5], we have made the paper at least logically almost 
self-contained in that the results needed from these papers are 
summarized in \S 1.

\section{Background Material.}
\setcounter{equation}{0}

\addtocounter{theorem}{-1}

 Throughout the paper, we assume knowledge of the Cheeger-Gromov theory 
of convergence and collapse of Riemannian manifolds, [12], [7,8], as 
well as the $L^{2}$ Cheeger-Gromov theory on 3-manifolds, [2, \S 3], 
c.f. also [5, \S 2] for a summary.

 The following geometric quantities will be used frequently.

\begin{definition} \label{d 1.1.}
 Let (N, g) be a complete Riemannian manifold, and $x\in N.$

{\bf (i).}
 The $L^{2}$ curvature radius at $x$ is the largest radius $\rho (x)$ 
such that for any geodesic ball $B_{y}(s) \subset B_{x}(\rho (x)), s 
\leq  \rho (x),$ one has
\begin{equation} \label{e1.1}
\frac{s^{4}}{volB_{y}(s)}\int_{B_{y}(s)}|r|^{2}dV \leq  c_{o}; 
\end{equation}
where $r$ is the Ricci curvature and $c_{o}$ is a fixed small positive 
parameter. Throughout the paper, we set $c_{o} = 10^{-3}$ for 
convenience.

{\bf (ii).}
 The volume radius $\nu (x)$ at $x$ is given by
\begin{equation} \label{e1.2}
\nu (x) = \sup\{r:\frac{vol(B_{y}(s))}{s^{3}} \geq  \mu , \forall 
B_{y}(s) \subset B_{x}(r), s \leq  r\}, 
\end{equation}
where again $\mu $ is a fixed small parameter, which measures the 
degree of the volume collapse near x. As above, to be concrete, we 
assume $\mu  = 10^{-1}$ throughout the paper.
\end{definition}

 There are similar definitions for the $L^{k,p}$ curvature radius,  $k 
\geq $ 0, 1 $<  p <  \infty ,$ where the $L^{2}$ norm of $r$ in (1.1) 
is replaced by its $L^{k,p}$ norm, and the power $s^{4}$ is replaced by 
the power making the expression analogous to (1.1) scale-invariant. It 
is important to note that the $L^{k,p}$ curvature radius is continuous 
under convergence in the strong $L^{k,p}$ topology, c.f. [2, \S 3] 
and further references therein.

\medskip

 A domain $\Omega $ is {\it  weakly embedded}  in a 
closed 3-manifold $M$ if every subdomain $K \subset  \Omega $ with 
compact closure and with smooth boundary in $\Omega $ may be smoothly 
embedded as such a domain in $M$; in this case, we write
\begin{equation} \label{e1.3}
\Omega  \subset\subset  M .
\end{equation}
A domain is defined to be an open 3-manifold, not necessarily connected.

\medskip

  The Euler-Lagrange equations $\nabla I_{\varepsilon}^{~-}(g) =$ 0 for 
the functional $I_{\varepsilon}^{~-}$ from (0.11) at a unit volume 
metric $g$ are the following system of elliptic PDE, c.f. [5, \S 3]:
\begin{equation} \label{e1.4}
\varepsilon\nabla\mathcal{Z}^{2} + L^{*}\tau  + \phi\cdot  g = 0, 
\end{equation}
\begin{equation} \label{e1.5}
2\Delta (\tau +\frac{\varepsilon s}{12}) + \frac{1}{4}s\tau  = 
-\frac{1}{2}\varepsilon|z|^{2} + 3c. 
\end{equation}
Here $L^{*}$ is the $L^{2}$ adjoint of the linearization of the scalar 
curvature, given by
\begin{equation} \label{e1.6}
L^{*}u = D^{2}u -  \Delta u\cdot  g -  u\cdot  r,
\end{equation}
where $D^{2}$ is the Hessian and $\Delta  = trD^{2}$ is the Laplacian 
w.r.t. the metric. The term $\nabla\mathcal{Z}^{2}$ in (1.4) is the 
gradient of $\mathcal{Z}^{2} = \int|z|^{2}dV,$ and is of the form
$$\nabla\mathcal{Z}^{2} = D^{*}Dz + \frac{1}{3}D^{2}s -  
2\stackrel{\circ}{R} \circ z + \frac{1}{2}(|z|^{2} -  \frac{1}{3}\Delta 
s)\cdot  g,$$ 
c.f. [6,\S 4H]. The lower order terms in (1.4)-(1.5) are given by $\phi 
 = -\frac{1}{4}s\tau  + c$, $c = \frac{1}{12\sigma}\int 
(s^{-})^{2}+\frac{\varepsilon}{6}\int|z|^{2},$ $s$ is the scalar 
curvature of $g$ and
\begin{equation} \label{e1.7}
\tau  = \frac{s^{-}}{\sigma}, 
\end{equation}
with $\sigma  = (\int (s^{-})^{2})^{1/2}.$ Note that $\tau$ is 
non-positive and has $L^{2}$ norm equal to 1.

 The following result from [5,Thm.3.8] proves the existence and basic 
geometric properties of minimizers of the functional 
$I_{\varepsilon}^{~-}$.

\begin{theorem} \label{t 1.1.}
  {\bf (Geometric Decomposition for $I_{\varepsilon}^{~-}$).}
 Suppose $\sigma (M) \leq $ 0. For any $\varepsilon  > $ 0, there is a 
complete, $L^{3,p}\cap C^{2,\alpha}$ Riemannian metric 
$g_{\varepsilon},$ defined on a domain $\Omega_{\varepsilon} 
\subset\subset $ M, which realizes $\inf_{{\mathbb M}_{1}} 
I_{\varepsilon}^{~-},$ in the sense that
\begin{equation} \label{e1.8}
 I_{\varepsilon}^{~-}(g_{\varepsilon}) = 
\varepsilon\int_{\Omega_{\varepsilon}}|z_{g_{\varepsilon}}|^{2} + 
\bigr(\int_{\Omega_{\varepsilon}}(s_{g_{\varepsilon}}^{-})^{2}\bigl)^{1/
2} = \inf_{{\mathbb M}_{1}} I_{\varepsilon}^{~-}, 
\end{equation}
and
\begin{equation} \label{e1.9}
vol_{g_{\varepsilon}}\Omega_{\varepsilon} = 1. 
\end{equation}
The metric $g_{\varepsilon}$ weakly satisfies the Euler-Lagrange 
equations (1.4)-(1.5) and is $C^{\infty} smooth,$ in fact 
real-analytic, away from the locus where $s =$ 0.

 Further, the curvature $R$ of $g_{\varepsilon}$ is uniformly bounded 
on $\Omega_{\varepsilon}$ and $\Omega_{\varepsilon}$ consists of a 
finite number $Q = Q(\varepsilon ,$ M) of components. There is an 
exhaustion of $\Omega_{\varepsilon}$ by compact domains $K_{j}$ with 
$\partial K_{j}$ given by a finite collection of smooth tori, such that 
the complement $M\setminus K_{j}$ is a graph manifold embedded in $M$.
\end{theorem}

 As in the Introduction, $(\Omega_{\varepsilon}, g_{\varepsilon})$ is 
called a minimizing pair for $I_{\varepsilon}^{~-}.$ The domain 
$\Omega_{\varepsilon}$ is empty if and only if $M$ itself is a graph 
manifold. In fact,
\begin{equation} \label{e1.10}
\inf I_{\varepsilon}^{~-} = 0, \ \ {\rm for \ all} \ \ \varepsilon  >  
0, 
\end{equation}
if and only if $M$ is a graph manifold, c.f. [13, App.2] or [5, \S 
3.1]. Of course, this implies inf $\mathcal{S}  _{-}^{2} =$ 0 also for 
graph manifolds. (The converse is the main content of Conjecture II). 
Graph manifolds can be characterized geometrically as exactly the class 
of 3-manifolds which admit a volume collapse with bounded curvature, 
i.e. there exists a sequence of metrics $\{g_{k}\}$ on $M$ such that 
$|R_{g_{k}}| \leq $ 1 everywhere, while $vol_{g_{k}}M \rightarrow $ 0 
as $k \rightarrow  \infty .$ Equivalently, they admit an F-structure, 
in fact a polarized F-structure, in the sense of Cheeger-Gromov 
[7,8]. Briefly, a 3-manifold admits an F-structure if there is a 
partition of $M$ into domains, each of which admits either an $S^{1}$ 
or a $T^{2}$ action, for which the dimension of the orbits is positive 
everywhere. The F-structure is polarized if the group actions are 
locally free.

 It is interesting and important to note that the class of graph 
manifolds is closed under connected sums, c.f. [39, p.91], or also [34, 
Lemma 4], so that graph manifolds may well be reducible. Via (1.10), 
this shows in particular that tame 3-manifolds need not be irreducible, 
c.f. also \S 2.3 for further discussion.

\medskip

 Theorem 0.1 gives a geometric decomposition of $M$ w.r.t. the 
functional $I_{\varepsilon}^{~-}$, in that $M$ is the union of 
$\Omega_{\varepsilon}$, or more precisely $K_{j} \subset 
\Omega_{\varepsilon}$, and its graph manifold complement in $M$, with 
$(\Omega_{\varepsilon}, g_{\varepsilon})$ a solution to a natural 
geometric variational problem. However, it is very unlikely that this 
geometric decomposition corresponds to the geometric decomposition 
given in Conjecture I on arbitrary, in particular on reducible 
3-manifolds. In situations where $M$ is irreducible and Conjectures I 
and II are known to hold, then these two decompositions do coincide, 
c.f. Remark 2.10.

 On the other hand, without such knowledge, it is unknown if 
$\Omega_{\varepsilon}$ even has finite topological type, or if the 
number of components of $\Omega_{\varepsilon}$ remains uniformly 
bounded on a sequence $\varepsilon  = \varepsilon_{i} \rightarrow $ 0. 
As stated in the Introduction, it is also unknown if 
$(\Omega_{\varepsilon}, g_{\varepsilon})$ is unique, for $\varepsilon > 
0$ fixed.

  The behavior of the potential function $\tau = \tau_{\varepsilon}$ 
from (1.4)-(1.5) plays an important role throughout the paper. Let $T = 
T_{\varepsilon} = \sup_{\Omega_{\varepsilon}}|\tau_{\varepsilon}| = - 
\inf_{\Omega_{\varepsilon}}\tau_{\varepsilon}$. Then $T \geq 1$, since 
the $L^{2}$ norm of $\tau_{\varepsilon}$ over $\Omega_{\varepsilon}$ 
equals 1. By [5,(3.38)], one has the bound
\begin{equation} \label{e1.11}
1 \leq T \leq (1+\frac{2\varepsilon}{\sigma}\mathcal{Z}^{2}(g_{\varepsilon}))^{1/2}.
\end{equation}
The following result from [5, Thms. 3.10-3.11] summarizes some of the 
properties of $\tau_{\varepsilon}$.
\begin{proposition} \label{p 1.2.}
If $\sigma(M) < 0$, then the function $\tau_{\varepsilon}$ satisfies
\begin{equation} \label{e1.12}
\inf_{\Omega_{\varepsilon}}\tau_{\varepsilon} \rightarrow -1,
\end{equation}
and, for any $p < \infty$,
\begin{equation} \label{e1.13}
\int_{\Omega_{\varepsilon}}|\tau_{\varepsilon}+ 
1|^{p}dV_{g_{\varepsilon}} \rightarrow 0,
\end{equation}
as $\varepsilon \rightarrow 0$.
If $\sigma(M) = 0$, then
\begin{equation} \label{e1.14}
\inf_{\Omega_{\varepsilon}}s_{\varepsilon} \rightarrow 0, \ {\rm as} \ 
\varepsilon \rightarrow 0.
\end{equation}
Further, in both cases,
\begin{equation} \label{e1.15}
\int_{\Omega_{\varepsilon}}|\nabla 
\tau_{\varepsilon}|^{2}dV_{g_{\varepsilon}} \rightarrow 0.
\end{equation}
\end{proposition}

 We next summarize below two of the main results, namely Theorems B,C, 
of [5] concerning the structure of the blow-up limits $(N, g)$ 
discussed in \S 0. In analogy to the Euler-Lagrange equations 
(1.4)-(1.5), the {\sf $\mathcal{Z}_{c}^{2}$ equations} are an elliptic 
system of equations in a metric $g$, given by
\begin{equation} \label{e1.16}
\alpha\nabla\mathcal{Z}^{2} + L^{*}\tau  = 0, 
\end{equation}
\begin{equation} \label{e1.17}
\Delta (\tau +\frac{\alpha}{12}s) = -\frac{\alpha}{4}|z|^{2}, 
\end{equation}
for some constant $\alpha  > $ 0. As previously in (1.4)-(1.5), the 
function $\tau$ is viewed as a potential function.

\begin{theorem} \label{t 1.3.} {\bf (Structure of Blow-up Limits).} 
 Suppose $\sigma (M) \leq $ 0, and let $(\Omega_{\varepsilon}, 
g_{\varepsilon}), \varepsilon  = \varepsilon_{i} \rightarrow $ 0, be a 
sequence of minimizing pairs for $I_{\varepsilon}^{~-}.$ Suppose the 
sequence $\{(\Omega_{\varepsilon}, g_{\varepsilon})\}$ degenerates, in 
the sense that 
\begin{equation} \label{e1.18}
\int_{\Omega_{\varepsilon}}|z_{g_{\varepsilon}}|^{2}dV_{g_{\varepsilon}}
 \rightarrow  \infty , \ \ {\rm as} \ \ \varepsilon  \rightarrow  0. 
\end{equation}
Then there exist points $\{y_{\varepsilon}\}\in (\Omega_{\varepsilon}, 
g_{\varepsilon}),$ with $\rho (y_{\varepsilon}) \rightarrow $ 0, such 
that the blow-up metrics
\begin{equation} \label{e1.19}
g_{\varepsilon}'  = \rho (y_{\varepsilon})^{-2}\cdot  g_{\varepsilon}, 
\end{equation}
based at $y_{\varepsilon},$ have a subsequence converging in the strong 
$L^{2,2}$ topology to a limit $(N, g' , y)$. The limit $(N, g' )$ is a 
complete, non-flat Riemannian manifold, with uniformly bounded 
curvature, and non-negative scalar curvature s. 

  Further, $(N, g' )$ minimizes the $L^{2}$ norm of the curvature $z$ 
over all metrics $\bar g$ of non-negative scalar curvature satisfying 
$vol_{\bar g}K \leq  vol_{g'}K$ and $\bar g|_{N\setminus K} = 
g'|_{N\setminus K},$ for some compact set $K \subset  N$. The metric 
$g' $ is $C^{2,\beta} \cap L^{3,p}$ smooth, for any $\beta  <$ 1 and $p 
<  \infty ,$ and is an $L^{3,p}$ weak solution of the 
$\mathcal{Z}_{c}^{2}$ equations (1.16)-(1.17).

 The potential function $\tau $ is non-positive, and $\tau $ and $s$ 
are locally Lipschitz functions on $N$ with disjoint supports, in the 
sense that
$$s\cdot \tau  \equiv  0. $$
The metric $g' $ is $C^{\infty}$ smooth, in fact real-analytic, and the 
convergence to the limit is $C^{\infty}$ smooth, uniformly on compact 
subsets, in any region where $\tau  < $ 0 or $s > $ 0 on the limit. The 
complete manifold $N$ is weakly embedded in $\Omega_{\varepsilon},$ in 
the sense that any smooth compact domain $K \subset  N$ embeds in 
$\Omega_{\varepsilon},$ provided $\varepsilon $ is sufficiently small, 
depending on $K$. Consequently, $N$ weakly embeds in $M$.
\end{theorem}

\medskip

 The blow-up limits $(N, g' )$ model the small-scale degeneration of 
the sequence $(\Omega_{\varepsilon}, g_{\varepsilon})$ in neighborhoods 
of the base points $y_{\varepsilon}.$ Theorem 1.3 implies in particular 
that the metrics $g_{\varepsilon}'$ in (1.19) do not collapse, in the 
sense of Cheeger-Gromov, near the base points $y_{\varepsilon}$. The 
potential function $\tau $ in (1.16)-(1.17) is a limit of the functions 
$\tau_{\varepsilon}, \varepsilon  = \varepsilon_{i},$ in (1.7), or 
suitable renormalizations thereof.

\begin{theorem} \label{t 1.4.} {\bf (Asymptotically Flat Ends).} 
 Let $(N, g, \tau )$ be a complete non-flat $\mathcal{Z}_{c}^{2}$ 
solution, i.e a solution of (1.16)-(1.17). Suppose there exists a 
compact set $K \subset  N$ and a constant $\omega_{o} < $ 0 such that 
the potential function
$$\omega  = \tau +\frac{\alpha}{12}s: N \rightarrow  {\mathbb R} , $$
in (1.17) satisfies
\begin{equation} \label{e1.20}
\omega  \leq  \omega_{o} <  0, 
\end{equation}
on $N\setminus K$, and that the level sets of $\omega $ in $N$ are 
compact.

 Then $N$ is an open 3-manifold of the form
\begin{equation} \label{e1.21}
N = P\# (\#_{1}^{q}{\mathbb R}^{3}), 
\end{equation}
where $P$ is a closed 3-manifold, (possibly empty), admitting a metric 
of positive scalar curvature, i.e. $\sigma (P) > $ 0, and $1 \leq q <  
\infty .$

 Each end $E = E_{k}$, $1 \leq  k \leq  q$, of $N$ is asymptotically 
flat in the sense of Definition 0.3, and the potential $\omega $ has 
the expansion
\begin{equation} \label{e1.22}
\omega  = \omega_{E} + \frac{m_{E}|\omega_{E}|}{r} + O(r^{-2}), 
\end{equation}
where $\omega_{E} < $ 0 is a constant depending on $(E, g)$ and $m_{E} 
> $ 0 is the mass of the end $E$.
\end{theorem}

\medskip

 Theorem 1.4 characterizes the complete $\mathcal{Z}_{c}^{2}$ solutions 
for which all ends are asymptotically flat, (except possibly when the 
potential $\omega $ goes to 0 at infinity in some end).

 The simplest example of a metric satisfying the conclusions of Theorem 
1.4 is the Schwarzschild metric, (on the space-like hypersurface), 
\begin{equation} \label{e1.23}
g_{S} = (1 -  \frac{2m}{r})^{-1}dr^{2} + r^{2}ds^{2}_{S^{2}}, 
\end{equation}
defined on $[2m,\infty )\times S^{2},$ and isometrically doubled across 
the horizon $\Sigma  = \{r = 2m\}$; observe that $\Sigma $  is a 
totally geodesic 2-sphere, of constant curvature $(2m)^{-2}.$ The 
metric $g_{S}$ is asymtotically flat at each end. We refer to 
[3, Prop. 5.1] for the exact form of the potential $\tau $ for this metric.

 Note that if blow-ups $\{g_{\varepsilon}'\}$ of $\{g_{\varepsilon}\}$ 
converge to the Schwarzschild metric, then the metrics 
$\{g_{\varepsilon}\}$ themselves are collapsing or crushing the core 
$S^{2}$ in $g_{S}$ to a point. Similarly, this occurs for the 2-spheres 
in any asymptotically flat end of a blow-up limit. If such $S^{2}$'s 
are essential, the metrics $\{g_{\varepsilon}\}$ are, in effect, then 
performing or carrying out a process analogous to the sphere 
decomposition (0.1) of $M$; c.f. the proof of Theorem 2.14 for a 
concrete illustration of this. 

\medskip

  We close this background section with a discussion of the regularity 
of the metrics $(\Omega_{\varepsilon}, g_{\varepsilon})$ as 
$\varepsilon \rightarrow 0$ when $M$ is tame. Thus, the $L^{2}$ bound 
on the curvature (0.15) may not apriori control the $L^{\infty}$ norm 
of the curvature or even the $L^{2}$ curvature radius $\rho$ of 
$(\Omega_{\varepsilon}, g_{\varepsilon})$ from (1.1). As an example, 
consider the following family of 2-dimensional metrics suggested by 
Gallot [10]: for any $\varepsilon  > $ 0,
$$h_{\varepsilon} = dr^{2} + (\varepsilon +r)^{2k}d\theta^{2}, $$
on $(- 1, 1)\times S^{1}.$ For $k > $ 3, the metric $h_{\varepsilon}$ 
satisfies the bound (0.15), with $\Lambda  = k^{2}(k- 1)^{2}/(k- 3).$ 
However, the curvature blows up in $L^{\infty}$ and the $L^{2}$ 
curvature radius goes to 0, in that
$$\rho_{\varepsilon}(0,\theta ) \rightarrow  0 \ {\rm as} \  
\varepsilon  \rightarrow  0, $$
for any $\theta\in S^{1}.$ The metrics $h_{\varepsilon}$ collapse at $r 
=$ 0 as $\varepsilon  \rightarrow $ 0, in the sense that the volume 
radius goes to 0, (but do not collapse at any $r \neq $ 0). Observe 
that this transition from non-collapse to collapse takes place within 
bounded distance, so that the limit metric is incomplete. Of course 
there are similar metrics in dimensions $\geq 3$.

 However, for the ``special'' metrics $g_{\varepsilon},$ this 
phenomenon does not occur, as shown in the following Lemma.

\begin{lemma} \label{l 1.5.}
  For any minimizing pair $(\Omega_{\varepsilon}, g_{\varepsilon}),$ 
there is a constant $\rho_{o} = \rho_{o}(\Lambda ) > $ 0 such that if 
(0.15) holds, then, for all $x_{\varepsilon} \in (\Omega_{\varepsilon}, 
g_{\varepsilon})$,
\begin{equation} \label{e1.24}
\rho(x_{\varepsilon}) \geq \rho_{o}.
\end{equation}
\end{lemma}

\noindent
{\bf Proof:} 
 The proof is by contradiction, so suppose (0.15) holds, but (1.24) does 
not, on some sequence $\varepsilon = \varepsilon_{j} \rightarrow 0$. 
Then by [5, Thm.4.5], (a slightly stronger version of Theorem 1.3 
above), there is a sequence of base points 
$y_{\varepsilon}\in\Omega_{\varepsilon},$ for some sequence of 
minimizing pairs $(\Omega_{\varepsilon}, g_{\varepsilon}), \varepsilon  
= \varepsilon_{j} \rightarrow $ 0, with 
\begin{equation} \label{e1.25}
\rho (y_{\varepsilon}) \rightarrow  0, 
\end{equation}
such that the blow-up metrics $g_{\varepsilon}'  = \rho 
(y_{\varepsilon})^{-2}\cdot  g_{\varepsilon}$ based at 
$y_{\varepsilon},$ have a subsequence converging to a complete, 
non-flat limit $(N, g')$ having the properties in Theorem 1.3. 
We point out again that this result implies in particular that 
the metrics $g_{\varepsilon}'$ do not collapse near $y_{\varepsilon}$ 
in this situation. The convergence is in the strong $L^{2,2}$, (in fact 
$C^{2, \alpha}$), topology on compact subsets of $N$, (c.f (1.26) below).

  Now by the scaling properties of curvature and volume, one has
$$\int_{\Omega_{\varepsilon}}|z_{g_{\varepsilon}'}|^{2}dV_{g_{\varepsilon}'} 
= \rho (y_{\varepsilon})\cdot 
\int_{\Omega_{\varepsilon}}|z_{g_{\varepsilon}}|^{2}dV_{g_{\varepsilon}}
 \rightarrow  0, $$
by (0.15) and (1.25). It follows that the limit $(N, g' )$ has $z_{g'} 
=$ 0, and hence $(N, g' )$ is flat, a contradiction. 
{\endproof}

 With regard to higher order regularity, by [5, Thm.4.2/Rmk.4.3] 
there is a constant $d_{o} > $ 0, independent of $\varepsilon ,$ such 
that
\begin{equation} \label{e1.26}
\rho_{\varepsilon}^{1,p}(x) \geq  d_{o}\cdot \rho_{\varepsilon}(x), 
\end{equation}
for all $x\in\Omega_{\varepsilon};$ here $d_{o}$ depends only on $p$ 
and $\rho^{1,p}$ is the $L^{1,p}$ curvature radius. This means the 
metric is controlled locally in $L^{3,p}$, in harmonic coordinates, in 
any region where it is controlled in $L^{2,2}$, (away from the boundary). 
By Sobolev embedding, $L^{3,p} \subset C^{2,\alpha}$, for any $\alpha = 
\alpha(p) < 1$. In particular, (1.26) implies that the curvature $|r|$ 
satisfies
\begin{equation} \label{e1.27}
|r_{\varepsilon}|(x) \leq  K\cdot \rho_{\varepsilon}(x)^{-2}, 
\end{equation}
where $K$ depends only on the choice of $c_{o}$ in (1.1).  
Hence for tame 3-manifolds, there is a constant $\lambda  = 
\lambda (\Lambda ) <  \infty $ such that $|r_{\varepsilon}| \leq  
\lambda $ on some, (in fact any), sequence of minimizing pairs 
$(\Omega_{\varepsilon}, g_{\varepsilon}),$ with 
$\varepsilon  = \varepsilon_{i} \rightarrow $ 0. Further, again by 
[5, Rmk.4.3], $|\nabla^{k}r_{\varepsilon}|$ is then bounded as 
$\varepsilon \rightarrow 0$, for any $k < \infty$, in regions where 
the potential $\tau = \tau_{\varepsilon} < 0$ or where 
$s_{g_{\varepsilon}} > 0$. Of course these estimates come from the fact 
that $g_{\varepsilon}$ satisfies the elliptic system (1.4)-(1.5); the 
main point is then to show that the estimates are independent of 
$\varepsilon$.

\section{Geometrization of Tame 3-Manifolds.}
\setcounter{equation}{0}

 In this section, we prove Theorem 0.2. Thus, throughout this section, 
it is assumed that $M$ is a tame, closed oriented 3-manifold with $\sigma 
(M) \leq $ 0, so that there is a sequence $\varepsilon  = 
\varepsilon_{i} \rightarrow $ 0, and a sequence of minimizing pairs 
$(\Omega_{\varepsilon}, g_{\varepsilon})$, $\varepsilon  = 
\varepsilon_{i}$, $\Omega_{\varepsilon} \subset\subset  M$, such that the 
$L^{2}$ curvature is uniformly bounded, i.e.
\begin{equation} \label{e2.1}
\int_{\Omega_{\varepsilon}}|z_{\varepsilon}|^{2}dV_{g_{\varepsilon}} 
\leq  \Lambda , 
\end{equation}
for some $\Lambda  <  \infty .$ It follows from Lemma 1.5 that
\begin{equation} \label{e2.2}
\rho_{\varepsilon}(x) \geq  \rho_{o}, 
\end{equation}
for some $\rho_{o} = \rho(\Lambda) >  0$.

  In \S 2.1, we prove various convergence and collapse results needed 
for the proof of Theorem 0.2; based on these, the proof will be 
completed in \S 2.2. A detailed discussion of the geometrization of 
graph manifolds from this point of view is then given in \S2.3.

\medskip

{\bf \S 2.1.}
We apply the $L^{2},$ (or in view of (2.2) and (1.27), the 
$L^{\infty}),$ Cheeger-Gromov theory as described in [2, \S 3] to the 
sequence $(\Omega_{\varepsilon}, g_{\varepsilon}),$ with $\varepsilon  
= \varepsilon_{i} \rightarrow 0.$ There are exactly three possible 
behaviors for the sequence $(\Omega_{\varepsilon}, g_{\varepsilon});$ 
namely a subsequence either converges, collapses or forms cusps. These 
cases are distinguished by the possible behaviors of the volume radius 
$\nu_{\varepsilon}$ as $\varepsilon  \rightarrow $ 0, in that either:

\begin{itemize}
\item
$\nu_{\varepsilon}(x) \geq  \nu_{o}, \forall x\in\Omega_{\varepsilon}$, 
\ \ for some $\nu_{o} > 0$.
\item
$\nu_{\varepsilon}(x) \rightarrow  0, \forall 
x\in\Omega_{\varepsilon}$, \ \  as $\varepsilon  \rightarrow  0$. 
\item
There exist points $x_{\varepsilon}$ and $y_{\varepsilon}$ in 
$\Omega_{\varepsilon}$ such that
\end{itemize}
$$\nu_{\varepsilon}(x_{\varepsilon}) \geq  \nu_{o} \ \ {\rm and} \ \ 
\nu_{\varepsilon}(y_{\varepsilon}) \rightarrow  0, \ \ {\rm as} \ \ 
\varepsilon  \rightarrow  0.$$

 We treat each of these cases separately.
\begin{proposition} \label{p 2.1.}
{\bf (Convergence).}
  Let $(\Omega_{\varepsilon}, g_{\varepsilon})$ be a sequence of 
minimizing pairs satisfying (2.1), for some sequence $\varepsilon  = 
\varepsilon_{i} \rightarrow $ 0. Suppose that there is a uniform 
constant $\nu_{o} > $ 0 such that
\begin{equation} \label{e2.3}
\nu_{\varepsilon}(x) \geq  \nu_{o} >  0, 
\end{equation}
for all $x \in  (\Omega_{\varepsilon}, g_{\varepsilon}).$ Then
$$\Omega_{\varepsilon} = M, $$
for all $\varepsilon  > $ 0 small, and a subsequence of $(M, 
g_{\varepsilon})$ converges to a constant curvature metric $g_{o}$ on 
$M$, of scalar curvature $\sigma (M)/6$ and volume 1.
\end{proposition}

\noindent
{\bf Proof:}
 By Theorem 1.1, if $\Omega_{\varepsilon} \subset\subset  M$ but 
$\Omega  \neq  M$, then $g_{\varepsilon}$ collapses, (for $\varepsilon$ 
fixed), with bounded curvature everywhere at infinity in 
$(\Omega_{\varepsilon}, g_{\varepsilon}).$ This implies that 
$\nu_{\varepsilon}(x) \rightarrow $ 0, as $x$ diverges to infinity in 
$\Omega_{\varepsilon},$ contradicting (2.3). Thus (2.3) implies that 
$\Omega_{\varepsilon} = M$, and $g_{\varepsilon}$ is globally defined 
on $M$.

 Note also that there is a uniform bound on the diameter of $(M, 
g_{\varepsilon}),$ i.e.
\begin{equation} \label{e2.4}
diam_{g_{\varepsilon}}M \leq  D, 
\end{equation}
where $D$ depends only on $\nu_{o}.$ This follows since 
$vol_{g_{\varepsilon}}M =$ 1, and by (2.3), there is a uniform lower 
bound on $vol_{g_{\varepsilon}}B_{x}(\nu_{o}),$ for any $x\in M.$ 
Hence, there is a uniform bound on the number of disjoint geodesic 
balls of radius $\nu_{o}$ in $(M, g_{\varepsilon}).$

 Since the curvature of $(M, g_{\varepsilon})$ is uniformly bounded in 
$L^{2},$ it follows by the $L^{2}$ Cheeger-Gromov theory that a 
subsequence converges, (modulo diffeomorphisms of $M$), in the weak 
$L^{2,2}$ topology to a limit $L^{2,2}$ metric $g_{o}$ defined on $M$, 
(c.f. [2,Thm.3.7]). In fact, by (2.2) and (1.26), the 
curvature of $g_{\varepsilon}$ is uniformly bounded in $L^{1,p},$ for 
any fixed $p <  \infty ,$ and hence the convergence is actually in the 
weak $L^{3,p}$ topology. By the compactness of the embedding $C^{2,\alpha} 
\subset  L^{3,p},$ it follows that the convergence of $g_{\varepsilon}$ 
to $g_{o}$ is in the $C^{2,\alpha}$ topology. In particular, the 
limit metric $g_{o}$ is at least $L^{3,p}$ smooth.

 The limit metric $g_{o}$ is a weak solution to the limit equations
\begin{equation} \label{e2.5}
L^{*}\tau  -  (\frac{1}{4}s\tau  + \frac{1}{12}\sigma (M))\cdot  g = 0, 
\end{equation}
\begin{equation} \label{e2.6}
2\Delta\tau  + \frac{1}{4}s\tau  = \frac{1}{4}|\sigma (M)|. 
\end{equation}
obtained by taking the limit as $\varepsilon_{i} \rightarrow $ 0 in 
(1.4)-(1.5), c.f. also (0.13)-(0.14). The equations (2.5)-(2.6) are 
just the Euler-Lagrange equations for a critical point of the 
functional $\mathcal{S}  _{-}^{2}.$ 

 By definition, i.e. from (1.7), on $(M, g_{\varepsilon})$ one has 
\begin{equation} \label{e2.7}
\int_{M}\tau_{\varepsilon}^{2}dV_{\varepsilon} = 1. 
\end{equation}
Since the local and global geometry of $(M, g_{\varepsilon})$ is 
uniformly bounded, standard elliptic theory applied to the trace 
equation (1.5), c.f. [11, \S 8], implies that $\tau_{\varepsilon}$ is 
uniformly bounded in $L^{1,p}, p < \infty$. Hence, the limit function 
$\tau$ is in $L^{1,p}$, and is non-positive since $\tau_{\varepsilon}$ 
is. In particular, (2.7) also holds for the limit function $\tau$.

 Suppose first that $\sigma (M) < $ 0. We may apply the minimum 
principle to the trace equation (2.6) in a neighborhood of a point $q$ 
realizing the minimal value of $\tau .$ Since $\tau  \leq $ 0 and the 
$L^{2}$ average of $\tau $ over $M$ is 1, $\tau (q) \leq  - 1.$ Since 
then $(s\tau )(q) \geq  |\sigma (M)|,$ the minimum principle implies 
that
\begin{equation} \label{e2.8}
\tau  \equiv  -1, s \equiv  \sigma (M). 
\end{equation}
Alternately, (2.8) can be derived from (1.13) and (1.15). From (2.5) 
and the definition of $L^{*}$ in (1.6), it follows that
$$- r + \frac{1}{3}\sigma (M) = 0, $$
and hence
\begin{equation} \label{e2.9}
z = 0. 
\end{equation}
Thus, $(M, g_{o})$ is a smooth unit volume metric, of constant 
sectional curvature $\sigma (M)/6.$ The convergence to the limit is 
then in the $C^{\infty}$ topology, by the remarks following (1.27).

 Next suppose $\sigma (M) =$ 0. From (1.14), one sees that the limit 
scalar curvature $s$ is non-negative. Since the limit function $\tau$ 
is non-positive, it follows that $s \cdot \tau \equiv 0$ on the limit. 
Hence the trace equation (2.6) gives $\Delta\tau  =$ 0 on $(M, g_{o}),$ 
so that again by the minimum principle, $\tau $ is constant. Hence, by 
(2.7), we must have $\tau  \equiv  - 1, s = \sigma (M) =$ 0, and (2.8) 
holds in this situation also. Thus (2.5) again implies that $(M, 
g_{o})$ is a flat 3-manifold of unit volume.

{\endproof}

 In case $\sigma (M) < $ 0, the metric $g_{o}$ is unique, (up to 
isometry), by the Mostow rigidity theorem [25]. Of course, if $\sigma 
(M) =$ 0, i.e. $g_{o}$ is a flat metric, one cannot expect uniqueness, 
since the moduli space of flat metrics on $M$ may well be non-trivial.

\begin{remark} \label{r 2.2.}
  It follows that under the assumption (2.3), $\inf \mathcal{S}_{-}^{2}$ 
is realized by a constant curvature metric $g_{o}$ on $M$. Clearly 
$I_{\varepsilon}^{~-}(g) \geq  \mathcal{S}  _{-}^{2}(g),$ for all 
$g\in{\mathbb M} $ and $\varepsilon  \geq $ 0, with equality if and only 
if $z_{g} =$ 0. In particular,
\begin{equation} \label{e2.10}
I_{\varepsilon}^{~-}(g_{o}) = I_{o}^{-}(g_{o}) = \mathcal{S}  
_{-}^{2}(g_{o}) = |\sigma (M)|. 
\end{equation}
In case $\sigma (M) < $ 0, it follows from (0.4) and Mostow rigidity 
that the constant curvature metric $g_{o}$ uniquely realizes $\inf 
I_{\varepsilon}^{~-}$, for all $\varepsilon  \geq $ 0, among metrics 
on $M$. Thus, the ``family'' of metrics $g_{\varepsilon}$ satisying 
(2.3) is constant, i.e. 
\begin{equation} \label{2.11}
g_{\varepsilon} = g_{o},
\end{equation}
for all $\varepsilon  \geq $ 0. More precisely, the isometry class of 
the metrics $g_{\varepsilon}$ is constant, in the sense that for any 
$\varepsilon  > $ 0, there is a diffeomorphism $\psi_{\varepsilon}$ of 
$M$ such that $\psi_{\varepsilon}^{*}g_{\varepsilon} = g_{o}.$ This 
means of course that the sequence $(\Omega_{\varepsilon}, 
g_{\varepsilon})$ in Proposition 2.1 is in fact a constant sequence, 
modulo diffeomorphisms. Note however that this discussion does not 
automatically preclude the existence of other minimizers 
$(\Omega_{\varepsilon}, g_{\varepsilon})$ of $I_{\varepsilon}^{~-}$ 
with $\Omega_{\varepsilon} \subset\subset  M$, for which (2.3) does not 
hold.

 Similarly, if $\sigma (M) =$ 0, the metrics $g_{\varepsilon}$ 
satisfying (2.3) are all flat metrics on $M$, for all $\varepsilon  
\geq $ 0. Since the moduli space $\mathcal{M}_{F}$ of flat metrics on $M$ 
may be non-trivial, the metrics $g_{\varepsilon}$ may not be unique. 
However, the volume assumption (2.3) implies that the metrics 
$g_{\varepsilon}$ remain within a compact set of the moduli space 
$\mathcal{M}_{F}.$

\end{remark}

\begin{proposition} \label{p 2.3.}
{\bf (Collapse).}
  Let $(\Omega_{\varepsilon}, g_{\varepsilon})$ be a sequence of 
minimizing pairs satisfying (2.1), for some sequence $\varepsilon  = 
\varepsilon_{i} \rightarrow $ 0. Suppose that
\begin{equation} \label{e2.12}
\nu_{\varepsilon}(x) \rightarrow  0, 
\end{equation}
for all $x \in  (\Omega_{\varepsilon}, g_{\varepsilon}).$ Then $M$ is a 
graph manifold with
$$\sigma (M) = 0. $$
\end{proposition}

\noindent
{\bf Proof:}
 Under the assumption (2.12), $(\Omega_{\varepsilon}, g_{\varepsilon})$ 
becomes arbitrarily thin at every point. Recall from \S 1 that the 
complement of a sufficiently large compact set $K_{\varepsilon}$ in 
$\Omega_{\varepsilon}$ admits a polarized F-structure 
$\mathcal{F}_{\infty} = \mathcal{F}_{\infty}(\varepsilon)$, 
i.e. a graph manifold structure, and that this structure extends 
to the rest of $M$, so that $V_{\varepsilon} = M\setminus K_{\varepsilon}$ 
admits a polarized F-structure, also called $\mathcal{F}_{\infty}.$ 

 Now from Lemma 1.5, (1.27) and the Cheeger-Gromov theory [7,8], 
there is a $\delta_{1} = \delta_{1}(\Lambda )$ such that if 
\begin{equation} \label{e2.13}
\nu_{\varepsilon}(x) \leq  \delta_{1}, 
\end{equation}
for all $x \in  (\Omega_{\varepsilon}, g_{\varepsilon}),$ then 
$\Omega_{\varepsilon}$ admits a polarized F-structure 
$\mathcal{F}_{\varepsilon}.$ It is clear that the F-structures 
$\mathcal{F}_{\infty}$ and $\mathcal{F}_{\varepsilon}$ are compatible on their 
intersections, and thus define a global F-structure $\mathcal{F} $ on $M$.

 Thus, (2.12), or even the weaker condition (2.13), implies that $M$ 
itself is a graph manifold. We have already noted in (1.10) that in 
this case, for any $\varepsilon  > $ 0, 
\begin{equation} \label{e2.14}
\inf I_{\varepsilon}^{~-} = \inf \mathcal{S}  _{-}^{2} = \sigma (M) = 0,
\end{equation}
which proves the result. In fact, as noted in \S 0, any minimizing 
sequence for $I_{\varepsilon}^{~-}$, with $\varepsilon > 0$ fixed, 
collapses $M$ along an F-structure, (or possibly a sequence of 
F-structures), unless $M$ is a flat manifold. In particular, the 
estimate (2.12) implies then that either
\begin{equation} \label{e2.15}
\Omega_{\varepsilon} = \emptyset  , \ \ \forall\varepsilon  >  0,
\end{equation}
or $M = \Omega_{\varepsilon}$ and $g_{\varepsilon}$ is a collapsing 
sequence of flat metrics on $M$.

{\endproof}

\begin{remark} \label{r 2.4.}
  We point out that this discussion shows that any closed graph 
manifold $G$ necessarily satisfies $\sigma (G) \geq $ 0. Hence if 
$\sigma(M) < 0$, then $M$ cannot be a graph manifold.
\end{remark}

\begin{proposition} \label{p 2.5.}
{\bf (Cusps I).}
  Suppose $\sigma (M) < $ 0, and let $(\Omega_{\varepsilon}, 
g_{\varepsilon})$ be a sequence of minimizing pairs satisfying (2.1), 
for some sequence $\varepsilon  = \varepsilon_{i} \rightarrow $ 0. 
Suppose that there exist $x_{\varepsilon}$ and $y_{\varepsilon}$ in 
$\Omega_{\varepsilon}$ such that
\begin{equation} \label{e2.16}
\nu_{\varepsilon}(x_{\varepsilon}) \geq  \nu_{o} \ \ {\rm and} \ \ 
\nu_{\varepsilon}(y_{\varepsilon}) \rightarrow  0, \ \ {\rm as} \ \ 
\varepsilon  \rightarrow  0, 
\end{equation}
for some $\nu_{o} > $ 0. Then there is a subsequence, also denoted 
$\{\varepsilon\}$, of $\{\varepsilon_{i}\}$, such that
$$\Omega_{\varepsilon} = \Omega $$
is independent of $\varepsilon ,$ and is given by a finite union of 
complete, non-compact hyperbolic manifolds of finite volume. There is 
an embedding 
$$\Omega  \subset  M, $$
for which the complement $M\setminus \Omega $ is a graph manifold. The 
metrics $g_{\varepsilon}$ converge to the complete constant curvature 
metric $g_{o},$ of scalar curvature $\sigma (M)/6,$ and
$$vol_{g_{o}}\Omega  = 1, $$
while collapsing the graph manifold $M\setminus \Omega $ to a lower 
dimensional space. 
\end{proposition}

\noindent
{\bf Proof:}
 Suppose the sequence $\{(\Omega_{\varepsilon}, g_{\varepsilon})\}, 
\varepsilon  = \varepsilon_{i},$ has basepoints $\{x_{\varepsilon}\}$ 
and $\{y_{\varepsilon}\}$ satisfying (2.16), and consider the pointed 
sequence $\{(\Omega_{\varepsilon}, g_{\varepsilon}, 
x_{\varepsilon})\}.$ As in the proof of Proposition 2.1, from Lemma 1.5 
and (1.27), the curvature of $(\Omega_{\varepsilon}, g_{\varepsilon})$ 
is uniformly bounded, independent of $\varepsilon .$ It follows from 
the $L^{\infty}$ Cheeger-Gromov theory, (c.f. [2,\S 2] and 
[3,Rmk.5.5]), that there are diffeomorphisms $\psi_{\varepsilon}$ of 
$\Omega_{\varepsilon},$ and a maximal connected open domain $\Omega_{o} 
\subset\subset  M$ such that a subsequence of 
$(\psi_{\varepsilon})^{*}g_{\varepsilon}, \varepsilon  = 
\varepsilon_{i}$, converges weakly in the $L^{2,p}$ topology based at 
$x_{\varepsilon},$ and uniformly on compact subsets, to a limit 
$L^{2,p}$ metric $g_{o}$ on $\Omega_{o},$ with
\begin{equation} \label{e2.17}
0 <  vol_{g_{o}}\Omega_{o} \leq  1. 
\end{equation}
By (1.26), the convergence is actually in the weak $L^{3,p}$ and strong 
$C^{2,\alpha}$ topologies, and the limit metric $g_{o}$ is $L^{3,p}.$ 
The curvature of $(\Omega_{o}, g_{o})$ is uniformly bounded in 
$L^{\infty}$, and since $(\Omega_{\varepsilon}, g_{\varepsilon})$ is 
complete for all $\varepsilon ,$ so is $(\Omega_{o}, g_{o}).$ As in 
Proposition 2.1, the metric $g_{o}$ is a (weak) solution of the 
Euler-Lagrange equations (2.5)-(2.6) for $\mathcal{S}  _{-}^{2}.$

 Since $\sigma (M) < $ 0, it follows from (1.13) that 
$\tau_{\varepsilon} \rightarrow  - 1$ almost everywhere on 
$\Omega_{\varepsilon}$ as $\varepsilon  \rightarrow $ 0. Hence the 
limit function $\tau  =$ lim $\tau_{\varepsilon}$ satisfies
\begin{equation} \label{e2.18}
\tau  \equiv  - 1, s \equiv  \sigma (M), 
\end{equation}
as in (2.8). The same reasoning as in (2.8)-(2.9) shows that
\begin{equation} \label{e2.19}
z_{g_{o}} = 0, 
\end{equation}
so that $(\Omega_{o}, g_{o})$ is a complete, non-compact manifold of 
constant sectional curvature $\sigma (M)/6$ and volume at most 1.

 We now fix the subsequence of $g_{\varepsilon}, \varepsilon  = 
\varepsilon_{i}$ above, but do not change the notation for the 
subsequence. Suppose there is another sequence of points 
$\{x_{\varepsilon}^{1}\}$ in $\Omega_{\varepsilon}$, with
\begin{equation} \label{e2.20}
\nu_{\varepsilon}(x_{\varepsilon}^{1}) \geq  \nu_{o}, 
\end{equation}
for some constant $\nu_{o} > $ 0, and $\{x_{\varepsilon}^{1}\}$ does 
not (sub)-converge to a point in $\Omega_{o}.$ Then one may repeat the 
above process to obtain another smooth maximal open domain $\Omega_{1} 
\subset\subset  M$, with smooth complete metric $g_{o},$ which again is 
of constant curvature, and scalar curvature $s = \sigma (M) < $ 0.

 This process may be repeated as many times as necessary until there 
are no sequences left satisfying (2.20), and not subconverging to a 
previously defined domain. There is clearly a uniform upper bound $Q$ 
on the number of domains $\Omega_{j}, j =$ 0, 1, ..., $Q$, obtained in 
this way, since $vol \Omega_{j} \geq  10^{-1}\nu_{o}^{3}$ and the total 
volume of $\Omega  = \cup\Omega_{j}$ is at most 1. In fact, the bound 
$Q$ is independent of the choice of $\nu_{o}$ in (2.20) or (2.16), and 
depends only on the 3-manifold $M$. For if $(U, g_{o})$ is any 
complete, constant curvature, open 3-manifold, with scalar curvature 
$\sigma (M),$ then
\begin{equation} \label{e2.21}
vol_{g_{o}}U \geq  v_{o} >  0, 
\end{equation}
where $v_{o}$ depends only on an upper bound for $|\sigma (M)|; v_{o}$ 
is essentially the Margulis constant, c.f. [36] or [37]. Thus, 
choosing $\nu_{o}$ above small, i.e. $\nu_{o} <<  v_{o},$ it follows 
that there is a fixed bound $Q = Q(M)$ on the number of domains 
$\Omega_{j}$ obtained by the process above.

 Let
\begin{equation} \label{e2.22}
\Omega  = \bigcup _{o}^{Q}\Omega_{j}. 
\end{equation}
Then $(\Omega , g_{o})$ is a collection of a finite number of 
connected, noncompact, complete Riemannian manifolds, each with 
constant curvature metric $g_{o}$ with $s_{g_{o}} = \sigma (M) < $ 0. 
In particular, each component $\Omega_{j}$ is of finite topological 
type, having a finite number of ends, i.e. cusps, each diffeomorphic to 
$T^{2}\times {\mathbb R}^{+}.$ It follows that not only is each component 
$\Omega_{j}$ weakly embedded in $M$, but $\Omega_{j}$ is actually 
embedded in $M$ as an open domain. Consequently, there is an embedding 
$\Omega  \subset M$. The complement $M\setminus \Omega $ is a compact 
manifold with boundary consisting of a finite number of tori, which has 
a sequence of F-structures $\mathcal{F}_{\varepsilon}$ along which a 
subsequence of the metrics $g_{\varepsilon}, \varepsilon  = 
\varepsilon_{i},$ collapse $M\setminus \Omega $ to a lower dimensional 
space. In particular, $M\setminus \Omega $ is a graph manifold with 
toral boundary.

 By the construction of $\Omega$,
\begin{equation} \label{e2.23}
vol_{g_{o}}\Omega  \leq  1. 
\end{equation}
We claim that
\begin{equation} \label{e2.24}
vol_{g_{o}}\Omega  = 1. 
\end{equation}
The proof of this is exactly the same as the proof of (1.9), given in 
[2, Thm. 5.7]. Briefly, suppose $vol_{g_{o}}\Omega  < $ 1, so that 
$vol_{g_{o}}\Omega  \leq  1-\mu ,$ for some $\mu  > $ 0. Using the 
structure of graph manifolds, for any $\delta  > $ 0, one may construct 
a metric $\bar g = \bar g_{\delta}$ on $M$, agreeing with $g_{o}$ on a 
prescribed sufficiently large compact set $K \subset  \Omega ,$ such 
that $vol_{\bar g}M\setminus K \leq  \delta $ and $|\bar s| \leq  C$ on 
$M\setminus K$, where $C$ is independent of $K$ and $\delta .$ Thus, 
choosing $\delta $ much smaller than $\mu $ implies that
$$v(\bar g)^{1/3}\int_{M}(s_{\bar g}^{-})^{2}dV_{\bar g} <  
v(g_{o})^{1/3}\int_{\Omega}(s_{g_{o}}^{-})^{2}dV_{g_{o}} = \sigma 
(M)^{2}, $$
which contradicts (0.4). Finally, the arguments proving (2.11) in 
Remark 2.2 also prove that $\Omega = \Omega_{\varepsilon}$ and $g_{o} = 
g_{\varepsilon}$ in this situation, by means of the Mostow-Prasad 
rigidity theorem [27]. This completes the proof.

{\endproof}

\begin{remark} \label{r 2.6.}
  As in Remark 2.2, one sees that 
$\inf I_{\varepsilon}^{~-} = \inf \mathcal{S}_{-}^{2} = |\sigma (M)|$ 
is realized by the union of the complete 
constant negative curvature metrics $g_{o}$ on $\Omega\subset M.$ 
Apriori however, as before, it may be possible that $(\Omega, g_{o})$ 
is not unique; this will be discussed further in \S 2.2.

\end{remark}

 We now turn to the analogue of Proposition 2.5 in case $\sigma (M) =$ 
0.

\begin{proposition} \label{p 2.7.}
{\bf (Cusps II).}
   Suppose $\sigma (M) =$ 0 and let $(\Omega_{\varepsilon}, 
g_{\varepsilon})$ be a sequence of minimizing pairs satisfying (2.1), 
for some sequence $\varepsilon  = \varepsilon_{i} \rightarrow $ 0. Then 
$(\Omega_{\varepsilon}, g_{\varepsilon})$ cannot form cusps, i.e. 
either there exists $\nu_{o} > $ 0 such that
$$\nu_{\varepsilon}(x) \geq  \nu_{o} >  0, \ \ {\rm for \ all} \ \ 
x\in\Omega_{\varepsilon} ,$$
or
$$\nu_{\varepsilon}(x) \rightarrow  0, \ \ {\rm for \ all} \ \ 
x\in\Omega_{\varepsilon}. $$

\end{proposition}

\noindent
{\bf Proof:}
 The proof is by contradiction and so we assume $(\Omega_{\varepsilon}, 
g_{\varepsilon})$ does form cusps as $\varepsilon = \varepsilon_{i} 
\rightarrow 0$, i.e. there are base points $x_{\varepsilon}, 
y_{\varepsilon}$ satisfying (2.16). The same arguments as in the proof 
of Proposition 2.5 prove that a subsequence of 
$\{(\Omega_{\varepsilon}, g_{\varepsilon}, x_{\varepsilon})\}$ 
converges to a complete connected maximal limit $(\Omega_{o}, g_{o}, 
x)$ with
\begin{equation} \label{e2.25}
\int_{\Omega_{o}}|z|^{2}dV_{o} \leq  \Lambda , \ \ 
vol_{g_{o}}\Omega_{o}  \leq  1. 
\end{equation}
By (1.14), $s_{g_{o}} \geq 0$ everywhere and by (1.27), $g_{o}$ has 
uniformly bounded curvature. Since $(\Omega_{o}, g_{o})$ is of finite 
volume, it collapses everywhere along an F-structure at infinity. In 
particular, a neighborhood of infinity of $\Omega_{o}$ is a graph 
manifold. Further, by continuity, $(\Omega_{o}, g_{o})$ minimizes the 
$L^{2}$ norm of $z$ among all compact perturbations of $g_{o}$ with 
non-negative scalar curvature and with volume at most that of $g_{o}.$

 As in the proof of Proposition 2.5, this process may be repeated for 
any other sequence of base points $\{x_{\varepsilon}^{j}\}$ satisfying 
(2.16) giving rise to disjoint cusps $\Omega_{j}$ in the limit. Each 
such $\Omega_{j}$ satisfies (2.25) with the same $\Lambda$. Again 
however this process terminates after a finite number $Q$ of steps. 
This is because if $j$ is sufficiently large, then $\Omega_{j}$ 
satisfies (2.25) and $vol\Omega_{j} \leq \delta_{1}$, where $\delta_{1} 
= \delta_{1}(j)$ may be made arbitrarily small by choosing $j$ 
sufficiently large. However, such a manifold is a graph manifold by 
(2.13). This implies that the minimizing pairs $(\Omega_{\varepsilon}, 
g_{\varepsilon}), \varepsilon = \varepsilon_{i}$, will have collapsed 
$\Omega_{j}$ for $\varepsilon$ sufficiently small, i.e. in effect 
$\Omega_{j} = \emptyset$, as in (2.15). Thus, as in (2.22), the maximal 
domain $\Omega = \cup \Omega_{j}$ consists of a finite number of 
connected components. For similar reasons, the same arguments proving 
(2.24), proves that
\begin{equation} \label{e2.26}
vol_{g_{o}}\Omega  = 1. 
\end{equation}

 As in Proposition 2.1, the limit metric $g_{o}$ and limit potential 
$\tau  =$ lim $\tau_{\varepsilon}$ are a solution of the Euler-Lagrange 
equations (2.5)-(2.6) with $s \geq $ 0. Since the region where $s > 0$ 
is disjoint from the region where $\tau < 0$, $s \cdot \tau \equiv 0$, 
and hence the limit function $\tau $ is harmonic on $(\Omega, g_{o})$. 
By (2.7), $\tau$ is in $L^{2}(\Omega, g_{o}).$ This is easily seen to 
imply, by a standard cutoff argument, that $\tau $ is constant on each 
component $\Omega_{j}$ of $\Omega$; (this can also be seen from the 
estimate (1.15)). 

 If the limit potential $\tau  \neq $ 0 on some $\Omega_{j},$ then as 
in (2.8)-(2.9), again one has $z =$ 0, and hence the complete limit 
$(\Omega_{j}, g_{o})$ is flat. However, any complete non-compact flat 
manifold must have infinite volume, contradicting (2.26). Thus, such 
cusps cannot form.

 It follows we must have $\tau  \equiv $ 0 on the limit $(\Omega , 
g_{o}).$ Although this situation is very special and unlikely to occur, 
it will take some further arguments to rule it out. (In case 
$\Omega_{j}$ is irreducible, for some $j$, these arguments can be 
bypassed, c.f. Remark 2.8 below).

\medskip

 Recall that the $L^{2}$ norm of $\tau_{\varepsilon}$ on 
$(\Omega_{\varepsilon}, g_{\varepsilon})$ is 1 by (2.7). Hence in this 
situation, all of $\tau_{\varepsilon}$ is concentrating at infinity in 
$\Omega_{\varepsilon}$ as $\varepsilon  \rightarrow $ 0, and so by 
(2.26), concentrating on a set of measure converging to 0. In 
particular, $\tau_{\varepsilon}$ must be unbounded in this region, so 
that $T_{\varepsilon} = \sup|\tau_{\varepsilon}| \rightarrow  \infty ,$ 
as $\varepsilon  \rightarrow $ 0. By (1.11), this forces
\begin{equation} \label{e2.27}
\frac{\varepsilon}{\sigma} \rightarrow  \infty , \ \ {\rm as} \ \ 
\varepsilon  \rightarrow  0, 
\end{equation}
where $\sigma  = \sigma (g_{\varepsilon}) = \mathcal{S}  
_{-}^{2}(g_{\varepsilon}).$ It is then natural to consider the 
renormalized functional
\begin{equation} \label{e2.28}
\frac{1}{\varepsilon}I_{\varepsilon}^{~-} = \mathcal{Z}^{2} + 
\frac{1}{\varepsilon}\mathcal{S}  _{-}^{2} 
\end{equation}
on $(\Omega_{\varepsilon}, g_{\varepsilon}),$ together with the 
renormalized Euler-Lagrange equations (1.4)-(1.5), i.e.
\begin{equation} \label{e2.29}
\nabla\mathcal{Z}^{2} + L^{*}\bar \tau + \bar \phi \cdot  g = 0, 
\end{equation}
\begin{equation} \label{e2.30}
2\Delta (\bar \tau+\frac{s}{12}) + \frac{1}{4}s\bar \tau = 
-\frac{1}{2}|z|^{2} + 3\bar c, 
\end{equation}
where $\bar \tau = \tau /\varepsilon , \bar c = c/\varepsilon $ and 
$\bar \phi = \phi /\varepsilon .$ By (2.27) and (2.26), $\bar c$ 
converges to $\frac{1}{6}\int_{\Omega}|z|^{2}dV_{o} \leq \Lambda$ as 
$\varepsilon  \rightarrow $ 0, and by the regularity following (1.27), 
$\nabla \mathcal{Z}^{2}$, $|z|^{2}$ and $s$ are uniformly bounded as 
$\varepsilon \rightarrow 0$.

 We divide the discussion now into three cases, according to the 
behavior of $\bar \tau = \bar \tau_{\varepsilon}$ as $\varepsilon  
\rightarrow $ 0, on a sequence of points $y_{\varepsilon} \in  
\Omega_{\varepsilon}$ converging to a limit point $y \in  \Omega .$ 
These behaviors are either $\bar \tau(y_{\varepsilon}) \rightarrow 
-\infty$, $\bar \tau(y_{\varepsilon})$ remains bounded, or $\bar \tau$ 
tends to $0$ in neighborhoods of $y_{\varepsilon}$ as $\varepsilon 
\rightarrow 0$. (This discussion closely resembles, at least formally, 
that in [5, \S 4.1]). Let $\Omega_{y}$ be the component of $\Omega $ 
containing $y$.

{\bf Case (i).}
 Suppose $\bar \tau(y_{\varepsilon}) \rightarrow  -\infty $ as 
$\varepsilon  \rightarrow $ 0. In this case, divide the equations 
(2.29)-(2.30) further by $|\bar \tau(y_{\varepsilon})|.$ Let 
$\widetilde \tau = \widetilde \tau_{\varepsilon} = \bar \tau/|\bar 
\tau(y_{\varepsilon})|.$ Since the right hand side of (2.30) is 
uniformly bounded, elliptic regularity, c.f. [11, Ch.8] implies that 
$\widetilde \tau_{\varepsilon}$ is uniformly bounded at points a 
bounded distance to the base points $y_{\varepsilon}.$ Hence, a 
subsequence converges to a limit function $\widetilde \tau.$ The limit 
(renormalized) equations on $(\Omega_{y} , g_{o})$ are then
\begin{equation} \label{e2.31}
L^{*}\widetilde \tau = 0, \ \  \Delta\widetilde \tau = 0, 
\end{equation}
i.e. the static vacuum Einstein equations, (c.f. [4]). Since 
$\widetilde \tau \leq $ 0, the maximum principle implies that 
$\widetilde \tau < $ 0 everywhere. Hence $(\Omega_{y}, g_{o}, y)$ is a 
complete solution to the static vacuum equations, with non-vanishing 
potential. By [4, Thm.3.2], it follows that $(\Omega_{y}, g_{o})$ is 
flat, which gives a contradiction as before.

{\bf Case (ii).}
 Suppose $\bar \tau(y_{\varepsilon})$ is bounded as $\varepsilon  
\rightarrow $ 0. The limit equations then take the form
\begin{equation} \label{e2.32}
\nabla\mathcal{Z}^{2} + L^{*}\bar \tau + \bar c\cdot  g = 0, 
\end{equation}
\begin{equation} \label{e2.33}
2\Delta(\bar \tau + \frac{s}{12}) = -\frac{1}{2}|z|^{2} + 3\bar c, 
\end{equation}
with $\bar c = \mathcal{Z}^{2}(g_{o})/6.$

 Now on the one hand, the metric $g_{\varepsilon}$ is a critical point 
of the functional $\frac{1}{\varepsilon}I_{\varepsilon}^{~-}.$ On the 
other hand, the second term in (2.28) converges to 0 on 
$g_{\varepsilon},$ by (2.27). Thus let $\eta$ be a positive cutoff 
function on $\Omega_{y}$, with $\eta \equiv 1$ on $B_{y}(R)$, $\eta 
\equiv 0$ on $\Omega_{y} \setminus B_{y}(2R)$, with $|d\eta| \leq c/R$, 
and consider for instance the variation of $g_{\varepsilon}$ given by
$$g_{\varepsilon ,t} = g_{\varepsilon} -  t\eta r_{\varepsilon}. $$
Note that since $r_{\varepsilon}$ is uniformly bounded as $\varepsilon  
\rightarrow $ 0, the metrics $g_{\varepsilon ,t}$ are well defined for 
all $t \leq  t_{o} = t_{o}(R),$ independent of $\varepsilon .$ 

 A straightforward computation shows that $\frac{d}{dt}\mathcal{S}  
_{-}^{2}(g_{\varepsilon ,t}) \leq $ 0, for $R$ sufficiently large. 
(Alternately, if $\frac{d}{dt}\mathcal{S}  _{-}^{2}(g_{\varepsilon ,t})$ 
were positive, replace $g_{\varepsilon ,t}$ by $g_{\varepsilon} + t\eta 
r$ in the ensuing argument). Hence 
$$0 \leq  \frac{1}{\varepsilon}\mathcal{S}  _{-}^{2}(g_{\varepsilon ,t}) 
\leq  \frac{1}{\varepsilon}\mathcal{S}  _{-}^{2}(g_{\varepsilon}),$$
for all $t$ small, and so both terms converge to 0 as $\varepsilon  
\rightarrow $ 0 by (2.27). Since the metrics converge smoothly to the 
limit, it follows that
\begin{equation} \label{e2.34}
\lim_{\varepsilon \rightarrow 0} \frac{d}{dt}\mathcal{S}_{-}^{2}
(g_{\varepsilon,t})|_{t=0} = 0.
\end{equation}
Thus, on the limit $(\Omega_{y}, g_{o})$, one has
\begin{equation} \label{e2.35}
0 = \int_{\Omega_{y}}< L^{*}\bar \tau, \eta z>  = \int_{\Omega_{y}}< 
D^{2}\bar \tau, \eta z>  -\int_{\Omega_{y}}\eta \bar \tau|z|^{2}.
\end{equation}
We note that $\int_{\Omega_{y}}< D^{2}\bar \tau, \eta z> \rightarrow 0$ 
as $R \rightarrow \infty$. This follows by integrating by parts, i.e. 
applying the divergence theorem, twice, together with the Bianchi 
identity and the fact that $s \cdot \bar \tau \equiv 0$

 Hence, (2.35) implies that either $z \equiv $ 0, in which case 
$(\Omega_{y}, g_{o})$ is flat, giving a contradiction as before, or the 
limit $\bar \tau \equiv $ 0. This is treated in the last case.

{\bf Case (iii).}
 Suppose $\bar \tau \equiv $ 0. In this case, the limit equations are
\begin{equation} \label{e2.36}
\nabla\mathcal{Z}^{2}+ \bar c\cdot  g = 0, 
\end{equation}
\begin{equation} \label{e2.37}
\frac{1}{6} \Delta s + \frac{1}{2}|z|^{2} = 3\bar c. 
\end{equation}
We will prove below that necessarily $\bar c = 0$. This then implies 
that $\int_{\Omega_{y}}|z|^{2}dV_{o} = 0$, so that $(\Omega_{y}, 
g_{o})$ is flat, and one has a contradiction as before.

 Now $(\Omega , g_{o})$ is of finite volume and collapsing at infinity. 
Hence, for any divergent sequence $\{x_{i}\}\in\Omega_{y} $ and any $R 
<  \infty ,$ the pointed sequence $(B_{x_{i}}(R), g_{o}, x_{i})$ in 
$(\Omega_{y} , g_{o})$ collapses with bounded curvature along an 
injective F-structure. We may then unwrap this collapse by passing to 
sufficiently large finite covers; choosing a sequence $R_{j} 
\rightarrow  \infty ,$ and a suitable diagonal subsequence gives rise 
to a complete limit solution $(N, g_{\infty}, x_{\infty})$ of 
(2.36)-(2.37) with a free isometric $S^{1}$ action. The convergence to 
such limits is smooth, by elliptic regularity applied to the equations 
(2.36)-(2.37), (c.f. [4, \S 3]). The constant $\bar c$ of course 
remains the same in passing to this geometric limit. Thus it suffices 
to evaluate $\bar c$ on these simpler manifolds $(N, g_{\infty})$.

 Let $V$ be the orbit space of the $S^{1}$ action, so that $V$ is a 
complete Riemannian surface. Let $f: V \rightarrow  {\mathbb R} $ denote 
the length of the $S^{1}$ fibers and $A$ the curvature form of the 
$S^{1}$ bundle. In [4, Prop. 4.1] it is proved, (via a Gauss-Bonnet 
type argument), that
$$\int_{V}|\nabla \log f|^{2} <  \infty , \ {\rm  and} \  \int_{V}|A|^{2} 
<  \infty . $$

  Suppose first that there exists $v_{o} > $ 0 and a sequence of points 
$p_{i}\in V$ such that $area D_{p_{i}}(1) \geq  v_{o},$ where 
$D_{p_{i}}(1)$ is the geodesic disc about $p_{i}$ of radius 1 in $V$. 
Then $\nabla \log f \rightarrow $ 0 and $A \rightarrow $ 0 in 
$D_{p_{i}}(R)$ as $i \rightarrow  \infty $, for any fixed $R$. This 
implies that any geometric limit $(N', g_{\infty}', p_{\infty})$ of 
$(N, g_{\infty}, p_{i})$ is a complete product metric of the form $V' 
\times S^{1}$, of non-negative scalar curvature, and satisfying 
(2.36)-(2.37), again with the same $\bar c$. It follows that the Gauss 
curvature of $V'$ is non-negative and hence there are points $q_{i} \in 
V'$ such that the metric $g_{\infty}'$ on $D_{q_{i}}(1) \subset V'$ 
converges to the flat metric, (in covers if necessary). This of course 
implies $\bar c = 0$, as required.

  On the other hand, suppose that $area D_{p_{i}}(1) \rightarrow 0$ for 
any divergent sequence $p_{i}$ in $V$, so that $V$ itself collapses 
everywhere at infinity. Then repeating the construction above on a 
divergent sequence, unwrapping the collapse as above, gives rise to a 
further complete limit $(N' , g_{\infty}' , x_{\infty}' )$ which is 
still a solution of (2.36)-(2.37), and which now has a free isometric 
$S^{1}\times S^{1}$ action. The second $S^{1}$ action arises from the 
unwrapping of the collapse of $V$ at infinity. However, for instance by 
[15, Thm.8.4], the only complete metrics of non-negative scalar 
curvature with such an action are flat. Hence again $\bar c = 0$, as 
claimed.

  Thus in all cases the existence of cusps leads to a contradiction, 
which proves the result.

{\endproof}

\begin{remark} \label{r 2.8}
A much simpler proof of Proposition 2.7 is possible if it is assumed 
that some (non-empty) component $\Omega_{o}$ of $\Omega$ is 
irreducible. Namely, there exists an exhaustion of $\Omega_{o}$ by 
compact sets $K_{j}$, such that $\Omega_{o} \setminus K_{j}$ is a graph 
manifold and $\partial K_{j}$ is a collection of tori. If any such 
torus is incompressible in $\Omega_{o}$, then as above, [15, Thm.8.4] 
implies that $\Omega_{o}$ is flat, and one has a contradiction as 
before. Thus, all such tori must be compressible. If now $\Omega_{o}$ 
is irreducible, then standard arguments in 3-manifold topology, (c.f. 
also the proof of Theorem 2.9 below), imply that $\Omega_{o}$ must be a 
solid torus $D^{2} \times S^{1}$. But this means that $\Omega_{o}$ is a 
graph manifold, and so the discussion in (2.14)-(2.15) holds. This 
means that $\Omega_{o}$ is either empty or flat, either of which is a 
contradiction.

  With some further topological arguments, this argument can be 
extended to the situation where it is assumed that $M$ is irreducible 
in place of $\Omega_{o}$.
\end{remark}

\medskip

{\bf \S 2.2.}
 In this subsection, we assemble the results of \S 2.1 to complete the 
proof of Theorem 0.2.

\medskip

\noindent
{\bf Proof of Theorem 0.2: $\sigma (M) < $ 0.}
 
 Suppose $M$ is a closed, oriented tame 3-manifold with
$$\sigma (M) <  0. $$
Again, let $\{(\Omega_{\varepsilon}, g_{\varepsilon})\}, \varepsilon  = 
\varepsilon_{i}$ be a sequence of minimizers of $I_{\varepsilon}^{~-}$ 
satisfying (2.1), so that as discussed in \S 2.1, a subsequence of 
$\{g_{\varepsilon}\}$ either converges, collapses, or forms cusps. By 
Proposition 2.3, the condition $\sigma (M) < $ 0 implies that no 
subsequence of $\{g_{\varepsilon}\}$ can collapse.

 If a subsequence of $\{g_{\varepsilon}\}$ converges, (modulo 
diffeomorphisms of $M$), then Proposition 2.1 shows it converges to a 
unit volume constant curvature metric $g_{o}$ on $M$, with scalar 
curvature $s_{g_{o}} = \sigma (M).$ In particular, (by rescaling 
$g_{o}$ by the factor $\sigma (M)/6), M$ admits a hyperbolic structure, 
and if $vol_{-1}M$ is the volume of $M$ in the hyperbolic metric, (of 
constant curvature -1), then
\begin{equation} \label{e2.38}
|\sigma (M)| = 6(vol_{-1}M)^{2/3}, 
\end{equation}
giving (0.9).

 If a subsequence of $\{g_{\varepsilon}\}$ forms cusps, then by 
Proposition 2.5, there is a maximal open set $\Omega\subset M,$ 
consisting of a bounded number of components, on which the subsequence 
converges, (modulo diffeomorphisms of $M$), to a limit metric $g_{o}.$ 
The metric $g_{o}$ is a complete, constant curvature metric satisfying 
\begin{equation} \label{e2.39}
s_{g_{o}} = \sigma (M), \ \  vol_{g_{o}}\Omega  = 1, 
\end{equation}
so that the pair $(\Omega , g_{o})$ realizes the Sigma constant of $M$ 
in this generalized sense. In particular, as in (2.38), this proves 
(0.8).

 The open manifold $\Omega $ has a finite number of components 
$\Omega_{i};$ each $\Omega_{i}$ has a finite number of ends, each 
diffeomorphic to $T^{2}\times {\mathbb R}^{+}.$ The domain $\Omega $ 
embeds as an open domain in $M$, and the complement $G = M\setminus 
\Omega $ has the structure of a manifold with boundary, with a finite 
number of components $G_{j}.$ Each $G_{j}$ is a graph manifold, with a 
finite number of boundary components, each diffeomorphic to $T^{2}.$ 
Thus the decomposition
\begin{equation} \label{e2.40}
M = \Omega  \cup  G 
\end{equation}
gives a decomposition of $M$ into hyperbolic and graph manifold 
regions, as in (0.7).

\medskip

 With the above understood, it remains to prove that each torus $T^{2}$ 
in a hyperbolic cusp $T^{2}\times {\mathbb R}^{+}\subset\Omega $ is 
incompressible in $M$. This is of course a crucial issue, since without 
it the geometric decomposition (2.40) may have no topological 
significance. In the same vein, one also needs to establish the 
uniqueness of the decomposition of $M$ into $\Omega $ and $M\setminus 
\Omega .$ Otherwise, some subsequences of $(\Omega_{\varepsilon}, 
g_{\varepsilon})$ may converge everywhere, while others may form cusps. 
We first prove the incompressibility of the tori; the uniqueness then 
follows easily after this.

 In the following, we assume for simplicity that the metric $g_{o}$ 
above is scaled to give a hyperbolic metric, i.e. a metric of constant 
sectional curvature $- 1,$ on $\Omega .$ 

\begin{theorem} \label{t 2.9.}
  Each torus $T^{2}$ in a hyperbolic cusp $T^{2}\times {\mathbb 
R}^{+}\subset\Omega $ is incompressible in $M$.
\end{theorem}

\noindent
{\bf Proof:}
  Let $T = T^{2}$ and suppose $D$ is a compressing disc in $M$, with 
$\partial D \subset  T$. By Dehn's Lemma, c.f. [18], we may assume that 
$D$ is embedded. Since $T$ is incompressible in $\Omega , D$ intersects 
$G$, and we may assume that $D \subset  G$. Let $G' $ be the component 
of $G$ containing $D$. Perturb $D$ to obtain another disjoint, isotopic 
disc $D' , \partial D'  \subset  T$, so that $\partial D\cup\partial D' 
$ bounds a small annulus $A$ in $T$. Note that $A\cup D\cup D' $ is the 
boundary of a small 3-ball in $M$, while the surface $(T\setminus 
A)\cup (D\cup D' )$ is an embedded 2-sphere $S^{2} \subset  G' .$

 If $G' $ is irreducible, or more precisely if the $S^{2}$ above bounds 
a 3-ball in $G' ,$ then it is clear that $T$ is the boundary of a solid 
torus $D^{2}\times S^{1}$, (c.f. [19, II.2.4]). In this case, the 
component $G' $ is $D^{2}\times S^{1}.$

 If the $S^{2}$ above does not bound a 3-ball in $G' ,$ then we may 
write $G' $ as
\begin{equation} \label{e2.41}
G'  = G_{1}\#... \#G_{q}, 
\end{equation}
where each $G_{i}$ is irreducible. By definition of connected sum, one 
factor in (2.41), say $G_{1}$, is then a solid torus $D^{2}\times 
S^{1}$ glued onto $T$. Thus, one has
\begin{equation} \label{e2.42}
G'  = (D^{2}\times S^{1})\#G'' , 
\end{equation}
where $G'' $ is a graph manifold. The factor $G'' $ may either be 
closed or have non-empty boundary; in the latter case, $\partial G'' $ 
again a union of tori, each contained in distinct ends of the 
hyperbolic part $\Omega $ of $M$.

 It follows that part of the original manifold $M$ is obtained by 
glueing on a solid torus, i.e. performing a Dehn surgery, to the torus 
$T$ in a cusp of $\Omega ,$ and possibly adding on other graph manifold 
components by connected sum. 

 We first prove Theorem 2.9 in the case $G''  = \emptyset  $ in (2.42), 
so that $G' $ is a solid torus. The general case is then easily 
obtained from this.

\medskip

  The idea here is to metrically glue on a solid torus explicitly and 
in such a way as to decrease $\mathcal{S}_{-}^{2}$ a definite amount below 
the value $|\sigma(M)|$, which of course gives a contradiction to 
(0.4). This turns out to be possible because the most natural metrics 
on a solid torus are of positive scalar curvature, and so contribute 
nothing to the value of $\mathcal{S}_{-}^{2}$. 

 To begin, choose an $S^{1}\times S^{1}$ product structure for $T^{2}$ 
in the hyperbolic cusp $T^{2}\times {\mathbb R}^{+}.$ This choice is of 
course not unique; it can be changed by an automorphism of $T^{2},$ 
i.e. element of $SL(2,{\mathbb Z} ).$ Thus, we may choose a basis 
$S^{1}\times S^{1},$ so that the Dehn surgery glues on a disc $D^{2}$ 
onto the first factor, while the second factor is left fixed.

 In this basis, the flat metric on $T^{2}$ may be written as
\begin{equation} \label{e2.43}
d_{1}^{2}d\theta_{1}^{2} + d_{2}^{2}d\theta_{2}^{2} + 2ad_{1}d_{2} 
d\theta_{1}d\theta_{2}, 
\end{equation}
where $a = \cos \alpha$, $\alpha$ the angle between the two $S^{1}$ 
factors, and $d_{1}, d_{2}$ are constants. The hyperbolic metric on the 
cusp is then given by
\begin{equation} \label{e2.44}
dt^{2} + e^{-2t}(d_{1}^{2}d\theta_{1}^{2} + d_{2}^{2}d\theta_{2}^{2} + 
2ad_{1}d_{2} d\theta_{1}d\theta_{2}). 
\end{equation}
Since we are only concerned with the end behavior, assume (2.44) to 
hold for $t \geq  t_{o},$ for $t_{o}$ a free large parameter, and 
assume the torus $T$ corresponds to the slice at $t_{o}.$ Now 
metrically glue on a disc $D^{2}$ to the first $S^{1}$ factor, with 
metric of the form,
\begin{equation} \label{e2.45}
dr^{2} + f_{1}^{2}d\theta_{1}^{2} + f_{2}^{2}d\theta_{2}^{2} + 
2af_{1}f_{2} d\theta_{1}d\theta_{2}. 
\end{equation}
Here, we choose functions $f_{1}, f_{2}$ depending only on $r$, where 
$r\in [0, \frac{\pi}{2}].$ For the moment, it is required that the metric 
is piecewise smooth (at least $C^{2}),$ and is $C^{1}$ at the seam 
where the disc is glued to $T$. Comparing the forms (2.44) and (2.45), 
this amounts to the requirement that $f_{1}, f_{2}$ satisfy the 
boundary conditions 
\begin{equation} \label{e2.46}
f_{1}(0) = 0, \ f_{1}'(0)(1-a^{2})^{1/2} = 1, \  f_{1}(\frac{\pi}{2}) = 
d_{1}e^{-t_{o}}, \ f_{1}' (\frac{\pi}{2}) = d_{1}e^{-t_{o}}, 
\end{equation}
\begin{equation} \label{e2.47}
f_{2}(0) = a >  0, \ f_{2}' (0) = 0, \ f_{2}(\frac{\pi}{2}) = 
d_{2}e^{-t_{o}}, \  f_{2}' (\frac{\pi}{2}) = d_{2}e^{-t_{o}}, 
\end{equation}
 It is simpler for computations to follow to orthogonalize the basis 
$\theta_{1}, \theta_{2}.$ Thus, set
$$d\bar \theta_{2} = d\theta_{2} + a \frac{f_{1}}{f_{2}}d\theta_{1}. $$
A simple substitution shows that the metric (2.45) may be rewritten as
\begin{equation} \label{e2.48}
dr^{2} + f_{1}^{2}(1 -  a^{2})d\theta_{1}^{2} + f_{2}^{2}d\bar 
\theta_{2}^{2}. 
\end{equation}

 In the following, specific warping functions $f_{1}$ and $f_{2}$ 
are chosen to satisfy (2.46)-(2.47), giving a specific metric glueing 
of $D^{2}$ onto $S^{1}.$ We are not interested in any optimal choices, 
and so just choose specific forms for $f_{1}$ and $f_{2}$ that suffice 
for the needs of the argument.

 As a first approximation to the glueing, set
\begin{equation} \label{e2.49}
f_{1}(r) = c_{1}\tan \frac{r}{2}, \ \   f_{2}(r) = c_{2}e^{-\cos r}, 
\end{equation}
where $c_{1} = d_{1}e^{-t_{o}}, c_{2} = d_{2}e^{-t_{o}}.$ A simple 
calculation shows that this metric is $C^{1}$ at $t = t_{o},$ i.e. 
where $r = \frac{\pi}{2}.$ Further, $f_{2}$ satisfies (2.47) and 
$f_{1}$ satisfies $f_{1}(0) =$ 0, but $f_{1}' (0)\cdot  (1-a^{2})^{1/2} 
= \frac{1}{2}c_{1}(1-a^{2})^{1/2} \neq $ 1. Note that $c_{1} << $ 1 for 
$t_{o}$ sufficiently large. In other words, the glueing of $D^{2}\times 
S^{1}$ is a $C^{o}$ glueing, but not $C^{1}.$ The glued solid torus is 
metrically a cone manifold, with cone angle $\pi c_{1}(1-a^{2})^{1/2}$ 
along the core geodesic $\gamma  = \{r=0\}.$ In particular, there is a 
concentration of positive curvature along $\gamma ,$ for $t_{o}$ large. 
This cone singularity along the core curve will be smoothed later.

 For metrics of the form (2.48), it is easily computed that the scalar 
curvature is given by
\begin{equation} \label{e2.50}
\frac{1}{2}s = -  \frac{f_{1}''}{f_{1}} -  \frac{f_{2}''}{f_{2}} -  
\frac{f_{1}' f_{2}'}{f_{1}f_{2}}. 
\end{equation}
For the choices (2.49), one has
$$\frac{f_{1}''}{f_{1}} = \frac{1}{2}(\cos (r/2))^{-2}, \ \  
\frac{f_{1}'}{f_{1}} = (\sin r)^{-1} ,   $$
$$\frac{f_{2}''}{f_{2}} = \cos r + \sin^{2} r, \ \  \frac{f_{2}'}{f_{2}} 
= \sin r , $$
which gives
$$\frac{1}{2}s = - 1 -  \cos r -  \sin^{2}r -  \frac{1}{2}\cos^{-2}(r/2). 
$$

 An exercise in calculus gives the bound
\begin{equation} \label{e2.51}
s \geq  - 6, 
\end{equation}
on the full solid torus $r^{-1}[0, \frac{\pi}{2}],$ with $s >  - 6$ on 
$r^{-1}[0, \frac{\pi}{2}).$ It is worth pointing out that this metric 
does not have sectional curvature $K \geq  - 1$ everywhere.

 Next, we estimate the volume of the glueing, compared with the volume 
of the hyperbolic cusp. First, the volume of the hyperbolic cusp, cut 
off at $t = t_{o}$ is given by
\begin{equation} \label{e2.52}
V_{C}(t_{o}, \infty ) = d_{1}d_{2}\int_{t_{o}}^{\infty}e^{-2t}dt = 
\frac{1}{2}d_{1}d_{2}e^{-2t_{o}}. 
\end{equation}
On the other hand the volume of the solid torus glued in at $t = t_{o}$ 
is given by
$$V_{D}(0, \frac{\pi}{2}) = \int_{0}^{\pi /2}f_{1}f_{2} dr = 
d_{1}d_{2}e^{-2t_{o}}\int_{0}^{\pi /2}e^{-\cos r}\tan (r/2) dr,$$
and, for instance, a numerical evaluation shows that
\begin{equation} \label{e2.53}
\int_{0}^{\pi /2}e^{-\cos r}\tan (r/2)dr < .464 <  \frac{1}{2}. 
\end{equation}
Thus, the volume of the glued solid torus is less than the volume of 
the hyperbolic cusp.

 We now turn to the smoothing of the cone singularity along the core 
curve. Since there is a concentration of positive curvature at the core 
curve, we will smooth it to have positive scalar curvature nearby. 

 In detail, the function $f_{2}$ satisfies the boundary condition 
(2.47) at $r =$ 0, and so need not be changed. Now recall that 
$f_{1}(0) =$ 0, while 
\begin{equation} \label{e2.54}
f_{1}' (0) = c_{1}/2 <<  (1- a^{2})^{-1/2}.  
\end{equation}
The latter estimate follows since we are free to choose $c_{1}$ 
arbitrarily small, by going sufficiently far down the hyperbolic cusp. 
Given then a fixed small choice of $c_{1},$ choose $r_{o} <<  c_{1}.$ 
One may then bend $f_{1}$ on the interval $[r_{o}/2, r_{o}],$ to a new 
function $\bar f_{1}$ satisfying the boundary conditions (2.46) at 
$r_{o}/2$ in place of 0, i.e.
\begin{equation} \label{e2.55}
\bar f_{1}(r_{o}/2) = 0, \bar f_{1}' (r_{o}/2) = (1-a^{2})^{-1/2} ,
\end{equation}
while $\bar f_{1}$ agrees $C^{1}$ with $f_{1}$ at the value $r_{o}.$ 
This bending then obviously has the property that
\begin{equation} \label{e2.56}
\frac{\bar f_{1}''}{\bar f_{1}} <<  0. 
\end{equation}
One may then readjust $f_{2}$ so that it satisfies the boundary 
conditions (2.47) at $r_{o}/2$ in place of 0. This requires only a 
small $C^{2}$ perturbation of $f_{2},$ so that its contribution to the 
curvature only changes slightly. Together with (2.50) and (2.51), 
(2.56) implies that the scalar curvature of the resulting metric is 
(very) positive in the small solid torus given by the region 
$r^{-1}[r_{o}/2, r_{o}].$ Obviously the change to the volume by this 
bending can be made arbitrarily small.

 It follows that the metric on $D^{2}\times S^{1} \cup  T^{2}\times 
{\mathbb R}^{+},$ given by (2.45) on $[r_{o}/2, \frac{\pi}{2}],$ for the 
above choices for $f_{1}$ and $f_{2},$ and by (2.44) on $[t_{o}, 
t_{o}+1]$ is a $C^{1}$ smooth metric, which is $C^{\infty}$ off the 
seams at $t_{o} = \{r = \frac{\pi}{2}\},$ and $r_{o}.$ This metric has 
smaller volume than the volume of the hyperbolic cusp and has scalar 
curvature $s \geq  - 6$ everywhere, except at the seams, where $s$ is 
not defined. Note that $s$ is a piecewise smooth function, with only 
jump discontinuities at the seams $t_{o}, r_{o},$ and that $s >  - 6$ 
for $r <  \frac{\pi}{2}.$ We may take a smooth approximation $\bar g$ 
to this metric by smoothing the warping functions. The scalar curvature 
$\bar s$ of $\bar g$ then interpolates the values of $s$, i.e. smooths 
the jump discontinuity. Clearly, one may thus choose a smoothing so 
that $\bar s \geq  - 6$ everywhere.

 The metric $\bar g$ is a complete smooth metric on a domain $\bar M 
\subset  M$ satisfying
\begin{equation} \label{e2.57}
s_{\bar g} \geq  - 6, \ \ {\rm and} \ \  vol_{\bar g}\bar M <  
vol_{g_{o}}\Omega . 
\end{equation}
In case $\bar M \neq  M$, so that there are graph manifold components 
of $M$ remaining in $M\setminus \bar M,$ one may use [2, Thm 5.7] as 
above in the proof of (2.24), to obtain a metric $\widetilde g$ on $M$ 
satisfying
\begin{equation} \label{e2.58}
s_{\widetilde g} \geq  - 6 \ {\rm on} \ M\setminus K, \ |s_{\widetilde 
g}| \leq  C, \ vol_{\widetilde g}M\setminus K \leq  \delta , 
\end{equation}
where $K$ is an arbitrarily prescribed compact domain in $\bar M, C$ is 
independent of $K$, and $\delta $ may be made arbitrarily small by 
choosing $K$ sufficiently large.

 Now from (2.39), $\sigma (M)$ is realized by the complete hyperbolic 
metric $g_{o}$ on $\Omega  \subset  M$. It follows from the last 
estimates in (2.57) and (2.58) that, for $\delta $ sufficiently small,
\begin{equation} \label{e2.59}
vol_{\widetilde g}M <  vol_{g_{o}}\Omega , 
\end{equation}
while by the remaining estimates in (2.57)-(2.58),
\begin{equation} \label{e2.60}
\int_{M}(s_{\widetilde g}^{-})^{2}dV_{\widetilde g} <  
\int_{\Omega}s_{g_{o}}^{2}dV_{g} = 
\int_{\Omega}(s_{g_{o}}^{-})^{2}dV_{g}. 
\end{equation}
This implies that 
\begin{equation} \label{e2.61}
\mathcal{S}  _{-}^{2}(\widetilde g) <  |\sigma (M)|, 
\end{equation}
which is of course impossible by (0.4). 

 This completes the proof of Theorem 2.9 in case $G' $ is a solid 
torus. Now suppose that the component $G' $ of $G$ is of the form
\begin{equation} \label{e2.62}
G'  = (D^{2}\times S^{1})\#G'' , 
\end{equation}
with $G'' $ a graph manifold. 

 Suppose first that $G'' $ is closed. Then one has the decomposition
$$M = M' \#G'' , \ \ {\rm where} \ \ M'  = (M\setminus G' 
)\cup_{T^{2}}(D^{2}\times S^{1}), $$
and $M' $ is closed. As above, construct the metric $\widetilde g$ 
on the closed manifold $M' ,$ where $M' $ is obtained by performing 
Dehn surgery on the given $T^{2}.$ As in (2.61), this gives
\begin{equation} \label{e2.63}
\mathcal{S}  _{-}^{2}(\widetilde g) <  |\sigma (M')|. 
\end{equation}

 If however $G'' $ is not closed, then one may write
$$G''  = \#_{1}^{k}(D^{2}\times S^{1})\#G''' , $$
where $G''' $ is a closed graph manifold and $G'' $ has $k \geq $ 1 
toral boundary components. Each such torus is (isotopic to) a torus in 
a hyperbolic cusp in $\Omega$. Hence,
$$M = M' \#G''' , $$
where $M' $ is the closed manifold obtained from $M\setminus G' $ by 
glueing on solid tori to its boundary components. Thus, in the same 
manner as before, construct the metric $\widetilde g$ on $M' $ 
satisfying (2.63).

 Now since $G'' $ or $G''' $ is a closed graph manifold, it admits 
metrics $h = h_{\delta}$ with arbitrarily small volume, and arbitrarily 
small $L^{2}$ norm of scalar curvature. In particular, $\mathcal{S}  
_{-}^{2}(h_{\delta}) \leq  \delta ,$ for any prescribed $\delta  > $ 0. 
As discussed in [2,\S 7], [5,\S 6.1], or also in Case (i) of 
Theorem 2.14 below, one may then form a metric $g^{*}$ on the sum $M' 
\#G'' ,$ (or $M' \#G''' ),$ agreeing with $\widetilde g$ outside a very 
small ball in $M' ,$ and with $h_{\delta}$ outside a very small ball in 
$G'' ,$ and so that on the neck $S^{2}\times I, g^{*}$ has uniformly 
bounded scalar curvature and arbitrarily small volume. In fact one may 
choose $g^{*}|_{S^{2}\times I}$ to be a truncation of the isometrically 
doubled Schwarzschild metric $g_{S}$ in (1.23) with mass $m <<  \delta 
.$ Thus we have constructed metrics $g^{*}$ on $M$ satisfying
\begin{equation} \label{e2.64}
\mathcal{S}  _{-}^{2}(g^{*}) \leq  \mathcal{S}  _{-}^{2}(\widetilde g) + 
\delta , 
\end{equation}
for any given $\delta  > $ 0; of course $g^{*}$ depends on $\delta$. 
Choosing $\delta $ sufficiently small, from (2.63) one again obtains
$$\mathcal{S}  _{-}^{2}(g^{*}) <  |\sigma (M)|, $$
which is impossible. It follows that each torus $T_{i}$ in the 
hyperbolic cusps of $\Omega $ is incompressible in $M$.

{\endproof}

 This result is similar in spirit, (although the proof is quite 
different), to Thurston's cusp closing theorem, c.f. [36].

\medskip

 Theorem 2.9 and the preceding work imply the existence of the 
decomposition (0.7) of $M$. In particular, the tori $T_{i}^{2}$ in the 
hyperbolic cusps give a partial torus decomposition of $M$, in the 
sense of Jaco-Shalen-Johannson, [19], [21]. We discuss the full torus 
decomposition of $M$, corresponding to the further decomposition of $G$ 
into Seifert fibered spaces $S_{k},$ in \S 2.3 below. The manifold $G = 
\cup S_{k}$ is called the {\it  characteristic variety}  of $M$, c.f. 
[19], [21]; Theorem 0.2 gives a new (geometric) proof of its existence 
for tame 3-manifolds.

 Finally, we discuss the issue of uniqueness of this decomposition. 
Thus, given one decomposition
\begin{equation} \label{e2.65}
M = \Omega\cup G, 
\end{equation}
where the union is along incompressible tori $\{T_{i}^{2}\},$ let $T' $ 
be any other incompressible torus. We claim that $T' $ can be isotoped 
into $G$; my thanks to Yair Minsky for assistance with the argument 
below. First, since $\Omega $ is hyperbolic, any incompressible torus 
$T' $ embedded in $\Omega $ is boundary parallel, and thus may be 
isotoped into $G$. If $T' $ intersects $\Omega $ and $G$, then there is 
a torus $T\in\{T_{i}^{2}\}$ such that $T\cap T'  \neq  \emptyset  .$ 
Now $T\cap T' $ is a collection of embedded essential circles, bounding 
annuli $\{A_{i}\}$ in $T$. It follows that one may form a new torus 
$\bar T$ by matching the annuli $T'\cap\Omega $ with the annuli 
$\{A_{i}\}.$ The torus $\bar T$ remains incompressible, and lies in the 
hyperbolic manifold $\Omega .$ Thus, it may be isotoped into $G$, which 
induces an isotopy of $T' $ into $G$.

 Hence, if $\{T_{j}'\}$ is another torus decomposition of $M$ into 
$\Omega' $ and $G' ,$ then $\{T_{j}'\} \subset  G$. Further, one may 
assume that the collection $\{T_{j}'\}$ is disjoint from the collection 
$\{T_{i}\}.$ Now if $\{T_{j}'\}$ is not isotopic to $\{T_{i}\},$ then 
it follows that one of the tori $T_{k}\in\{T_{i}\}$ is contained in one 
of the hyperbolic manifolds $\Omega_{k}' $ and is not isotopic to a 
boundary torus in a cusp of $\Omega_{k}' .$ Since $T_{k}$ is 
incompressible, this is impossible. Thus, the collections $\{T_{i}\}$ 
and $\{T_{j}'\}$ of tori are isotopic.

 This completes the proof of Theorem 0.2 when $\sigma (M) < $ 0.

{\endproof}

\noindent
{\bf Proof of Theorem 0.2: $\sigma (M) =$ 0.}

  Suppose $M$ is a closed, oriented tame 3-manifold with
$$\sigma (M) = 0. $$
Consider a sequence $\{(\Omega_{\varepsilon}, g_{\varepsilon})\}$ of 
minimizers of $I_{\varepsilon}^{~-}, \varepsilon  = \varepsilon_{i},$ 
satisfying (2.1). Thus, as discussed in \S 2.1, a subsequence of 
$\{(\Omega_{\varepsilon}, g_{\varepsilon})\}$ either converges, 
collapses, or forms cusps.

 If $\{g_{\varepsilon}\}$ converges, (modulo diffeomorphisms of $M$), 
then by Proposition 2.1, it converges everywhere on $M$ to a smooth 
limit metric $g_{o}$ of constant curvature, and scalar curvature 
$s_{g_{o}} = \sigma (M) =$ 0. Thus, $g_{o}$ is a flat metric on $M$, so 
that $(M, g_{o})$ is a flat 3-manifold, and in particular a graph 
manifold with infinite $\pi_{1}.$ In this case, the metric $g_{o}$ on 
$M$ realizes $\sigma (M).$

 If $\{g_{\varepsilon}\}$ collapses, then as in Proposition 2.3, it 
collapses along a sequence of F-structures on $M$. In particular, $M$ 
is a graph manifold. In this case, if $M$ is not a flat 3-manifold, no 
smooth metric on $M$ realizes $\sigma (M).$

 We claim that $|\pi_{1}(M)| = \infty .$ For it is standard, c.f. 
[26, \S5], that a closed oriented graph manifold $M$ with $\pi_{1}(M)$ 
finite must be a spherical space form $S^{3}/\Gamma ,$ where $\Gamma $ 
is a finite subgroup of $SO(4)$. Such manifolds have Yamabe metrics of 
positive scalar curvature, so that $\sigma (M) > $ 0. Since $\sigma (M) 
=$ 0, this gives the claim.

 By Proposition 2.7, the sequence $\{g_{\varepsilon}\}$ cannot form 
cusps. This completes the proof of Theorem 0.2 in case $\sigma (M) =$ 0.

{\endproof}

\begin{remark} \label{r 2.10} {\bf (i).} 
  Referring to Remarks 2.2 and 2.6, the uniqueness results above prove 
that if $M$ is tame, then minimizing pairs $(\Omega_{\varepsilon}, 
g_{\varepsilon})$ are unique, and independent of $\varepsilon ,$ for 
all $\varepsilon  > $ 0, (with the exception of flat manifolds). In 
particular, the decomposition (0.7) is the same as the decomposition 
given in Theorem 1.1.

  Note also that the proof of Theorem 0.2 does not require that $M$ is 
irreducible.

{\bf (ii).}
 For completeness, we make a few remarks on the Sigma constant of graph 
manifolds. By Remark 2.4, an arbitrary closed oriented graph manifold 
$G$ necessarily satisfies
$$\sigma (G) \geq  0. $$

 Standard 3-manifold topology, c.f. [26, \S 6] or [38], implies that an 
irreducible graph manifold of infinite $\pi_{1}$ necessarily has a 
${\mathbb Z}\oplus{\mathbb Z}  \subset\pi_{1}.$ From results of Schoen-Yau, 
[31], this implies that $M$ has no metric of positive scalar 
curvature, and thus in particular, no Yamabe metric of positive scalar 
curvature. Thus, we have the implication
\begin{equation} \label{e2.66}
|\pi_{1}(N)| = \infty  \Rightarrow  \sigma (N) = 0, 
\end{equation}
for irreducible graph manifolds $N$. The opposite implication
$$\sigma (M) = 0 \Rightarrow  |\pi_{1}(M)| = \infty , $$
for arbitrary closed graph manifolds follows from the proof of Theorem 
0.2 above.

\end{remark}

\begin{remark} \label{r 2.11.}
  As noted in the Introduction, the geometrization of closed oriented 
irreducible 3-manifolds which are tame and of infinite $\pi_{1}$ is now 
a simple consequence of Theorem 0.2. Namely, as noted following (0.1), 
such a manifold $M$ must be a $K(\pi ,1),$ and hence by (0.6), $\sigma 
(M) \leq $ 0.

 If in addition $M$ is atoroidal, then by Remark 2.10, $M$ cannot be an 
irreducible graph manifold, and so the $\sigma (M) =$ 0 case of Theorem 
0.2 implies that $\sigma (M) < $ 0. Hence Theorem 0.2 implies that $M$ 
is hyperbolic.

 If $M$ contains essential tori, then Theorem 0.2 implies that $M$ is a 
union along essential tori of complete hyperbolic manifolds and graph 
manifolds. The hyperbolic manifolds are of course geometric, and graph 
manifolds are also geometric, since they are unions along tori of 
Seifert fibered spaces; the geometrization of the graph manifolds will 
be discussed next in \S 2.3.

\end{remark}

{\bf \S 2.3.}
 We complete this section with an analysis of how the torus 
decomposition and the geometrization of graph manifolds can be obtained 
from the near-limiting behavior of suitable minimizing sequences for 
$\mathcal{S}_{-}^{2}$ or $I_{\varepsilon}^{~-}$. Thus we describe in 
particular how the other Seifert fibered geometries arise in this 
context. On the other hand, this section is not logically necessary for 
any later work.

 The graph manifold structures in Theorem 0.2 in either case $\sigma 
(M) < $ 0 or $\sigma (M) =$ 0 are obtained from the near-limiting 
behavior of a minimizing sequence $\{g_{i}\}$ for 
$I_{\varepsilon}^{~-}, \varepsilon  >$ 0; as noted in Remark 2.10, the 
parameter $\varepsilon $ no longer plays any role and the conclusions 
of Theorem 0.2 may be obtained for any choice of $\varepsilon  > $ 0. 
The results are independent of $\varepsilon .$

 In the situation $\sigma (M) < $ 0, on the hyperbolic domain $H = \cup 
H_{j},$ the metrics $g_{i}$ on $M$ converge to the unique limit metric 
$g_{o}$ on $H$. On the graph manifold domain $G = \cup G_{k},$ the 
sequence $\{g_{i}\}$ collapses with bounded curvature (in $L^{2}).$ Of 
course, there are many possible choices for the minimizing sequence 
$\{g_{i}\},$ and so apriori there may be many possible ways to collapse 
$G$.

 Thus, throughout \S 2.3, let $G$ be an oriented graph manifold, either 
closed and satisfying $\sigma (G) =$ 0, or with boundary consisting of 
a finite number of incompressible tori corresponding to the situations 
$\sigma (M) =$ 0 or $\sigma (M) < $ 0 in Conjectures II and I 
respectively. In the latter case, it follows from a well-known result 
of Gromov-Lawson [15, Thm.8.4], that $G$ admits no complete metrics of 
non-negative scalar curvature which are not flat. Since such a $G$ does 
admit complete metrics with $\mathcal{S}  _{-}^{2}$ arbitrarily small, we 
will say that $\sigma (G) =$ 0 in this case also. Note that we do not 
assume that $G$ is irreducible. 

\smallskip

 We begin with a general discussion of the decomposition of graph 
manifolds obtained from the geometry of a collapsing sequence of 
metrics. Since $G$ admits polarized F-structures, $\inf I_{\varepsilon}^{~-} 
= \inf \mathcal{S}  _{-}^{2} = 0$, and there are 
minimizing sequences $\{g_{i}\}$ for $I_{\varepsilon}^{~-}$ which have 
uniformly bounded curvature in $L^{\infty}$ and which volume collapse 
$G$ along a sequence $\mathcal{F}_{i}$ of polarized F-structures. (The 
uniqueness of this collapse will be discussed below). From the 
definition of polarized F-structure, it follows that $G$ is partitioned 
into a collection of regions, possibly depending on $i$, on which the 
orbits of $\mathcal{F}_{i}$ are circles and tori. Let $S_{j} = S_{j}(i)$ 
be the components of the regions with $S^{1}$ fibers, and $L_{k} = 
L_{k}(i)$ the components of regions with $T^{2}$ fibers; thus $\{S_{j}, 
g_{i}\}$ collapses along $S^{1}$ orbits, while $\{(L_{k}, g_{i})\}$ 
collapses along $T^{2}$ orbits. We assume w.l.o.g. that all domains 
$S_{j}$ and $L_{k}$ have only toral boundary components, invariant 
under the group actions. Assuming the collection $\{S_{j}\}$ is 
non-empty, the components $L_{k}$ are considered as regions where the 
$S_{j}$ components are glued together by toral diffeomorphisms. 

 Although this decomposition $D$ of $G$ into these types of regions, 
i.e. the placement of the tori separating the $S$ and $L$ factors, is 
not precisely determined by the geometry of $\{g_{i}\},$ it is a 
consequence of the fact that the F-structure is polarized that the 
transition between any two distinct components of this division 
necesssarily takes place over domains of larger and larger diameter as 
$i \rightarrow  \infty ,$ c.f. [9]. Thus, as $i \rightarrow  \infty ,$ 
the diameter of each component $S_{j}$ or $L_{k}$ is forced to go to 
$\infty $ and the metric $g_{i}$ restricted to any such component 
becomes more and more complete.

\smallskip

 The Seifert fibered components $S = S_{j}$ are $S^{1}$ fibrations over 
a surface $\Sigma  = \Sigma_{j}.$ Assuming that $\Sigma $ is not 
closed, $\Sigma $ is thus either an open hyperbolic surface, an annulus 
$S^{1}\times I$, a M\"obius band, or a disc $D^{2}$. The latter three 
cases are somewhat exceptional and we make some further remarks on them.

If $\Sigma $ is an annulus, so that $S_{j} = T^{2}\times I,$ one has 
two possibilities. If the $S^{1}$ structure on $S_{j}$ is compatible 
with the $S^{1}$ structure on some neighbor $S_{k}$, (near the 
boundary), we just extend the Seifert fibered structure on $S_{k}$ to 
include $S_{j}$ and drop $S_{j}$ from the list. If the $S^{1}$ 
structure on $S_{j}$ is not compatible with either of its neighbors in 
this sense, then the collapse theory [8] implies that $S_{j}$ must be 
an $L$ factor, i.e. $S_{j}$ is collapsed along both $S^{1}$ factors. 
Thus in either case, there are no such $S$ factors.

 If $\Sigma $ is a M\"obius band, then $S$ is the oriented $S^{1}$ 
bundle over $\Sigma .$ Note that $\partial S$ is connected, and is an 
incompressible torus in $S$. This manifold admits a complete flat 
geometry.

 The situation where $\Sigma  = D^{2},$ so that $S$ is a solid torus 
$D^{2}\times S^{1}$ is the most important, in that this is the only 
Seifert fibered space with compressible boundary. Although $D^{2}$ 
admits complete hyperbolic and complete flat metrics, it is the 
spherical metric on $D^{2}$ embedded as a hemisphere in $S^{2}$ which 
is geometrically the most natural, c.f. the proof of Theorem 2.9. Such 
a metric is of course incomplete.

\medskip

 In the next result we show that any decomposition $D$ of $G$ 
containing solid tori components in the Seifert fibered factors above 
may be altered to a new decomposition without such components.

\begin{lemma} \label{l 2.12.}
  Let $G$ be either a closed oriented graph manifold with $\sigma (G) 
=$ 0, or a compact oriented graph manifold with incompressible 
boundary, so that $\partial G$ is a finite union of tori. Let
\begin{equation} \label{e2.67}
G = G_{1}\#... \#G_{k} 
\end{equation}
be the sphere decomposition of $G$ into prime factors. 

 Then each $G_{p}$ is a prime graph manifold with $\sigma (G_{p}) \geq 
$ 0, and $\sigma (G_{j}) =$ 0, for some j. Further, there exists a 
decomposition $D_{p}$ of each $G_{p}$ such that no Seifert fibered 
component $S$ in $G_{p}$ above is a solid torus.
\end{lemma}

\noindent
{\bf Proof:}
 Each $G_{p}$ is a graph manifold, either closed or with incompressible 
torus boundary. Hence by Remark 2.4 and the discussion above, $\sigma 
(G_{p}) \geq $ 0 for each $p$. If $\sigma (G_{p}) > $ 0, then $G_{p}$ 
is either $S^{2}\times S^{1}$ or a spherical space form $S^{3}/\Gamma 
.$ If all $G_{p}$ satisfy $\sigma (G_{p}) > $ 0, then $\sigma (G) > $ 
0, since positive scalar curvature is preserved under connected sums, 
c.f. [14], [32]; thus $\sigma(G_{j}) = 0$, for some $j$.

 Let $D$ be any decomposition of $G$ as above into $S_{j}$ and $L_{k}$ 
factors. Suppose some $S = S_{j}$ factor is a solid torus $D^{2}\times 
S^{1}$. The factor $S$ must then be glued, (via an $L$ component), onto 
a boundary component $T$ of a Seifert fibered space $S'  = S_{j'}.$ We 
may assume that $T$ is incompressible in $S' ,$ for otherwise $G$ is a 
union of two solid tori, and hence a lens space or $S^{2}\times S^{1},$ 
either of which have positive Sigma constant, contradicting $\sigma(G) 
= 0$.

 There are exactly two possibilities for the resulting topology of 
$S'\cup_{T^{2}}(D^{2}\times S^{1}),$ as discussed by Waldhausen, c.f. 
[39].

 (i). First, if the fiber of $S' $ is not glued to the meridian, i.e. 
to $D^{2},$ then the resulting manifold is again a Seifert fibered 
space $\bar S,$ possibly with an exceptional fiber at the core of the 
solid torus, c.f. [39, \S 10]. Hence, in this case one may just alter 
the decomposition $D$ by eliminating the $L$ factor between $S' $ and 
$D^{2}\times S^{1}$ and enlarging $S' $ to the Seifert fibered space 
$S'\cup_{T^{2}}(D^{2}\times S^{1}).$ (Waldhausen's definition of graph 
manifold excludes this operation, but we do not do so here).

 (ii). On the other hand, if the fiber is glued to the meridian, then 
the resulting manifold is reducible, c.f. [39, p.90]. This implies that 
$G$ has a non-trivial sphere decomposition, $G = G^{1}\#G^{2}.$ The 
essential 2-sphere is formed as in the proof of Theorem 2.9. The 
factors $G^{1}$ and $G^{2}$ may be disconnected by glueing in a 3-ball 
on each side. The original decomposition $D$ is then altered or split 
to give a decomposition $D^{i}$ of each $G^{i},$ c.f. [39, p.90]. By 
the Kneser finiteness theorem [22], or by [39, p.92], there are only 
finitely many isotopy classes of essential 2-spheres in $G$, so that 
this process terminates after a finite number of iterations.

 The result of this process is then the decomposition (2.67) with each 
$G_{p}$ prime and a corresponding decomposition $D_{p}$ of $G_{p}$ 
without any solid torus components.

{\endproof}

 Note that the proof above shows that any initial graph manifold 
decomposition of $G$ may be altered in a canonical way to a 
decomposition $D$ satisfying the conclusions of Lemma 2.12. In 
particular, this process carries out the sphere decomposition of $G$. 
The union of the decompositions $D_{p}$ of each prime factor $G_{p}$ in 
(2.67) does not extend to a graph manifold decomposition of $G$. This 
corresponds to the fact that the necks $S^{2}\times I$ about essential 
2-spheres in $G$ are no longer decomposed into $S$ and $L$ factors. 

 Next we discuss the uniquenesss of the decomposition of each factor 
$G_{p}$ in (2.67).

\begin{lemma} \label{l 2.13.}
  Let $G$ be a closed oriented graph manifold with $\sigma (G) =$ 0, or 
a compact oriented graph manifold with incompressible toral boundary 
components, (and hence also $\sigma (G) =$ 0). Suppose $G$ is 
irreducible, and not a closed flat, Nil or Sol manifold, or the 
oriented $S^{1}$ bundle over a M\"obius band.

 Then the decomposition $D$ of $G$ as in Lemma 2.12, into Seifert 
fibered and toral factors without solid torus components, is unique up 
to isotopy of G. The tori in the $L_{k}$ factors give the JSJ torus 
decomposition of G. Further, the decomposition is injective in the 
sense that $\pi_{1}(S_{j})$ and $\pi_{1}(L_{k})$ inject in $\pi_{1}(G).$
\end{lemma}

\noindent
{\bf Proof:}
 The existence of the decomposition follows from Lemma 2.12 and its 
uniqueness up to isotopy is proved by Waldhausen [39, \S 8], c.f. also 
[26, p.132] and [28, Lemma 5.5]. The four exceptions arise from the fact 
that in these cases one may either choose $D$ to be empty, (the 
manifolds are geometric), or may one choose $D$ to be non-empty in 
inequivalent ways. The injectivity statement is proved in [28, Thm. 
4.2]. The fact that the tori in the $L_{k}$ factors give the JSJ 
decomposition is then immediate.
{\endproof}

 We recall that the prime factors $G_{p}$ of the sphere decomposition 
(2.67) with $\sigma (G_{p}) > $ 0 are either $S^{2}\times S^{1}$ or 
space forms $S^{3}/\Gamma .$

\medskip

 Given these results on the structure of graph manifolds, we now return 
to the geometric behavior of suitable minimizing sequences. Let $G$ be 
as in Lemma 2.12. We will say that a minimizing sequence $\{g_{i}\}$ 
for $\mathcal{S}  _{-}^{2}$ on $G$ {\it  performs}  the geometrization of 
$G$ if the following conditions hold:

 (i). Let $\{S_{q}^{2}\}$ be a fixed collection of essential 2-spheres 
defining the sphere decomposition (2.67) of $G$. The sequence 
$\{g_{i}\}$ crushes each $S^{2}\in\{S_{q}^{2}\}$ to a point $x = x_{q} 
\in  G$, in the sense that
\begin{equation} \label{e2.68}
diam_{g_{i}}S^{2} \rightarrow  0 \ \ {\rm and} \ \ \rho_{i}(x) 
\rightarrow  0,  \ \ {\rm as} \ \ i \rightarrow  \infty , 
\end{equation}
where $\rho $ is the $L^{2}$ curvature radius. Further, the blow-ups of 
the metrics $\{g_{i}\}$ by $\rho_{i}(x)^{-2}$ converge to the complete 
doubled Schwarzschild metric (1.23), with mass $m \sim $ 1.

 (ii).  The metrics $\{g_{i}\}$ collapse each irreducible factor 
$G_{p}$ with $\sigma (G_{p}) =$ 0, (except possibly factors which are 
closed flat manifolds), along a polarized F-structure, (unique in the 
non-exceptional cases), determined by Lemma 2.13, on a scale however 
much larger than the crushing of the essential 2-spheres, i.e.
\begin{equation} \label{e2.69}
diam_{g_{i}}\mathcal{O} (y) >>  diam_{g_{i}}S^{2},  
\end{equation}
for any $y$, where $\mathcal{O} $ is the orbit of the F-structure through 
$y$. 

 (iii). After passing to suitable covers of each $S$ component of 
$G_{p}$ with $\sigma (G_{p}) =$ 0, the metrics $g_{i}$ converge to one 
of the six Seifert fibered geometries on $S$, i.e. $SL(2,{\mathbb R})$, 
Nil, $H^{2}\times {\mathbb R}$, ${\mathbb R}^{3}$, $S^{2}\times {\mathbb R}$, 
$S^{3},$ or suitable covers of $G_{p}$ converge to the Sol geometry. 
On each factor $G_{p}$ with $\sigma (G_{p}) > 0$, the metrics $g_{i}$ 
converge to the constant curvature metric of curvature +1 on $S^{3}/ 
\Gamma$ or the product metric on $S^{2}\times S^{1}$ of the form 
$S^{2}(1)\times S^{1}(1)$.

\smallskip

 Note that in general, such a sequence $\{g_{i}\}$ is not tame, in that 
the $L^{2}$ norm of the curvature $z_{g_{i}}$ diverges to infinity, due 
to the crushing of the essential 2-spheres. On the other hand, any 
graph manifold $G$ as above does admit tame minimizing sequences which 
volume collapse $G$ with bounded curvature, as discussed in the 
beginning of \S 2.3. In particular, on reducible 3-manifolds satisfying 
the conclusions of Conjectures I and II, a tame minimizing sequence 
does not perform the geometrization of $M$.

 The following result is now however a simple consequence of this 
discussion.

\begin{theorem} \label{t 2.14.}
  Let $G$ be either a closed oriented graph manifold with $\sigma (G) 
=$ 0, or a compact oriented graph manifold with incompressible 
boundary, so that $\partial G$ is a finite union of tori.

 Then there exists a minimizing sequence $\{g_{i}\}$ for 
$\mathcal{S}_{-}^{2}$ on $G$, i.e.
$$\mathcal{S}_{-}^{2}(g_{i}) \rightarrow  0, $$
which performs the geometrization of $G$.

 If in addition $G$ is irreducible, then there exist minimizing 
sequences for $I_{\varepsilon}^{~-},$ for any $\varepsilon  > 0$, 
which perform the geometrization of $G$.
\end{theorem}
\noindent
{\bf Proof:} 
 We construct the minimizing sequence $\{g_{i}\}$ to satisfy each of 
the conditions (i)-(iii) in turn.

 (i). Let $\{S_{q}^{2}\}$ be a collection of embedded 2-spheres in $G$ 
giving the sphere decomposition of $G$. For each such $S_{q}^{2},$ 
associate the isometrically doubled Schwarzschild metric $g_{i} = 
g_{S}(m_{i})$ in (1.23), with mass $m_{i} \rightarrow $ 0, defined on 
$S^{2}\times {\mathbb R}  = {\mathbb R}^{3}\#{\mathbb R}^{3}.$ Observe that 
this sequence of metrics converges to the union of two copies of ${\mathbb 
R}^{3},$ glued together at one point; the 2-sphere connecting the 
${\mathbb R}^{3}$ factors is crushed to a point as $i \rightarrow  \infty 
.$ Further note that $\rho_{i}(x) \sim  m_{i},$ for any $x\in 
S_{q}^{2}$ and $g_{i}$ is scalar-flat.

 Choose a sequence $\delta_{i} \rightarrow $ 0, with $m_{i} <<  
\delta_{i},$ and let $g_{i}|_{B_{q}(\delta_{i})}$ be the metric 
$g_{S}(m_{i})$ restricted to the $\delta_{i}$ tubular neighborhood of 
the horizon $\Sigma  = S^{2}(m_{i}).$ Thus, $g_{i}$ in the annular 
region $A_{q}(\delta_{i}/2, \delta_{i})$ differs from the flat metric 
on the order of $m_{i}/\delta_{i} << $ 1. The blow-ups $g_{i}'  = 
\rho_{i}(x)^{-2}\cdot  g_{i}$ based at $x$ converge to the complete 
isometrically doubled Schwarzschild metric $g_{S}$ of mass $m \sim $ 1.

 (ii). Next consider minimizing sequences for $\mathcal{S}  _{-}^{2}$ on 
the prime factors $G_{p}$ of $G$ first with $\sigma (G_{p}) =$ 0. Any 
minimizing sequence $\{g_{i}\}$ for $\mathcal{S}  _{-}^{2}$ with bounded 
curvature will volume collapse these factors, except possibly when 
$G_{p}$ is flat. In accordance with case (i) of Lemma 2.12, any such 
sequence may be altered to a sequence, still called $\{g_{i}\}$ and 
volume collapsing with bounded curvature, which has no solid tori 
components in its Seifert fibered factors. Lemma 2.13 implies that such 
a decomposition is unique up to isotopy, and induces the JSJ torus 
decomposition.

 For the exceptional cases where $G_{p}$ is a closed flat, Nil or Sol 
manifold, or the $S^{1}$ bundle over the M\"obius band, (so that 
$G_{p}$ has an empty torus decomposition), choose $\{g_{i}\}$ to be 
geometric, volume collapsing in the last three cases.

 Similarly on the factors with $\sigma (G_{p}) > 0$, choose $\{g_{i}\}$ 
to be geometric, and so $\{g_{i}\}$ is the constant sequence, c.f. 
(iii) above.

 We then glue together the metrics $g_{i}$ on the prime factors $G_{p}$ 
to form a metric on $G$. First, consider the factors $G_{p}$ with 
$\sigma (G_{p}) =$ 0. On such factors, although the metric $g_{i}$ is 
very collapsed for $i$ large, on scales much smaller than the collapse 
scale the metric $g_{i}$ on each $G_{p}$ is close to a flat metric on a 
ball. Thus, choose $\delta_{i}$ in (i) above so that $\delta_{i} <<  
diam_{g_{i}}\mathcal{O} (y),$ for all $y\in G_{p}.$ The geodesic balls in 
$(G_{p}, g_{i})$ are then almost flat, and in particular are 
topological balls. It follows that we may remove such balls from each 
$G_{p},$ and metrically glue in the Schwarzschild metrics $g_{i}$ from 
(i). This can be done so that the full curvature, and so in particular 
the scalar curvature, of the glueing remains bounded in the glueing 
region $A_{q}(\delta_{i}/2, \delta_{i})$ as $\delta_{i} \rightarrow $ 
0. This construction is the same as the end of the proof of Theorem 2.9 
and is carried out in detail in [4, \S 6.2] or [5, \S 6.1]; we 
refer there for further details.

 Similarly, the prime factors with $\sigma (G_{p}) > 0$, i.e. $S^{3}/ 
\Gamma$ and $S^{2}\times S^{1}$ with geometric metric, are glued to the 
other prime factors by Schwarzschild necks. Note that for $g$ a 
geometric metric on $S^{3}/ \Gamma$ or $S^{2}\times S^{1}$, 
$\mathcal{S}_{-}^{2}(g) = 0$. However, for any metric on $S^{2}\times S^{1}$, 
$I_{\varepsilon}^{~-} > 0$, for any fixed ${\varepsilon} > 0$; one can 
make $I_{\varepsilon}^{~-} \rightarrow 0$ only by collapsing the 
$S^{1}$ factor with bounded curvature.

 (iii). It remains to geometrize the Seifert fibered factors in each 
prime factor $G_{p}.$ This has already been done in (ii) for the 
factors with $\sigma (G_{p}) > $ 0, or for factors with flat, Nil or 
Sol geometries.

 Thus assume $\sigma (G_{p}) =$ 0, and so $\pi_{1}(G_{p})$ is infinite. 
It follows that $\pi_{1}(S_{j})$ is infinite, for all $j$. For if 
$\pi_{1}(S)$ is finite, for some $S = S_{j},$ then $S = S^{3}/\Gamma .$ 
Hence, since $S$ is closed, $S = G_{p}$, and $\sigma (G_{p}) > 0$, a 
contradiction.

 By standard results, c.f. [26, Ch.5],  the $S^{1}$ orbits inject in 
$\pi_{1}(S_{j}),$ for each $j$, (and thus inject in $\pi_{1}(G)$ by 
Lemma 2.13). Thus, the $S^{1}$-collapsed geometry of $(S_{j}, g_{i})$ 
may be unwrapped by passing to sufficiently large, (depending on $i$), 
finite covers $\bar S_{j},$ so that the geometry of $(\bar S_{j}, 
g_{i})$ is bounded; thus the injectivity radius of $\{g_{i}\}$ at a 
fixed base point $x_{j} \in  \bar S_{j}$ is bounded away from 0 and 
$\infty .$ The metrics $g_{i}$ then (sub)-converge to a limit metric 
$\bar g$ which is a complete $S^{1}$ invariant metric on $\bar S_{j}.$ 

 Now of course for an arbitrary minimizing sequence $\{g_{i}\}$ 
constructed as in (ii) above, the limit metric $\bar g$ on $\bar S_{j}$ 
will not be geometric, i.e. will not be given by one of the Seifert 
fibered geometries. However, we simply choose $\{g_{i}\}$ so that its 
limit $\bar g$ on $\bar S_{j}$ is a complete Seifert fibered geometry.

 For the same reasons as above, the metrics $g_{i}$ on the $L_{k}$ 
factors converge, when lifted to suitable ${\mathbb Z}\oplus{\mathbb Z} $ 
covers, to complete $T^{2}$ invariant metrics on $T^{2}\times {\mathbb R} 
.$ Observe that the Seifert fibered geometries are rigid within their 
class, in that one cannot continuously deform one Seifert fibered 
geometry class to another class. Similarly, distinct $S_{j}$ components 
with the same {\it  complete}  non-compact Seifert geometry cannot be 
glued together within any fixed geometry. Thus, the $L_{k}$ factors, 
which serve to glue together the Seifert fibered factors $S_{j},$ 
cannot be made geometric. 

 This completes the description of the minimizing sequence $\{g_{i}\}$ 
for $\mathcal{S}  _{-}^{2}.$ Finally, note that if $G$ is irreducible, so 
that Step (i) above is vacuous, one may choose $\{g_{i}\}$ to be volume 
collapsing with uniformly bounded curvature and hence $\{g_{i}\}$ is a 
minimizing sequence for $I_{\varepsilon}^{~-},$ for any $\varepsilon  > 
$ 0.

{\endproof}

\section{Metric Surgery on Spheres in Asymptotically Flat Ends.}
\setcounter{equation}{0}

 In this section, we prove Theorem 0.4. Thus, throughout this section, 
we assume that $M$ is a closed, oriented 3-manifold, which is not tame, 
but which is spherically tame. Hence, as discussed in Theorem 1.3, 
there is a sequence $\varepsilon = \varepsilon_{i} \rightarrow 0$, 
and a sequence of minimizing pairs $(\Omega_{\varepsilon}, g_{\varepsilon}), 
\Omega_{\varepsilon} \subset\subset  M$, together with a complete 
blow-up limit $(N, g')$ which has an asymptotically flat end $E$, 
topologically of the form $S^{2}\times {\mathbb R}^{+}$ outside a compact 
subset of $E$. Note that since $N \subset\subset  
\Omega_{\varepsilon},$ for $\varepsilon $ sufficiently small, one has 
$N \subset\subset  M$, so that in particular the 2-sphere $S^{2}$ in 
$E$ embeds in $M$, i.e.
\begin{equation} \label{e3.1}
S^{2} \subset  M. 
\end{equation}

\begin{theorem} \label{t 3.1.}
  Suppose $M$ is a spherically tame, but not tame, closed and oriented 
3-manifold, with $\sigma (M) \leq $ 0. Then the 2-sphere $S^{2} \subset 
 M$ in (3.1) is essential in $M$. In particular, $M$ is reducible.
\end{theorem}

\noindent
{\bf Proof:}
  The proof is by contradiction, so we assume that the 2-sphere $S^{2}$ 
in (3.1) bounds a 3-ball in $M$. The idea is similar to the proof of 
Theorem 2.9, in that we glue on a specific comparison 3-ball onto 
$S^{2}$, analogous to the specific metric Dehn surgery carried out in 
Theorem 2.9. The details of this construction are somewhat more 
involved however.

 Since $(E, g' )$ is asymptotically flat, there is a compact set $K 
\subset  E$ and a ball $D \subset {\mathbb R}^{3}$ such that $E\setminus 
K$ is diffeomorphic to ${\mathbb R}^{3}\setminus D,$ and in a suitable 
chart on $E\setminus K$, the metric $g' $ has the expansion
\begin{equation} \label{e3.2}
g'_{ij} = (1+\frac{2m}{r})\delta_{ij} + h, 
\end{equation}
where 
\begin{equation} \label{e3.3}
|h| = O(r^{-2}), |D^{k}h| = O(r^{-2- k}), k = 1,2 
\end{equation}
with $m > $ 0.

 As observed following (1.26) or in the proof of Proposition 2.1, the estimate 
(1.26) implies that the metrics $g_{\varepsilon}'$ from (1.19) converge 
in the weak $L^{3,p}$ and (strong) $C^{2, \alpha}$ topology to the limit 
metric $g'$, uniformly on compact subsets of $(N, g')$. Thus for any given 
$R$ large but fixed, the metric $g_{\varepsilon}'$ from (1.19) also has the 
form (3.2)-(3.3) on $B_{y_{\varepsilon}}(2R),$ 
i.e. for large $r \leq 2R$,
\begin{equation} \label{e3.4}
(g_{\varepsilon}' )_{ij} = (1+\frac{2m}{r})\delta_{ij} + 
h_{\varepsilon}, 
\end{equation}
where 
\begin{equation} \label{e3.5}
|h_{\varepsilon}| = O(r^{-2}), |D^{k}h_{\varepsilon}| = O(r^{-2- k}), k 
= 1,2, 
\end{equation}
on $B(2R)\setminus D \subset  {\mathbb R}^{3}.$ The expressions (3.4)-(3.5) 
are valid for all $\varepsilon $ sufficiently small depending on the 
choice of $R$, and $h_{\varepsilon}$ converges in the $C^{2, \alpha}$ 
topology to $h$ as $\varepsilon  \rightarrow 0$.

 The smooth domain $B(2R)\setminus D$  embeds smoothly in 
$\Omega_{\varepsilon}$ for $\varepsilon$ sufficiently small, and hence 
also in $M$. In particular, for all $R$ large but fixed, the canonical 
round 2-spheres $S^{2}(R) \subset  {\mathbb R}^{3}\setminus D  \cong 
E\setminus K$, centered at $0 \in {\mathbb R}^{3}$, embed in $M$.

 By the assumption above, $S^{2}(R)$ bounds a 3-ball $B^{3}$ in $M$, so 
that in particular $S^{2}(R)$ disconnects $M$ into two components, 
consisting of the {\it inside} component containing the base point $y$ 
of $N$, and the {\it outside} component containing the unbounded part 
of the end $E \subset  N$. We will say that $S^{2}(R)$ bounds a 3-ball 
$B^{3}$ on the inside, (or outside), if the inside (or outside) 
component is a 3-ball. Obviously this notion is independent of $R$, for 
all $R$ large.

\medskip

{\bf Case I.}
 Suppose that $S^{2}(R)$ bounds $B^{3}$ on the inside. 

\noindent
First, for clarity, we glue in a comparison 3-ball $\bar B$ onto the 
spherically symmetric metric
\begin{equation} \label{e3.6}
g_{S} = (1+\frac{2m}{r})\delta_{ij}, 
\end{equation}
at $S^{2}(R).$ Following this, it is shown that essentially the same 
construction and estimates are valid when glueing onto the limit metric 
$g' $ in (3.2) or the approximations $g_{\varepsilon}' $ in (3.4). 
A simple computation shows that $g_{S}$ has positive scalar curvature, 
$s_{g_{S}} > 0$, so that $(s_{g_{S}})^{-} \equiv 0$.

 In the metric $g_{S}$, the sphere $S^{2}(R)$ is isometric to a 
constant curvature sphere of Euclidean radius 
$R(1+\frac{2m}{R})^{1/2}.$ Let $(B_{o}, g_{o})$ be the flat 
Euclidean ball of the same radius, so that $S^{2}(R)$ is filled in on 
the ``inside'' with a flat Euclidean metric $g_{o}.$ The metrics 
$g_{S}$ and $g_{o}$ agree $C^{o}$ at the boundary, but not $C^{1}.$ Let 
$\bar g$ denote the union of these two metrics, as a metric on ${\mathbb 
R}^{3}.$

 The main point now is that $S^{2}(R)$ is more convex in the Euclidean 
metric than in the $g_{S}$ metric. Thus, let $A$ denote the $2^{\rm nd}$ 
fundamental form of $S^{2}(R)$, w.r.t. the outward normal, with 
$A_{g_{o}}$ denoting the $2^{\rm nd}$ fundamental form w.r.t. the Euclidean 
metric, and $A_{g_{S}}$ the form w.r.t. the $g_{S}$ metric; the sign 
conventions are such that $A > 0$ for spheres in Euclidean space. Then 
a short computation using (3.6) gives
\begin{equation} \label{e3.7}
A_{g_{o}}=\frac{1}{R(1+\frac{2m}{R})^{1/2}} I\sim 
\bigl(\frac{1}{R}-\frac{m}{R^{2}}+O(R^{-3})\bigr)I, \  A_{g_{S}} = 
\frac{1-\frac{m}{R}(1+\frac{2m}{R})^{-1}}{R(1+\frac{2m}{R})} I\sim 
\bigl(\frac{1}{R}-\frac{2m}{R^{2}}+O(R^{-3})\bigr)I, 
\end{equation}
where $I$ denotes the identity matrix. Thus
\begin{equation} \label{e3.8)}
A_{g_{o}} -  A_{g_{S}} =  \frac{m}{R^{2}}I + O(R^{-3}) > 0, 
\end{equation}
since $m > $ 0 and $R$ is assumed sufficiently small.

 This means that there is a concentration of positive curvature for 
$\bar g$ at the seam $S^{2}(R).$ It follows that the metric $\bar g$ 
may be smoothed to a metric $\widetilde g$ in the annulus $A(R- 1, 
R+1)$ so that it agrees with the flat metric at $S^{2}(R- 1)$ and with 
$g_{S}$ at $S^{2}(R+1),$ and satisfies
\begin{equation} \label{e3.9}
(\widetilde s)^{-} = (s_{g_{S}})^{-} = 0, 
\end{equation}
in $A(R- 1, R+1).$ Of course $\widetilde g$ agrees with $g_{S}$ outside 
the ball $B(R+1).$

 To verify (3.9), one may write $\bar g$ explicitly as a conformally 
flat warped product
\begin{equation} \label{e3.10}
\bar g  = dt^{2} + f^{2}(t)ds^{2}_{S^{2}}; 
\end{equation}
here $f(t) = t$, for $t \leq  \bar R \equiv  R(1+\frac{2m}{R})^{1/2},$ 
while $f(t) = r(1+\frac{2m}{r})^{1/2},$ for $r = r(t) \geq  \bar R.$ 
The scalar curvature of a metric of the form (3.10) is given by
$$\frac{1}{2}s = - 2\frac{f''}{f} + \frac{1- (f' )^{2}}{f^{2}}, $$
while $A = (f'/f) \cdot I$. The estimate (3.8) means that when 
smoothing the Lipschitz function $f$ at the seam $\bar R$ to obtain a 
smooth function $\widetilde f,$ one has $\widetilde f''  \leq $ 0 near 
the seam. Hence the smoothing $\widetilde g$ increases the scalar 
curvature, which gives (3.9).

 Next, we make a similar estimate for the change in volume, as above 
working with $g_{S}.$ Clearly $vol_{g_{o}}B(\bar R) = 
\omega_{3}\bar R^{3} = \nu_{2}\bar R^{3}/3,$ where $\omega_{3}$ 
($\nu_{2}$) is the volume of the Euclidean unit 3-ball (2-sphere). 
Another short computation using (3.6) shows that the volume of the 
region $B(R)$ in the end $(E, g_{S})$ satisfies
$$vol_{g_{S}}B(R) \geq  (\frac{1}{3}R^{3} + \frac{3}{2}mR^{2} -  
O(R))\nu_{2} = (R^{3} + \frac{9}{2}mR^{2} - O(R))\omega_{3}. $$
From the definition of $\bar R,$ it follows that
$$vol_{g_{o}}B_{o} -  vol_{g_{S}}B(R) \leq  
-\frac{3}{2}mR^{2}\omega_{3} + O(R), $$
so that, for $R$ sufficiently large, 
\begin{equation} \label{e3.11}
vol_{g_{o}}B_{o} <  vol_{g_{S}}B(R). 
\end{equation}

 Hence, one sees that $\widetilde g$ has non-negative scalar curvature 
and less volume than $g_{S},$ and so is an ``admissible'' comparison 
metric to $g_{S},$ in the sense of the minimizing properties stated in 
Theorem 1.3. We now claim that essentially the same construction and 
estimates may be carried out w.r.t. the metric $g' $ in (3.2) in place 
of the metric $g_{S}.$

 Thus, consider the sphere $S^{2}(R)$ w.r.t. the metric $g' ,$ for $R$ 
large. Since $h$ in (3.2) is sufficiently small in the $C^{2}$ 
topology, $(S^{2}(R), g' )$ has positive Gauss curvature. Hence, by the 
Weyl-Alexandrov embedding theorem, c.f. [35], $(S^{2}(R), g' )$ may be 
isometrically embedded in ${\mathbb R}^{3}$ as a convex surface. In 
particular, $S^{2}(R)$ bounds a convex domain $\bar B \subset  {\mathbb 
R}^{3},$ where ${\mathbb R}^{3}$ is given the flat metric $g_{o}$ as 
before. One may then glue in $\bar B$ to the boundary $S^{2}(R),$ to 
obtain a new complete metric $\bar g' $ on $\bar N,$ diffeomorphic to 
${\mathbb R}^{3}.$ This metric is $C^{o},$ and piecewise $C^{\infty},$ but 
is not $C^{1}$ at the seam $S^{2}(R).$

 Using the decay estimates on $h$, it follows as in (3.7)-(3.8) that
\begin{equation} \label{e3.12)}
A_{g_{o}} -  A_{g'} =  \frac{m}{R^{2}}I + o(R^{-2}) > 0, 
\end{equation}
for $R$ sufficiently large. As before, the metric $\bar g' $ may be 
smoothed to a metric $\widetilde g$ near the seam $S^{2}(R),$ in a 
manner similar to the smoothing following (3.10). More precisely, at 
the seam $(S^{2}(R), g')$ consider exponential normal coordinates 
w.r.t. $g'$ on the outside and w.r.t. $g_{o}$ on the inside. The metric 
$\bar g$ then has the form $\bar g = dt^{2} + g_{t}$, where $t$ is the 
signed distance to $S^{2}(R)$ and $g_{t}$ is a curve of metrics on 
$S^{2}$. The family $g_{t}$ is piecewise $C^{\infty}$ and Lipschitz 
through the seam $S^{2}(R)$. The second fundamental form $A$ in these 
coordinates is given by $A = dg/dt$. As above with (3.10), the estimate 
(3.12) implies there is a smoothing $\widetilde g$ near $S^{2}(R)$ so 
that $\widetilde s \geq 0$ near $S^{2}(R)$. Hence
\begin{equation} \label{e3.13}
(\widetilde s) ^{-} \geq  (s_{g'}) ^{-}, 
\end{equation}
as in (3.9). For the same reasons, the volume estimate (3.11) holds w.r.t. 
$\widetilde g$ and $g' .$

 Finally, we claim that these same estimates hold on the approximations 
$g_{\varepsilon}' $ in place of the limit $g' .$ This holds since by 
(3.4)-(3.5), $g_{\varepsilon}' $ has the same form as $g' ,$ i.e. the 
lower order terms $h_{\varepsilon}$ obey the same estimates as the 
lower order term $h$ for $g' .$ 

  It follows that there are metrics $\widetilde g_{\varepsilon}' ,$ 
(close to $\widetilde g ),$ such that for $\varepsilon $ sufficiently 
small,
\begin{equation} \label{e3.14}
(s_{\widetilde g_{\varepsilon}'}) ^{-} \geq  (s_{g_{\varepsilon}'}) 
^{-}, 
\end{equation}
in $A(R- 1, R+1),$ while $s_{\widetilde g_{\varepsilon}'} \equiv $ 0 in 
$B(R- 1)$ and $\widetilde g_{\varepsilon}'  = g_{\varepsilon}' $ 
outside $B(R+1).$ In addition, one has
$$vol_{\widetilde g_{\varepsilon}'}B(R+1) <  
vol_{g_{\varepsilon}'}B(R+1), $$
and hence
\begin{equation} \label{e3.15}
vol_{\widetilde g_{\varepsilon}'}\widetilde \Omega_{\varepsilon} <  
vol_{g_{\varepsilon}'}\Omega_{\varepsilon}, 
\end{equation}
where $\widetilde \Omega_{\varepsilon}$ is $\Omega_{\varepsilon}$ with 
the component bounding $S^{2}$ on the inside replaced by a 3-ball.

 We now compare the values of the functional $I_{\varepsilon}^{~-}$ on 
$\widetilde g_{\varepsilon}' $ and $g_{\varepsilon}' .$ From (3.14), 
and the fact $\widetilde g_{\varepsilon}' $ is scalar-flat in $B(R- 
1)$,
\begin{equation} \label{e3.16}
\int_{\widetilde \Omega_{\varepsilon}}(s_{\widetilde 
g_{\varepsilon}'}^{-})^{2}dV_{\widetilde g_{\varepsilon}'} \leq  
\int_{\Omega_{\varepsilon}}(s_{g_{\varepsilon}'}^{-})^{2}dV_{g_{\varepsilon}'}.\end{equation}
Thus, the estimates (3.15) and (3.16) imply that $\mathcal{S}  
_{-}^{2}(\widetilde g_{\varepsilon}' ) <  \mathcal{S}  
_{-}^{2}(g_{\varepsilon}' ).$

 Next, to compare the $L^{2}$ norm of $z$ on both metrics, the flat 
metric and the metric $g_{\varepsilon}' $ differ on the order of 
$R^{-1}$ near $S^{2}(R),$ while their curvatures differ on the order of 
$R^{-3}.$ It follows that the smoothing $\widetilde g_{\varepsilon}' $ 
in $A(R- 1, R+1)$ may be done so that the curvatures of $\widetilde 
g_{\varepsilon}'$ near the seam $S^{2}(R)$ are on the order of at most 
$R^{-2}.$ In terms of the discussion above on (3.10), this arises from 
the fact that since $f = O(R)$ and the jump in $f'$ at the seam is on 
the order of $O(R^{-1})$, one may choose the smoothing $\widetilde f$ 
with $|\widetilde f''| = O(R^{-1}),$ so that $|\widetilde f'' /f| = 
O(R^{-2}).$ Thus,
\begin{equation} \label{e3.17}
\int_{A(R- 1, R+1)}|z_{\widetilde g_{\varepsilon}'}|^{2}dV_{\widetilde 
g_{\varepsilon}'} \leq  cR^{-2}, 
\end{equation}
since $volA(R-1, R+1) = O(R^{2})$; here $c$ is a constant independent 
of $R$.

 In $B(R- 1), \widetilde g_{\varepsilon}' $ is flat, so that of course 
$z =$ 0 in this region. On the other hand, $g_{\varepsilon}' $ has a 
definite amount of curvature in this region, (inside $S^{2}(R- 1)).$ 
For instance, since at the base point $y_{\varepsilon}, \rho' 
(y_{\varepsilon}) =$ 1, one has
$$\int_{B_{y_{\varepsilon}}(1)}|r_{g_{\varepsilon}'}|^{2}dV_{g_{\varepsilon}'}
\geq  c_{o}vol B_{y_{\varepsilon}}(1),$$ 
and the same for limit metric $g'$ based at $y$. Since the limit $(N, 
g')$ is not of constant curvature on any open set, (by analyticity), we 
then also have
\begin{equation} \label{e3.18}
\int_{B_{y_{\varepsilon}}(1)}|z_{g_{\varepsilon}'}|^{2}dV_{g_{\varepsilon}'} 
\geq  d_{o}vol B_{y_{\varepsilon}}(1),
\end{equation}
for $\varepsilon $ sufficiently small, where $d_{o} = d_{o}(c_{o}).$ 
Comparing (3.17) and (3.18), it follows that by choosing $R$ 
sufficiently large, and $\varepsilon $ sufficiently small, one obtains
\begin{equation} \label{e3.19}
\int_{\widetilde \Omega_{\varepsilon}}|z_{\widetilde 
g_{\varepsilon}'}|^{2}dV_{\widetilde g_{\varepsilon}'} <  
\int_{\Omega_{\varepsilon}}|z_{g_{\varepsilon}'}|^{2}dV_{g_{\varepsilon}'}. 
\end{equation}
Combining the estimates (3.15), (3.16), (3.19), giving
\begin{equation} \label{e3.20}
I_{\varepsilon}^{~-}(\widetilde g_{\varepsilon}' ) <  
I_{\varepsilon}^{~-}(g_{\varepsilon}' ). 
\end{equation}
For $\varepsilon > 0$ fixed, the metrics $\widetilde g_{\varepsilon}'$ 
and $g_{\varepsilon}'$ are not (necessarily) defined on $M$. However, 
both are limits of sequences $\{\widetilde g_{i}\}$, $\{g_{i}'\}$ of 
metrics on $M$ for which $I_{\varepsilon}^{~-}$ converges to its value 
on $\widetilde g_{\varepsilon}'$ and $g_{\varepsilon}'$. Hence, (3.20) 
contradicts the minimizing property of $g_{\varepsilon}'$ or 
$g_{\varepsilon}$ in (1.8); we recall that $I_{\varepsilon}^{~-}$ is 
scale invariant. This proves Theorem 3.1 in case $S^{2}(R)$ bounds a 
3-ball on the inside.

\medskip

{\bf Case II.}
 Suppose $S^{2}(R)$ bounds $B^{3}$ on the outside. 

\noindent
We proceed as above to construct a suitable comparison metric, although 
this case requires a little more delicate consideration. In particular 
we will need to make stronger use of the decay estimate (3.3) on $h$; 
this estimate plays only a minor role in Case I, (c.f. also Theorem 4.2 
below). 

 Observe that if $S^{2}(R)$ bounds on the inside, in place of glueing 
in a flat 3-ball $\subset  {\mathbb R}^{3}$ as in Case I above, one may 
instead glue in the 3-ball contained in a very large 3-sphere 
$S^{3}(\delta^{-1})$ of radius $\delta^{-1};$ for $\delta $ 
sufficiently small, all the estimates above remain valid. The idea now 
is that if instead $S^{2}(R)$ bounds on the outside, we glue in the 
complementary (very) large 3-ball in $S^{3}(\delta^{-1}).$ 

 Note that in this case, since the end $E$ has infinite volume,
\begin{equation} \label{e3.21}
vol_{g_{\varepsilon}'}B^{3} \rightarrow  \infty , \ \ {\rm as} \ \ 
\varepsilon  \rightarrow  0, 
\end{equation}
where $B^{3}$ is the 3-ball bounding $S^{2}(R)$ on the outside.

 We first construct the comparison metric $\widetilde g$ on the limit, 
and then approximate it to obtain a comparison metric $\widetilde 
g_{\varepsilon}' $ for $g_{\varepsilon}' ,$ as in Case I. The end $E$ 
of $(N, g')$ is asymptotically flat, so satisfies (3.2). It is 
convenient, although not necessary, to rescale the metric $g' $ so that 
the mass is normalized to $m = \frac{1}{2};$ observe that the mass 
scales as distance. This normalization eliminates the dependence of the 
estimates to follow on the mass $m$. In particular, assuming $m = 
\frac{1}{2}$ from now on, one has
\begin{equation} \label{e3.22}
\int_{A(R,\infty )}|z|^{2}dV_{g'} \sim  R^{-3}, 
\end{equation}
as $R \rightarrow  \infty .$ Now as before, write the metric $g' $ in 
(3.2) as 
\begin{equation} \label{e3.23}
g'  = g_{S} + h, 
\end{equation}
where 
\begin{equation} \label{e3.24}
g_{S} = (1+\frac{1}{r})g_{Eucl} \ \ {\rm and} \ \ h = O(r^{-2}). 
\end{equation}
Given $R$ large, let $\delta $ and $D$ be the solutions to the equations
\begin{equation} \label{e3.25}
\sin \delta D = \delta R(1+\frac{1}{R})^{1/2},
\end{equation}
$$\cos \delta D = 1 -  \frac{1}{2R}(1-\frac{1}{R})^{-1}. $$
These equations mean that the geodesic sphere $S^{2}(D)$ of radius $D$ 
in $S^{3}(\delta^{-1}),$ about some base point $x_{o},$ is isometric to 
the sphere $S^{2}(R)\subset (E, g_{S}),$ and the $2^{\rm nd}$ fundamental 
form of $S^{2}(D)\subset S^{3}(\delta^{-1})$ satisfies 
\begin{equation} \label{e3.26}
A = \frac{\delta \cos \delta D}{\sin \delta D}I = 
R^{-1}(1+\frac{1}{R})^{-1/2}(1 -  \frac{1}{2R}(1+\frac{1}{R})^{-1})I, 
\end{equation}
agreeing with the $2^{\rm nd}$ fundamental form of $S^{2}(R) \subset (E, 
g_{S}).$ Thus in the metric $g_{S},$ the boundary $S^{2}(R)$ is 
isometric, and has identical $2^{\rm nd}$ fundamental form, to 
$S^{2}(D)\subset S^{3}(\delta^{-1}).$ Let $g_{\delta}$ denote the 
(round) metric on $S^{3}(\delta^{-1}).$ Note that $\delta^{2} = 
O(R^{-3}),$ so that the curvatures of $g_{S}$ and $(S^{3}, g_{\delta})$ 
are on the order of $R^{-3}.$

 As in Case I, for the moment, we work with the metric $g_{S}.$ It 
follows that if one attaches the complementary geodesic ball $B = 
B_{y_{o}}(2\pi\delta^{-1}- D)\subset S^{3}(\delta^{-1})$, to 
$S^{2}(D)$, where $y_{o}$ is the antipodal point to $x_{o},$ then the 
resulting metric $\bar g$ consisting of $g_{S}$ on the inside and $(B, 
g_{\delta})$ on the outside is piecewise $C^{\infty},$ and is $C^{1}$ 
smooth at the seam $S^{2}(R).$ We then smooth the seam $S^{2}(R),$ in a 
conformally flat way as in (3.10), in a band $A = A(R- 1, R+1)$ about 
$S^{2}(R).$ Because the curvatures of $g_{S}$ and $g_{\delta}$ are on 
the order of $O(R^{-3})$ and the metrics agree $C^{1}$ at the seam, one 
obtains in this way a smooth metric $\widetilde g$ satisfying
\begin{equation} \label{e3.27}
\int_{A}|z|^{2}dV_{\widetilde g} \leq  c\cdot  R^{-4}.
\end{equation}

 Note that $z =$ 0 past $S^{2}(R+1),$ i.e. outside $S^{2}(R+1).$ 
Further, since the scalar curvature of $S^{3}(\delta^{-1})$ is positive 
and $g_{S}$ is scalar-flat, the smoothing may be done so that the 
metric $\widetilde g$ has $\widetilde s \geq $ 0 pointwise. It follows 
that $\widetilde g$ is an admissible comparison metric to $g_{S},$ in 
that it has less volume than $g_{S}$ and has non-negative scalar 
curvature. By comparing (3.22) and (3.27), one sees that
\begin{equation} \label{e3.28}
\int|z_{\widetilde g}|^{2}  <  \int|z_{g_{S}}|^{2}. 
\end{equation}

 Next, we deal with the lower order term $h$ in (3.2). First, the 
curvatures of $g' $ and $g_{S}$ in the glueing region above differ on 
the order of at most $O(R^{-4}),$ since the metrics differ by 
$O(R^{-2})$ and $|D^{k}h| = O(R^{-2-k}), k =$ 1,2. Now redefine 
$\delta $ and $D$ so as to solve the system
\begin{equation} \label{e3.29}
\sin \delta D = \delta R(1+\frac{1}{R})^{1/2}, 
\end{equation}
$$\cos \delta D = 1 -  \frac{1}{2R}(1-\frac{1}{R})^{-1} -  
\frac{1}{R^{\lambda}}, $$
where $\lambda $ is a fixed number in $(\frac{3}{2},$ 2). For this 
choice of $\delta$ and $D$, the geodesic sphere $S^{2}(D)\subset 
S^{3}(\delta^{-1})$ is isometric to $S^{2}(R)$ in $g_{S},$ while its 
$2^{\rm nd}$ fundamental form $A_{g_{\delta}}$ is smaller than the 
$2^{\rm nd}$ fundamental form $A_{g_{S}};$ in fact, 
\begin{equation} \label{e3.30}
A_{g_{S}} -  A_{g_{\delta}} \sim  \frac{1}{R^{1+\lambda}} I >  0. 
\end{equation}
Thus, as in (3.8), there is a concentration of positive curvature at 
the seam.

 As above, use the Weyl embedding theorem to isometrically embed 
$(S^{2}(R), g' )$ in $S^{3}(\delta^{-1}).$ Since $h$ satisfies the 
decay estimates (3.3), $(S^{2}(R), g' )$ is close to $(S^{2}(D), 
g_{\delta}),$ and has $2^{\rm nd}$ fundamental form $A_{g_{\delta}}$ still 
satisfying (3.30), i.e.
\begin{equation} \label{e3.31}
A_{g'} -  A_{g_{\delta}} \sim  \frac{1}{R^{1+\lambda}} I. 
\end{equation}
Thus, to the metric $g' $ at the boundary $S^{2}(R)$, attach on the 
outside the large domain $\bar B$ in $S^{3}(\delta^{-1})$ with 
$\partial\bar B = S^{2}(R),$ so that $\bar B$ is a small perturbation 
of $B_{y_{o}}(2\pi\delta^{-1}- D)$.

 The resulting metric $\bar g' $ is $C^{o},$ and satisfies the estimate 
(3.31) at the seam $S^{2}(R).$ Exactly as discussed earlier following 
(3.12), this metric may be smoothed within the annulus $A(R- 1, R+1)$ 
to a metric $\widetilde g' $ so that
\begin{equation} \label{e3.32}
(\widetilde s' )^{-} \geq  (s' )^{-} = 0 ,
\end{equation}
everywhere. From the estimate (3.21), it is obvious that $\widetilde g' 
$ has less volume than $g' .$

 Further, using (3.31), the curvature of $\widetilde g' $ in $A(R- 1, 
R+1)$ is on the order of
\begin{equation} \label{e3.33}
|\widetilde z'|  = O(R^{-1-\lambda}). 
\end{equation}
This follows from the fact that the curvatures of $g'$ and $g_{\delta}$ 
are on the order of $O(R^{-3})$ together with (3.31) and the 
Gauss-Codazzi equations. Thus,
\begin{equation} \label{e3.34}
\int_{A(R- 1, R+1)}|\widetilde z'|^{2}dV_{\widetilde g'} = 
O(R^{-2\lambda}) = o(R^{-3}). 
\end{equation}
Comparing (3.34) with (3.22), $\widetilde g'$ hence has 
non-negative scalar curvature, less volume and less $L^{2}$ norm of 
curvature than $g' .$

 We are now in position to construct a comparison metric $\widetilde 
g_{\varepsilon}' $ to $g_{\varepsilon}' $ and compare the values 
$I_{\varepsilon}^{~-}(\widetilde g_{\varepsilon}' )$ and 
$I_{\varepsilon}^{~-}(g_{\varepsilon}' ).$ As $\varepsilon  = 
\varepsilon_{i} \rightarrow $ 0, the metrics $g_{\varepsilon}' $ 
converge smoothly to the limit $g' ,$ uniformly on compact sets. One may 
then fix a choice of $R$ sufficiently large, and then choose 
$\varepsilon $ sufficiently small so that in the region $A(R-10, 
R+10),$ the metric $g_{\varepsilon}' $ is of the form (3.2). It follows 
that the estimates and constructions above on $g' $ are equally valid 
for $g_{\varepsilon}' .$ Now from (3.21), for $\varepsilon $ 
sufficiently small,
$$vol_{\widetilde g_{\varepsilon}'}(\widetilde \Omega_{\varepsilon}) <  
vol_{g_{\varepsilon}'}(\Omega_{\varepsilon}), $$
where as in Case I, $\widetilde \Omega_{\varepsilon}$ is 
$\Omega_{\varepsilon}$ with the component bounding $S^{2}$ on the 
outside replaced by a 3-ball. Further, from the construction above, one 
has $(\widetilde s_{\varepsilon}' ) ^{-} \geq  (s_{\varepsilon}' ) 
^{-}.$ It follows that
$$\int_{\widetilde \Omega_{\varepsilon}}(s_{\widetilde 
g_{\varepsilon}'}^{-})^{2}dV_{\widetilde g_{\varepsilon}'} <  
\int_{\Omega_{\varepsilon}}(s_{g_{\varepsilon}'}^{-})^{2}dV_{g_{\varepsilon}'}. $$
and similarly, 
$$\int_{\widetilde \Omega_{\varepsilon}}|z_{\widetilde g_{\varepsilon}'}|^{2}dV_{\widetilde g_{\varepsilon}'} <  
\int_{\Omega_{\varepsilon}}|z_{g_{\varepsilon}'}|^{2}dV_{g_{\varepsilon}'}. $$
Thus, again one has
\begin{equation} \label{e3.35}
I_{\varepsilon}^{~-}(\widetilde g_{\varepsilon}' ) <  
I_{\varepsilon}^{~-}(g_{\varepsilon}' ), 
\end{equation}
contradicting, as in (3.20), the minimizing property of 
$g_{\varepsilon}' .$

 It follows that $S^{2}(R)$ cannot bound a 3-ball either on the inside 
or the outside, which completes the proof.

{\endproof}

\begin{remark} \label{r 3.2.}
  We point out that the comparison argument above strongly makes use of 
the flexibility in the functional $\mathcal{S}  _{-}^{2}$ or 
$I_{\varepsilon}^{~-},$ in that, because of the cutoff $s^{-} =$ 
min$(s, 0)$, one may ignore regions of the manifold where the scalar 
curvature of the metric is positive. A similar comparison argument for 
the more rigid functional $\mathcal{S}^{2}$ in (0.2), or the associated 
$I_{\varepsilon}$ as in (0.10), (with $\mathcal{S}^{2}$ in place of 
$\mathcal{S}  _{-}^{2}),$ would be much more difficult to carry out. In this 
situation, the blow-up limit $(N, g')$ is necessarily scalar-flat, and 
the allowable comparison metrics on the limit must also be scalar-flat. 
However, it is not clear that $g' $ admits any compact scalar-flat 
perturbations.

 Similarly, the comparison argument in Case II does not hold if one 
uses $|r|^{2}$ in place of $|z|^{2}$ in the definition (0.10) of 
$I_{\varepsilon}^{~-}.$
\end{remark}

\section{Asymptotically Flat Ends and Annuli.}
\setcounter{equation}{0}

  Theorems 0.2 and 0.4 are the main results of this paper, and as 
discussed in \S 0, reduce Conjectures I and II to the Sphere 
conjecture. This final section of the paper presents some remarks and 
results related to the Sphere conjecture, and so serves as a bridge to 
the sequel paper.

   Regarding the Sphere conjecture, Theorem 1.4 of course gives a 
natural condition implying that all ends of complete $\mathcal{Z}_{c}^{2}$ 
solutions $(N, g')$ are asymptotically flat. Proposition 4.1 below 
gives a relatively simple characterization of such limits which admit 
at least one asymptotically flat end, thus generalizing in a sense 
Theorem 1.4. (The proof however is just a minor variation of that of 
Theorem 1.4). The main result of this section, Theorem 4.2, shows that 
one may carry out metric sphere surgeries, as in the proof of Theorem 
3.1, under much weaker conditions than the assumption of an 
asymptotically flat end, at least from the inside. Hence, this result 
extends the domain of validity of Theorem 0.4. We include Theorem 4.2 
in this paper since the main ideas of the proof are similar to those of 
Theorem 3.1, although the technical details are somewhat different. 
These more technical issues in fact serve as an introduction to methods 
used more extensively in the sequel.

 We consider complete, non-flat $\mathcal{Z}_{c}^{2}$ solutions $(N, g' , 
\omega ),$ i.e. complete metrics satisfying the conclusions of Theorem 
1.3, with $\omega  = \tau  + \alpha s/12.$ In particular, the metric 
satisfies the equations (1.16)-(1.17). For convenience, set
\begin{equation} \label{e4.1}
u = -\omega , 
\end{equation}
so that $u > 0$ in the interior of the region where $s =$ 0. We begin 
with the following:

\begin{proposition} \label{p 4.1.}
  Let $(N, g, u)$ be a complete, non-flat $\mathcal{Z}_{c}^{2}$ solution. 

{\bf (i).}
 Suppose $u_{o} = sup_{N} u <  \infty .$ Then there is a constant 
$\delta_{o} > $ 0, depending only on $\alpha / u_{o}$, such that if a 
compact subset $C_{o}$ of the level set $L_{o} = \{u = 
u_{o}(1-\delta_{o})\}$ bounds a component of the superlevel set $U^{o} 
= \{u \geq  u_{o}(1-\delta_{o})\}$ in N, then (N, g) has an 
asymptotically flat end.

{\bf (ii).}
 There exists $K <  \infty $, depending only on $\alpha$, such that if 
a compact subset $C_{1}$ of the level set $L_{1} = \{u =$ K\}, (assumed 
non-empty), bounds a component of $U^{1} = \{u \geq $ K\} in N, then 
(N, g) has an asymptotically flat end on which $u$ is bounded.
\end{proposition}

\noindent
{\bf Proof:}
 {\bf (i).}
 Let $D^{o}$ be a component of $U^{o}$ such that $C_{o}$ = $\partial 
D^{o} \subset  L_{o}$ is compact. By the maximum principle applied to 
the trace equation (1.17), $D^{o}$ must be non-compact, and hence 
defines an end $E$ of $N$.

 We claim that $E$, or possibly a sub-end of $E$, is asymptotically 
flat. This follows by assembling results from [5,\S 7]. Namely one has
\begin{equation} \label{e4.2}
u \geq  u_{o}(1-\delta_{o}) >  0, 
\end{equation}
everywhere on $E$. Further, the oscillation of $u$ on $E$ satisfies 
$osc_{E }u \leq  \delta_{o}\cdot  u_{o}.$ By [5, Lemma 7.2], this 
implies that if $\delta_{o}$ is sufficiently small, then the curvature 
of $(E, g' )$ is everywhere small, in the sense that
\begin{equation} \label{e4.3}
|r|(x) \leq  \delta_{1}; 
\end{equation}
the constant $\delta_{1}$ depends only on $\delta_{o}.$ 

 It follows from the estimates (4.2) and (4.3), together with 
[5, Prop. 7.17] that
\begin{equation} \label{e4.4}
\limsup_{t \rightarrow \infty} |r| \rightarrow  0 
\end{equation}
in $E$, where $t(x)$ = dist$(x, \partial E)$. Now as noted at the end 
of [5,Prop.7.17], the estimates (4.4) and $u \geq u_{1}$, for some 
$u_{1} > 0$, are all that are required to carry out the proof of 
Theorem C in [5], i.e. Theorem 1.4 here, the point being that the 
proof takes place on each end individually. Hence a subend of $E$ is 
asymptotically flat.

{\bf (ii).}
 As above, the component $D^{1}$ of $U^{1} = \{u \geq K\}$ bounding 
$C_{1}$ defines an end $E$ of $N$. Now divide the Euler-Lagrange 
equations (1.16)-(1.17) by $K$, to obtain the equations
$$\frac{\alpha}{K}\nabla\mathcal{Z}^{2} -  L^{*}\bar u = 0, $$
$$\Delta\bar u = \frac{\alpha}{4K}|z|^{2}, $$
where $\bar u = \frac{u}{K}.$ If, for a given $\alpha$, $K$ is 
sufficiently large, these equations are close to the static vacuum 
equations. Since there are no complete non-flat static vacuum solutions 
by [4, Thm.3.2], it follows that the curvature is small sufficiently 
far out in $D^{1},$ i.e.
$$|r|(x) \leq  \delta_{1}, $$
for all $x\in D^{1}$ of distance at least $T$ to $C_{1},$ where 
$\delta_{1}$ depends only on a sufficiently large choice of $K$ and 
$T$; c.f. the proof of [5, Prop.7.17]. It then follows as before 
that (4.4) holds, and the remainder of the proof follows as in (i) 
above from [5, Thm.C].

{\endproof}

 For the work in the sequel paper, we will need a generalization of the 
metric sphere surgery given in Theorem 3.1. It is clear that one does 
not need a ``complete'' asymptotically flat end to carry out the proof 
of this result. It suffices to have a suitable spherical annulus $A = 
S^{2}\times I$ where the glueing takes place on which the metric $g' $ 
is sufficiently close to the flat metric, i.e. has the form (3.2)-(3.3).

 While the glueing on the outside, i.e. Case II of the proof of Theorem 
3.1, requires a strong estimate on the deviation of $g' $ from the flat 
metric, i.e. that $g' $ have the form (3.2) with $h = O(r^{-2}),$ a 
much weaker estimate suffices for glueings of 3-balls on the inside; 
this is often the more important case anyway, c.f. Remark 4.3.

  To describe this quantitatively, let $(N,g',y)$ be a complete 
non-flat $\mathcal{Z}_{c}^{2}$ solution, arising as a blow-up limit of 
$(\Omega_{\varepsilon}, g_{\varepsilon}, y_{\varepsilon}), \varepsilon  
= \varepsilon_{i} \rightarrow $ 0, as in Theorem 1.3. Let $v_{o}$, 
$\kappa $ be (arbitrary) small positive constants less than 
$\frac{1}{2}$, let $u_{o}$, $d$ be any positive constants, and let 
$D$, $R$ be large positive constants, with $R$ sufficiently large, 
(depending only on $\frac{\alpha}{u_{o}}$). We suppose there exists a 
component $A$ of the geodesic annulus $A_{y}((1-d)R, (1+d)R)$ about 
$y \in N$, which is topologically of the form $S^{2}\times I$, and 
which satisfies the global size bounds
\begin{equation} \label{e4.5}
vol A \geq  v_{o}\cdot  R^{3}, \ \  diam A \leq  D\cdot  R. 
\end{equation}
Further, suppose the potential function $u$ satisfies the oscillation 
bounds
\begin{equation} \label{e4.6}
osc_{A}u = \delta_{o}\cdot  u_{o}, \ \ \sup_{A}u \equiv  u_{o} > 0. 
\end{equation}
and that there is some level set $L_{o} = \{u = u_{o}(1-\delta' )\}$ of 
$u$ in $A$ such that 
\begin{equation} \label{e4.7}
L_{o} \subset  A_{\kappa d} = \{x\in A: dist(x, \partial A) = \kappa d 
R\}, 
\end{equation}
for some $\kappa  > $ 0, and that the sub-level set $U_{o} = \{u \leq  
u_{o}(1-\delta' )\}$ contains the inner boundary $\bar A\cap S_{y}((1- 
d)R)$ of $A$.

  Any asymptotically flat end on which $u$ is bounded away from 0 at 
infinity satisfies these conditions; in fact these conditions are much 
weaker than such an assumption. We then have:

\begin{theorem} \label{t 4.2.}
  Let $(N, g' , y)$ be a complete non-flat $\mathcal{Z}_{c}^{2}$ solution, 
as above, arising as a blow-up limit of $(\Omega_{\varepsilon}, 
g_{\varepsilon}, y_{\varepsilon}), \varepsilon  = \varepsilon_{i} 
\rightarrow $ 0, and satisfying the size and potential hypotheses 
(4.5)-(4.7).
 
  Then if $\delta_{o}$ is sufficiently small, depending only on $v_{o}, 
d, \kappa $ and $D$, the essential 2-sphere $S^{2} \subset  A \subset  
N$ cannot bound a 3-ball in $M$ on the inside of $S^{2}$. In 
particular, $N$ itself is a reducible 3-manifold, while $M$ is 
reducible on the inside of $S^{2}.$
\end{theorem}

\noindent
{\bf Proof:}
  For convenience, throughout the proof, we work in the scale where $R 
=$ 1, i.e. rescale $g' $ by $R^{-2}.$ For simplicity, we do not 
change the notation for $g' ,$ and often ignore the prime in the 
equations to follow. Thus, in this scale, (4.5) becomes
\begin{equation} \label{e4.8}
vol A \geq  v_{o}, \ \  diam A \leq  D. 
\end{equation}

 As in the proof of Proposition 4.1, there exists $\delta_{1} = 
\delta_{1}(\delta_{o}, \kappa ),$ which may be made arbitrarily small 
by requiring that $\delta_{o}$ is sufficiently small, such that
\begin{equation} \label{e4.9}
|r|(x) \leq  \delta_{1}, 
\end{equation}
for all $x\in A_{\kappa d}.$ We assume that $\delta_{o}$ is 
sufficiently small so that (4.9) implies that $A_{\kappa d}$ is 
diffeomorphic to a flat manifold. Hence $A_{\kappa d}$ carries a flat 
metric $g_{o},$ which is $\delta_{2} = \delta_{2}(\delta_{1})$ close to 
$g' $ in the $C^{1,\alpha}$ topology.

 The distance function $t$ from the base point $y$ in $N$, or 
equivalently (by subtracting a constant) from the inner boundary 
$\partial_{i}A = S_{y}(1- d)$ of $A$, is close to a distance function 
$t_{o}$ w.r.t. the flat metric $g_{o}$ on $A_{\kappa d}.$ Since $A$ is 
simply connected, the global size bounds (4.8) then imply that 
$(A_{\kappa d}, g_{o})$ may be isometrically immersed into a bounded 
domain $C$ in ${\mathbb R}^{3}$; the volume and diameter of $C$ are 
bounded, independent of $\delta$. This immersion is a covering map when 
restricted to $S^{2} \subset A$ and hence $A_{\kappa d}$ embeds 
isometrically onto the domain $C \subset {\mathbb R}^{3}$. One has 
$dist_{{\mathbb R}^{3}}(\partial_{i}C, \partial_{o}C)$ approximately equal 
to $2d$; here $\partial_{i}$ and $\partial_{o}$ correspond to the inner 
and outer boundaries of $A$, where $\partial_{o}A = S_{y}(1+d).$ In 
particular, there is a 2-sphere $S^{2}$ embedded in $C \subset  {\mathbb 
R}^{3},$ and essential in $C$ in the sense that $S^{2}$ separates 
$\partial_{i}C$ from $\partial_{o}C.$ 

  We now analyse in detail the structure of the metric $g' $ in 
$A_{\kappa d}.$ To do this, we essentially linearize the 
$\mathcal{Z}_{c}^{2}$ equations at the ``limit'' metric $g_{o}$ and ``limit'' 
potential $u \equiv  u_{o},$ i.e. we analyse the first order deviation 
of $(g, u)$ from the flat pair $(g_{o}, u_{o})$ in $A_{\kappa d}$. For 
simplicity, assume that 
\begin{equation} \label{e4.10}
u_{1} \equiv \sup_{A_{\kappa d}}u =  1;
\end{equation}
this may be achieved, w.l.o.g, by renormalizing the defining equations 
(1.16)-(1.17), i.e. changing $\alpha $ to $\alpha /u_{1}.$ Of course 
$u_{1} \sim u_{o}$ in (4.6).

 Consider the $\mathcal{Z}_{c}^{2}$ equations (1.16)-(1.17) in this scale, 
i.e. (using (4.1)),
\begin{equation} \label{e4.11}
\alpha\nabla\mathcal{Z}^{2} = L^{*}u, 
\end{equation}
\begin{equation} \label{e4.12}
\Delta u = \frac{\alpha}{4}|z|^{2}. 
\end{equation}
The coefficient $\alpha $ scales as the square of the distance, (so 
that the equations (4.11)-(4.12) are scale-invariant, c.f. [5, \S 
4.1]). Thus, in the scale above, $\alpha  \sim  R^{-2} << $ 1 while 
$|\nabla\mathcal{Z}^{2}| < O(\delta_{1}^{2}) << $ 1. This latter estimate 
follows from uniform elliptic regularity estimates for the 
$\mathcal {Z}_{c}^{2}$ equations in regions where $u$ satisfies 
(4.6), c.f. [5, Thm. 4.2] together with (4.9). Let
\begin{equation} \label{e4.13}
\mu  = \bigl(\int_{A_{d/2}}|r|^{2}dV\bigr)^{1/2}, 
\end{equation}
so that $\mu $ is small, depending on $\delta_{1}.$ We formally 
linearize the equations (4.11)-(4.12) at the flat metric $g_{o},$ by 
dividing by $\mu ,$ to obtain
$$\frac{\alpha}{\mu}\nabla\mathcal{Z}^{2} = L^{*}(\frac{u}{\mu}), $$

$$\Delta (\frac{1- u}{\mu}) = \frac{\alpha}{4\mu}|z|^{2}. $$
By the remarks on regularity above, $\delta_{1} \leq c \cdot \mu$, for 
a constant $c$ indpendent of $\delta_{o}$. Hence 
$|\frac{\alpha}{\mu}\nabla\mathcal{Z}^{2}| << $ 1 and 
$\frac{\alpha}{\mu}|z|^{2} << 1$ for $\delta_{o}$ sufficiently small. 
This gives 
\begin{equation} \label{e4.14}
|L^{*}(\frac{u}{\mu})| <<  1 \ \ {\rm and} \ \ |\Delta (\frac{1- 
u}{\mu})| <<  1, 
\end{equation}
away from $\partial A$. Thus, the function
\begin{equation} \label{e4.15}
\psi  \equiv  \frac{1- u}{\mu} 
\end{equation}
is almost harmonic w.r.t. the $g' $ metric, and hence almost harmonic 
w.r.t. the $g_{o}$ metric. Further, expanding $L^{*}$ in (4.14) gives
$$|D^{2}\psi  -  u(\frac{r}{\mu})| <<  1, $$
and since $u$ is close to 1,
\begin{equation} \label{e4.16}
\frac{r}{\mu} = D^{2}\psi  + o(1), 
\end{equation}
away from $\partial A$, where $o(1)$ is (arbitrarily) small if 
$\delta_{o}$ is sufficiently small. Note that by construction, i.e. by 
(4.13) and the elliptic regularity, $\frac{r}{\mu}$ is uniformly 
bounded, independent of the smallness of $\delta_{o},$ and hence so is 
$D^{2}\psi ,$ i.e. there exists $K$, independent of $\delta_{o}$, such 
that in $A_{\kappa d/2}$,
\begin{equation} \label{e4.17}
|D^{2}\psi| \leq  K. 
\end{equation}

 We claim that $|\psi|$ is also uniformly bounded above, i.e. for $x 
\in A_{\kappa d}$,
\begin{equation} \label{e4.18}
|\psi|(x)  \leq  L, 
\end{equation}
for some constant $L$, independent of $\delta_{o}$ small. To see this, 
let
\begin{equation} \label{e4.19}
\delta  = osc_{A_{\kappa d}}u \leq \delta_{o}. 
\end{equation}
Then for $\lambda  = \frac{\delta}{\mu},$ the function $\bar \psi = 
\psi /\lambda  = (1- u)/\delta $ has oscillation equal to 1 on 
$A_{\kappa d}.$ Now if $\lambda  >> $ 1, it follows from (4.17) that
$$|D^{2}\bar \psi| <<  1, $$
so that $\bar \psi = (1- u)/\delta $ is almost an affine function on 
$A_{\kappa d}.$ For $\lambda$ sufficiently large, this however 
contradicts the fact that the level set $L_{o}$ in (4.7) is compact and 
contained in $A_{\kappa d}$. Hence, $\lambda$ is bounded above for all 
$\delta_{o}$ sufficiently small and so, via (4.10), (4.18) holds.

  By (4.13) and elliptic regularity again, the $L^{\infty}$ norm of 
$\frac{r}{\mu}$ is on the order of 1, and hence by (4.16), $\psi $ 
cannot be too close to 0, or to any affine function on $A_{\kappa d}$. 
In particular since (4.19) gives $osc_{A_{\kappa d}} \psi = 
\frac{\delta}{\mu}$, $\lambda$ is also bounded below away from 0. Thus 
the ratio $\mu /\delta $ is bounded away from 0 and $\infty ,$ 
independent of $\delta_{o}$ small, i.e.
\begin{equation} \label{e4.20}
\mu  \sim  \delta . 
\end{equation}
Hence, from (4.15) one may write
\begin{equation} \label{e4.21}
u = 1 -  \delta\nu , 
\end{equation}
where $\nu $ differs from $\psi $ by a bounded scale factor; in 
particular $\nu $ is almost a positive harmonic function on $(A_{\kappa 
d}, g' ),$ uniformly bounded away from 0 and $\infty ,$ independent of 
$\delta_{o}.$

 Similarly, since $\frac{r}{\delta} \sim $ 1, i.e. formally the 
linearization of the curvature $r$ is bounded, we may write
\begin{equation} \label{e4.22}
g = g_{o} + \delta h + o(\delta ), 
\end{equation}
where $|h| \sim $ 1 in the $C^{o}$ topology, and hence $|h| \sim $ 1 in 
the $C^{k}$ topology, by elliptic regularity for the $\mathcal{Z}_{c}^{2}$ 
equations as above.

 Now as in [5, \S 7] for instance, consider the conformally 
equivalent metric
\begin{equation} \label{e4.23}
\widetilde g = u^{2}\cdot  g. 
\end{equation}
A standard computation of the Ricci curvature under conformal changes, 
c.f. [6, Ch.1J], gives
$$\widetilde r = r -  u^{-1}D^{2}u + 2(d\log u)^{2} -  u^{-1}\Delta 
u\cdot  g = $$
$$= 2(d\log u)^{2} -  \frac{\alpha}{2u}|z|^{2}\cdot  g -  
\frac{\alpha}{u}\nabla\mathcal{Z}^{2}. $$
where the second equality follows from (4.11)-(4.12). Since 
$|\nabla\mathcal{Z}^{2}| = O(\delta^{2}), \alpha  = O(R^{-2})$ and, from 
(4.19), $|d\log u| = O(\delta ),$ it follows that
\begin{equation} \label{e4.24}
\widetilde r \sim  2(du)^{2} = O(\delta^{2}), 
\end{equation}
and hence $\widetilde g$ is flat to order $\delta^{2}$, i.e.
$$\widetilde g = g_{o} + O(\delta^{2}).$$ 
Since $u =$ 1 $-  \delta\nu ,$ we obtain the expansion
\begin{equation} \label{e4.25}
g'  = (1+2\nu\delta )g_{o} + o(\delta ) 
\end{equation}
in $A_{\kappa d}$. This improves the estimate (4.22), i.e. shows that 
$$h = 2\nu g_{o},$$ 
so that to first order in $\delta ,$ the metric $g_{o}$ differs only 
conformally from the flat metric $g_{o}.$ Note in particular that this 
form of the metric agrees with the form (0.17) of $g$ in an 
asymptotically flat end, with $\nu $ then corresponding to $m/r$ - a 
multiple of the Green's function on ${\mathbb R}^{3}.$

\medskip

 Having identified the form of $g' $ and $u$ to first order in $\delta 
,$ we now are in position to carry out the metric sphere surgery of 
Theorem 3.1 under these circumstances. Thus, as before, we argue by 
contradiction and assume that the 2-sphere $S^{2},$ essential in $A$, 
bounds a 3-ball in $M$. 

 Note first that we need only consider the situation where $S^{2} 
\subset  A$ does not bound an end $E \subset  N$ on which osc $u \leq  
\delta_{o}.$ For in this situation, the assumptions of Proposition 4.1 
are satisfied, and so $(N, g')$ has an asymptotically flat end; the 
proof of Theorem 4.2 then proceeds as in Theorem 3.1. Hence for the 
remainder of the proof, assume that $S^{2}$ does not bound such an 
end in $N$, but does bound in $M$.

 Thus, suppose the $S^{2}$ in $A$ bounds on the inside in $M$, i.e. the 
component containing the base point $y$ and bounding $S^{2}$ is a 
3-ball in $M$. This implies that (a smooth approximation to) the outer 
boundary $\partial_{o}A_{\kappa d}$ of $A_{\kappa d}$ bounds a compact 
domain $W$ in ${\mathbb R}^{3},$ with $W \subset M$. Hence the level set 
$L_{o}$ of $\nu $ from (4.7) bounds a compact subdomain, (not 
necessarily a ball), in $W$.

 The metric $g' $ is real-analytic in $A$, as is the potential function 
$u$, c.f. Theorem 1.3. Hence the level sets of $u$, and so $\nu ,$ are 
real-analytic. Assuming the level set $L_{o}$ in (4.7) is regular, i.e. 
there are no critical points of $\nu $ on $L_{o},$ view the level set 
$L_{o}$ as isometrically embedded in $({\mathbb R}^{3}, g_{o})$. Since 
$\nu$ is constant on $L_{o}$, the Riemannian surface $(L_{o}, g')$ 
embeds to first order in $\delta$ isometrically in ${\mathbb R}^{3}$; in 
fact $(L_{o}, g')$ to first order is just the dilation of $(L_{o}, 
g_{o})$ by the factor $(1+2\nu \delta)$ by (4.25). It then follows from 
a result of [20], that $(L_{o}, g')$ itself embeds isometrically in 
${\mathbb R}^{3}$, for $\delta$ sufficiently small. Similarly, all the 
other regular level sets $L$ of $\nu $ which are compactly contained in 
$A_{\kappa d}$ embed isometrically in $(W, g_{o}).$ Such level sets $L$ 
then bound a smooth compact domain $V = V(L, \delta ) \subset  (W, 
g_{o}),$ which as above vary with $\delta .$ 

 For a suitable choice of the level $L$, (near $L_{o}$), we use the 
metric 
\begin{equation} \label{e4.26}
\bar g = g_{o}\cup g'  , g_{o} = g_{o}|_{V}, 
\end{equation}
as a comparison metric to $g' ,$ as in Case I of the proof of Theorem 
3.1. This metric is piecewise smooth, but is only $C^{o}$ at the seam 
$L$.

 The two main ingredients in the comparison argument in Case I of 
Theorem 3.1 are the relations between the $2^{\rm nd}$ fundamental forms 
(3.8) and the volumes (3.11) of the flat metric $g_{o}$ and the metric 
$g' .$

 First, we prove an analogue of the estimate (3.8). From the form of 
the metric $g' $ in (4.25), an easy computation for conformal metrics, 
(c.f. [6, Ch.1J]), gives
\begin{equation} \label{e4.27}
A_{g_{o}} -  A_{g'} = \ - <\nabla\nu , X>\delta \cdot I + o(\delta ), 
\end{equation}
where $X$ is the unit outward normal at $L$ and $I = g'|_{L}.$ By 
construction, i.e. from the hypothesis following (4.7), we have $- 
<\nabla\nu , X>  = |\nabla\nu| > $ 0 on any regular level $L$. Further, 
since $\nu  \sim $ 1 and osc $\nu  \sim $ 1 in $A_{\kappa d},$ 
independent of $\delta_{o},$ and so $\nu $ is uniformly controlled 
independent of $\delta_{o},$ it follows that there exist regular levels 
$L$ near $L_{o}$ such that
\begin{equation} \label{e4.28}
 |\nabla\nu| _{L} \geq  \mu_{o} >  0, 
\end{equation}
for some positive constant $\mu_{o},$ independent of $\delta_{o}.$

 Thus, as in (3.8), $L$ is more convex in the flat metric $g_{o}$ than 
in the $g' $ metric, and so there is a concentration of positive 
curvature on $L$, provided $\delta_{o}$ is sufficiently small. As in 
(3.13), one may then smooth the metric $\bar g$ near $L$ to a metric 
$\widetilde g$ satisfying
\begin{equation} \label{e4.29}
(s_{\widetilde g}) ^{-} \geq  (s_{g'}) ^{-} = 0, 
\end{equation}
everywhere.

 Observe that the blow-up metrics $g_{\varepsilon}' $ limiting on $g' $ 
have the same form (4.25) on $A_{\kappa d}$ provided $\varepsilon $ is 
chosen sufficiently small, depending only on $\delta_{o}.$ Hence, one 
may carry out the construction above, (with the same $g_{o}$), to 
obtain comparison metrics $\bar g_{\varepsilon}$ to $g_{\varepsilon}' 
,$ for which there are smoothings $\widetilde g_{\varepsilon}$ 
satisfying
\begin{equation} \label{e4.30}
(s_{\widetilde g_{\varepsilon}}) ^{-} \geq  (s_{g_{\varepsilon}'}) 
^{-}, 
\end{equation}
for $\varepsilon $ sufficiently small.

 To estimate the difference in the volumes, suppose first that $S^{2}$ 
does not bound in $N$; of course it does bound a 3-ball $B^{3}$ in $M$. 
In this situation,
$$vol_{g_{\varepsilon}'} B^{3} \rightarrow  \infty ,  $$
and so obviously, since the volume of the limit comparison flat 3-ball 
is finite,
\begin{equation} \label{e4.31}
vol_{\widetilde g_{\varepsilon}} B^{3} <  vol_{g_{\varepsilon}'} B^{3}, 
\end{equation}
for $\varepsilon $ sufficiently small.

 Hence we may assume that $S^{2}$ bounds a compact 3-ball (on the 
inside) in $N$. It follows from the discussion above that the domain in 
$N$ bounded by the level set $L$ is diffeomorphic to the flat domain 
$V$ above. Thus we have two metrics, the metric $g' ,$ (rescaled by 
$R^{-2}),$ and the flat metric $g_{o},$ on $V$.

 We first claim that the expansion (4.25) is valid outside a subset of 
(arbitrarily) small $g_{o}$ volume in $V$, depending only on 
$\delta_{o}.$ To see this, the expression (4.25) is valid in the 
subdomain $V_{k}$ of $V$ on which $\nu $ satisfies
$$\nu  \leq  k, $$
where $k$ may be made arbitrarily large if $\delta_{o}$ is chosen 
sufficiently small. Now the function $\nu$ is approximately a harmonic 
function on $(V_{k}, g'),$ with boundary value $\nu  \equiv $ const. on 
$\partial V = L$; in fact from (4.14), $|\Delta_{g'} \nu| \leq cR^{-2} 
\cdot \delta$. From the expression (4.25), a short computation then 
shows that
$$|\Delta_{g_{o}}(e^{\nu} - 1)| = R^{-2}O(\delta).$$
Let $v = e^{\nu} - 1$ and $v_{k} = \min(v, e^{k})$, so that 
(essentially) $v_{k} \equiv e^{k}$ on $V \setminus V_{k}$. It follows 
that $v_{k}$ is close to a positive superharmonic function $\phi_{k}$ 
on $(V, g_{o})$. It is well known, c.f. [17, Thm.5.8] that the measure 
of the set where such a function is very large is small, and hence
\begin{equation} \label{e4.32}
vol_{g_{o}}\{\nu  \geq  k\} \leq  \delta_{2}, 
\end{equation}
where $\delta_{2} = \delta_{2}(\delta , k)$. This proves the claim 
above.

 Given these preliminaries, we now do the volume comparison. We have 
$vol_{g'}V \geq  vol_{g'}V_{k}$. From the expression (4.25) and from 
the construction of $V$, c.f. the discussion preceding (4.26),
\begin{equation} \label{e4.33}
vol_{g'}V_{k} = vol_{g_{o}}V_{k}^{o} + 3(\int_{V_{k}^{o}}\nu 
dV_{g_{o}})\delta  + o(\delta ). 
\end{equation}
Observe here that $V_{k}^{o} \neq V_{k}$, but by construction, 
$V_{k}^{o}$ is the domain in $({\mathbb R}^{3}, g_{o})$ with $\partial 
V_{k}^{o} = L$, where $L$, (and {\it not} $(L, g_{o})$) is 
isometrically embedded in $({\mathbb R}^{3}, g_{o})$; this relation is 
just the same as the relation between $R$ and $\bar R$ in Case I of 
Theorem 3.1. On the other hand, since to $1^{st}$ order in $\delta$, 
$(L, g' )$ is just the dilation of $(L, g_{o})$ by the factor 
$(1+2\nu|_{L}\delta ),$ one has
\begin{equation} \label{e4.34}
vol_{g_{o}}V_{k} = (1+3\nu|_{L} \delta )\cdot  vol_{g_{o}}V_{k}^{o} + 
o(\delta ). 
\end{equation}
Thus, modulo lower order terms in $\delta$,
\begin{equation} \label{e4.35}
vol_{g_{o}}V_{k}^{o} \sim (1-3\nu|_{L} \delta )\cdot  vol_{g_{o}}V_{k}.
\end{equation}
Since, again by the hypothesis following (4.7) and the maximum 
principle,
\begin{equation} \label{e4.36}
\frac{1}{vol_{g_{o}} V_{k}}\int_{V_{k}}\nu dV_{g_{o}} >  \nu|_{L}, 
\end{equation}
independent of the size of $\delta ,$ it follows from (4.33)-(4.36) that
\begin{equation} \label{e4.37}
vol_{g'}V >  vol_{g_{o}}V_{k}, 
\end{equation}
for $\delta$ sufficiently small. This estimate together with (4.32), 
gives
\begin{equation} \label{e4.38}
vol_{g'}V >  vol_{g_{o}}V, 
\end{equation}
for $\delta$ small. The estimate (4.38) will of course also hold for 
$g_{\varepsilon}' $ in place of $g' ,$ and $\widetilde g_{\varepsilon}$ 
in place of $g_{o},$ for $\varepsilon $ sufficiently small, and so we 
have the analogue of (3.11), i.e.
\begin{equation} \label{e4.39}
vol_{g_{\varepsilon}'}V >  vol_{\widetilde g_{\varepsilon}}V. 
\end{equation}

 The estimate for the comparison for $\mathcal{Z}^{2}$ is essentially the 
same as in Case I of Theorem 3.1. Thus the metric $\widetilde g$ has $z 
\equiv $ 0 inside the seam $L$, while in a small band $T$ about $L$ 
where $\bar g$ is smoothed,
\begin{equation} \label{e4.40}
\int_{T}|\widetilde z|^{2}dV_{\widetilde g} \leq  \delta_{3}, 
\end{equation}
where $\delta_{3}$ may be made arbitrarily small by choosing 
$\delta_{o}$ sufficiently small. On the other hand, the limit $(N, g' , 
y)$ has a definite amount of curvature inside $S^{2},$ exactly as in 
(3.18), since $R$ is large. Thus, for the same reasons as in 
(3.19)-(3.20), one obtains in this situation
\begin{equation} \label{e4.41}
I_{\varepsilon}^{~-}(\widetilde g_{\varepsilon}) <  
I_{\varepsilon}^{~-}(g_{\varepsilon}' ), 
\end{equation}
which is impossible by the minimizing property of $g_{\varepsilon}'$. 
This completes the proof.

{\endproof}

\begin{remark} \label{r 4.3.}
  In the context of Theorem 4.2, suppose the 2-sphere $S^{2}$ in $A$ 
bounds a 3-ball on the outside in $M$. As remarked in the proof above, 
if $S^{2}$ bounds in $U^{o} \subset  N$, then again Proposition 4.1 
proves that there is an asymptotically flat sub-end, and one may apply 
Theorem 3.1.

 Thus suppose $S^{2}$ bounds a 3-ball $B$ on the outside, but does not 
bound in $U^{o}.$ In this case, if there is another larger annulus $A' 
,$ i.e. a component of a larger geodesic annulus $A((1- d)R' , (1+d)R' 
)$ with $R'  >  2R$, satisfying the assumptions of Theorem 4.2, then 
$A'  \subset  M$ is topologically contained in the 3-ball $B$. Hence, 
the 2-sphere $(S^{2})'$ in $A'$ bounds a 3-ball on the inside, and one 
obtains a contradiction again from the proof of Theorem 4.2. Hence, 
under such circumstances, the (original) $S^{2}$ cannot bound a 3-ball 
in $M$ on either side.

 As an example of such a situation, one might (possibly) have blow-up 
limit $\mathcal{Z}_{c}^{2}$ solutions $(N, g')$ which for instance are 
topologically the double of ${\mathbb R}^{3}\setminus \cup B_{i},$ where 
$B_{i}$ is a countable collection of disjoint 3-balls in ${\mathbb 
R}^{3},$ i.e. $N$ is an infinite connected sum of ${\mathbb R}^{3}$'s. 
Such a manifold cannot be asymptotically flat, but might satisfy the 
preceding condition.

\end{remark}

\bibliographystyle{plain}

\begin{thebibliography}{WWW}

\footnotesize


\bibitem [1]{1} {\it M. Anderson},  Scalar curvature and geometrization 
conjectures for 3-manifolds, Comparison Geometry, MSRI Publ., {\bf 30}, 
(1997), 49-82.

\bibitem [2]{2} {\it M. Anderson},  Extrema of curvature functionals on the 
space of metrics on 3-manifolds, Calc. Var. \& P.D.E., {\bf 5}, (1997), 
199-269.

\bibitem [3]{3} {\it M. Anderson},  Extrema of curvature functionals on the 
space of metrics on 3-manifolds, II, Calc.Var. \& P.D.E., {\bf 12}, 
(2001), 1-58.

\bibitem [4]{4} {\it M. Anderson},  Scalar curvature, metric degenerations 
and the static vacuum Einstein equations on 3-manifolds, I, Geom. \& Funct. 
Anal., {\bf 9}, (1999), 855-967.

\bibitem [5]{5} {\it M. Anderson},  Scalar curvature, metric degenerations 
and the static vacuum Einstein equations on 3-manifolds, II, Geom. \& Funct. 
Anal., {\bf 11}, (2001), 273-381.

  The papers above are also available at  http://www.math.sunysb.edu/ 
$\sim$ anderson/

\bibitem [6]{6} {\it A. Besse}, Einstein Manifolds, Ergebnisse der Mathematik, 
3. Folge, Band 10, Springer Verlag, New York (1987).

\bibitem [7]{7} {\it J. Cheeger and M. Gromov}, Collapsing Riemannian 
manifolds while keeping their curvature bounded I, Jour. Diff. Geom., 
{\bf 23}, (1986), 309-346.

\bibitem[8]{8} {\it J. Cheeger and M. Gromov},  Collapsing Riemannian 
manifolds while keeping their curvature bounded II, Jour. Diff. Geom., 
{\bf 32}, (1990), 269-298.

\bibitem[9]{9} {\it K. Fukaya}, Collapsing Riemannian manifolds to ones of 
lower dimension, Jour. Diff. Geom., {\bf 25}, (1987), 139-156.

\bibitem[10]{10} {\it S. Gallot}, Isoperimetric inequalities based on integral 
norms of Ricci curvature, Asterisque {\bf 157-158}, (1988), 191-217.

\bibitem[11]{11} {\it D. Gilbarg and N. Trudinger}, Elliptic Partial 
Differential Equations of the Second Order, 2nd Edition, Springer 
Verlag, New York, 1983.

\bibitem[12]{12} {\it M. Gromov}, Metric Structures for Riemannian and 
Non-Riemannian Spaces, Progress in Math., {\bf 152}, Birkhauser Verlag, 
Boston, (1999).

\bibitem[13]{G13} {\it M. Gromov}, Volume and bounded cohomology, Publ. Math. 
I.H.E.S., {\bf 56}, (1982), 5-100.

\bibitem[14]{G14} {\it M. Gromov and H. Lawson, Jr.}, Spin and scalar 
curvature in the presence of a fundamental group, Annals of Math., 
{\bf 111}, (1980), 209-230.

\bibitem[15]{15} {\it M. Gromov and H. Lawson, Jr.}, Positive scalar 
curvature and the Dirac operator on complete Riemannian manifolds, 
Publ. Math. I.H.E.S., {\bf 58}, (1983), 83-196.

\bibitem[16]{16} {\it R. Hamilton}, Non-singular solutions of the Ricci flow 
on three-manifolds, Comm. Geom. Anal, {\bf 7}, (1999), 695-729.

\bibitem[17]{17} {\it W. Hayman}, Subharmonic Functions, Academic Press, 
London, (1976).

\bibitem[18]{18} {\it W. Jaco}, Lectures on Three-Manifold Topology, Conf. Bd. 
Math. Sci., A.M.S., {\bf 43}, (1980).

\bibitem[19]{19} {\it W. Jaco and P. Shalen}, Seifert fibered spaces in 
3-manifolds, Memoirs A.M.S., {\bf 220}, (1979).

\bibitem[20]{20} {\it H. Jacobowitz}, Local isometric embeddings of surfaces 
into Euclidean four space, Indiana Univ. Math. Jour., {\bf 21}, (1971), 
249-254.

\bibitem[21]{21} {\it K. Johannson}, Homotopy equivalence of 3-manifolds with 
boundaries, Lect. Notes in Math., {\bf 761}, Springer Verlag, (1979).

\bibitem[22]{22} {\it M. Kneser}, Geschlossene Fl\"achen in dreidimensionale 
Mannifaltigkeiten, Jahresber. Deutsch. Math. Verein, {\bf 38}, (1929), 248-260.

\bibitem[23]{23} {\it O. Kobayashi}, Scalar curvature of a metric of unit 
volume, Math. Annalen, {\bf 279}, (1987), 253-265.

\bibitem[24]{24} {\it J. Milnor}, A unique factorization theorem for 
3-manifolds, Amer. J. Math., {\bf 84}, (1962), 1-7.

\bibitem[25]{25} {\it G. Mostow}, Quasiconformal mappings in n-space and the 
strong rigidity of hyperbolic space forms, Publ. Math. I.H.E.S., {\bf 34}, 
(1968), 53-104.

\bibitem[26]{26} {\it P. Orlik}, Seifert Manifolds, Lect. Notes in Math., 
vol. 291, Springer Verlag, New York, (1972).

\bibitem[27]{27} {\it G. Prasad}, Strong rigidity of Q-rank 1 lattices, 
Inventiones Math., {\bf 21}, (1973), 255-286.

\bibitem[28]{28} {\it X. Rong}, The limiting eta invariant of collapsed 
3-manifolds, Jour. Diff. Geom., {\bf 37}, (1993), 535-568.

\bibitem[29]{29} {\it R. Schoen},  Conformal deformation of a metric to 
constant scalar curvature, Jour. Diff. Geom., {\bf 20}, (1985), 479-495.

\bibitem[30]{30} {\it R. Schoen},  Variational theory for the total scalar 
curvature functional for Riemannian metrics and related topics, Lect. Notes 
Math. {\bf 1365}, Springer Verlag, Berlin, (1987), 120-154.

\bibitem[31]{31} {\it R. Schoen and S.-T. Yau}, Existence of incompressible 
minimal surfaces and the topology of three dimensional manifolds 
with non negative scalar curvature, Annals of Math., {\bf 110}, (1979), 
127-142.

\bibitem[32]{32} {\it R. Schoen and S.-T. Yau}, On the structure of 
manifolds with positive scalar curvature, Manuscripta Math., {\bf 29}, 
(1979), 159-183.

\bibitem[33]{33} {\it P. Scott}, The geometries of 3-manifolds, Bull. London 
Math. Soc., {\bf 15}, (1983), 401-487.

\bibitem[34]{34} {\it T. Soma}, The Gromov invariant of links, Inventiones 
Math., {\bf 64}, (1981), 445-454.

\bibitem[35]{35} {\it M. Spivak}, A Comprehensive Introduction to Differential 
Geometry, V, Publish or Perish, Houston, 1979.

\bibitem[36]{36} {\it W. Thurston}, The Geometry and Topology of 3-Manifolds, 
(preprint), Princeton, 1979.

\bibitem[37]{37} {\it W. Thurston}, Three dimensional manifolds, Kleinian 
groups and hyperbolic geometry, Bulletin A.M.S., {\bf 6}, (1982), 357-381.

\bibitem[38]{38} {\it F. Waldhausen}, Eine Klasse von 3-dimensionalen 
Mannigfaltigkeiten I, Inventiones Math., {\bf 3}, (1967), 308-333.

\bibitem[39]{39} {\it F. Waldhausen}, Eine Klasse von 3-dimensionalen 
Mannigfaltigkeiten II, Inventiones Math., {\bf 4}, (1967), 87-117.

\end{thebibliography}

\begin{center}
September, 2000
\end{center}
\medskip
\address{Department of Mathematics\\
S.U.N.Y. at Stony Brook\\
Stony Brook, N.Y. 11794-3651}\\
\email{anderson@math.sunysb.edu}

\end{document}